\documentclass[twoside]{elsart}
\usepackage{amssymb,amsmath,graphicx}
\usepackage{amsfonts}
\newtheorem{corollary}{Corollary}
\newtheorem{definition}{Definition}
\newtheorem{lemma}{Lemma}
\newtheorem{proposition}{Proposition}
\newtheorem{remark}{Remark}
\newtheorem{theorem}{Theorem}
\newtheorem{problem}{Open Problem}
\newcounter{AD}

\newcommand{\pp}{\noindent {\em Proof. }}
\newcommand{\bx}{\hfill$\Box$}
\newcommand{\bee}[1]{\begin{equation}\label{#1}}
\newcommand{\ene}{\end{equation}}
\newcommand{\beq}[1]{\begin{eqnarray}\label{#1}}
\newcommand{\eqe}{\end{eqnarray}}
\newcommand{\bea}{\begin{eqnarray*}}
\newcommand{\eqa}{\end{eqnarray*}}
\newcommand{\cd}{\cdots}
\newcommand{\cir}{\,{\scriptstyle\circ}\,}
\newcommand{\ld}{\ldots}
\newcommand{\ra}{\rightarrow}
\newcommand{\vp}{\varphi}
\newcommand{\ve}{\varepsilon}
\newcommand{\G}{\Gamma}
\newcommand{\cg}{\mathcal{G}}
\newcommand{\cc}{\mathcal{C}}
\newcommand{\ct}{\mathcal{T}}
\newcommand{\card}{\#\,}

\newcommand{\df}{\mathrm{def}}
\newcommand{\Def}[2]{\mathrm{def}_{#1}(#2)} 
\newcommand{\dist}{\mathrm{dist}}

\newcommand{\lb}[1]{\mathrm{Lab}(#1)}

\newcommand{\Lab}[1]{\mathrm{Lab}(#1)}

\newcommand{\rk}[1]{\mathrm{rank}\,#1}

\newcommand{\bA}{\overline{A}}

\newcommand{\bU}{\overline{U}}
\newcommand{\bV}{\overline{V}}

\newcommand{\fa}{free associative algebra }

\newcommand{\fas}{free associative algebras }

\newcommand{\fac}{free associative algebra, }

\newcommand{\far}[1]{\mathcal{A}\langle x_1,\ld,x_{#1}\rangle}
\newcommand{\wax}{W(X)}
\newcommand{\fgax}{F(X)}
\newcommand{\fax}{\mathcal{A}\langle X\rangle}
\newcommand{\gax}{\mathcal{F}\langle X\rangle}
\renewcommand{\wr}{W_r}
\newcommand{\fr}{F_r}
\newcommand{\ar}{\mathcal{A}_r}
\newcommand{\fgr}{\mathcal{F}_r}
\newcommand{\mg}{maximal growth }

\newcommand{\mgc}{maximal growth, }
\newcommand{\mgp}{maximal growth. }
\newcommand{\nat}{\mathbb{N}}
\newcommand{\mb}[2]{\mathcal{B}(#1,#2)}
\newcommand{\ms}[2]{\mathcal{S}(#1,#2)}
\newcommand{\mcb}{\mathcal{B}}
\begin{document}
\begin{frontmatter}

\title{Actions of Maximal Growth}

\author[yb,ybao,yb1]{Yuri Bahturin} and \author[ao,ybao,ao1]{Alexander Olshanskii}
\address[yb]{Department of Mathematics and Statistics, Memorial University of Newfoundland, St. John's, NL, A1C 5S7, \textsc{Canada}}
\address[ao]{Department of Mathematics,
1326 Stevenson Center,
Vanderbilt University, Nashville, TN 37240, \textsc{USA}}
\address[ybao]{Department of Algebra, Faculty of Mathematics and Mechanics, 119899 Moscow, \textsc{Russia}}
\thanks[yb1]{Partially supported by NSERC grant \# 227060-04}
\thanks[ao1]{Partially supported by NSF grants DMS-0700811 and DMS-0455881 and by RFBR grant  08-01-00573}

\begin{abstract}
We study acts and modules of maximal growth over finitely generated free monoids and free associative algebras as well as free groups and free group algebras.  The maximality of the growth implies some other specific properties of these acts and modules that makes them close to the free ones; at the same time, we show that being a strong ``infiniteness'' condition, the maximality of the growth can still be combined with various finiteness conditions, which would normally make finitely generated acts finite and finitely generated modules finite-dimensional.
\end{abstract}

\begin{keyword}
Free algebra\sep module\sep free group\sep free monoid\sep $G$-set\sep act\sep growth

\MSC 17B01\sep 17B50\sep20F40 
\end{keyword}
\end{frontmatter}

\tableofcontents

\section*{Introduction}\addcontentsline{toc}{section}{Introduction}

Let us consider the following very general mathematical context. Suppose that a finite collection of linear operators $\mathcal{L}_1,\ld,\mathcal{L}_r$ acts on an (infinite-dimensional) linear space $V$ over a field $\Phi$. We assume that there is a finite-dimensional subspace $V(0)\subset V$ such that $V$ is generated by $V(0)$ with the help of the above operators. In other words, $V = \bigcup V(i)$, where, for $i>0$ we set 
\bee{e000}
V(i) = V(i-1)+ \mathcal{L}_1(V(i-1))+\cd+\mathcal{L}_r(V(i-1)).
\ene
In this case $g(n) = \dim V(n)$ is the growth function of the space $V$ with respect to the fixed set of linear operators and fixed generating subspace $V(0)$.

If we replace $V(0)$ by any other finite-dimensional subspace $V'(0)$ then there is some $c\ge 0$ such that $V'(0) \le  V(c)$, and so $g'(n) \le g(n+c)$ where $g'(n)$ is the growth function of $V$ defined by the subspace $V'(0)$. If we define the equivalence of two functions $g$ and $g'$ by the conditions  $g'(n) \le g(n+c)$ and $g(n) \le g'(n+c')$, for some positive constants $c$ and $c'$, then we will obtain an invariant of the space $V$ with the action of the given set of linear operators, which does not depend on the choice of a finite-dimensional generating subspace. It is natural to call the respective equivalence class the growth of $V$ with respect to $\mathcal{L}_1,\ld,\mathcal{L}_r$. (As an example, the functions $n$ and $n+100$ are equivalent while $2n$ and $4n$ are not; also $2^n$ and $3^n$ are not equivalent). 

Now the action of $\mathcal{L}_1,\ld,\mathcal{L}_r$ uniquely defines the action of all linear combinations of the composition of there operators. Thanks to the universal property of the free associative algebra $\mathcal{A}_r=\far{r}$, if we map $x_1\mapsto \mathcal{L}_1,\ld,x_r\mapsto \mathcal{L}_r$ then $V$ becomes a finitely generated module over $\mathcal{A}_r$. 

If the linear operators $\mathcal{L}_1,\ld,\mathcal{L}_r$ are invertible then it is natural to assume that $V$ is generated by $V(0)$ with the use of $\mathcal{L}_1^{\pm 1},\ld,\mathcal{L}_r^{\pm 1}$. Then we will have to replace (\ref{e000}) by
\bee{e001}
V(i) = V(i-1)+ \mathcal{L}_1^{\pm 1}(V(i-1))+\cd+\mathcal{L}_r^{\pm 1}(V(i-1)).
\ene
In this case the same mapping $x_1\mapsto \mathcal{L}_1,\ld,x_r\mapsto \mathcal{L}_r$ makes $V$ into a module over the group algebra $\mathcal{F}_r=\Phi F_r$ of the free group $F_r=F(x_1,\ld,x_r)$ of rank $r$ over $\Phi$. Again, the growth of the action of invertible operators $\mathcal{L}_1,\ld,\mathcal{L}_r$ translates into the growth of modules over $\mathcal{F}_r$.

In the same way, it is natural to speak about the growth of an arbitrary set $S$ with a fixed set of transformations of any nature $f_1,\ld,f_r$; if $W_r=W(x_1,\ld,x_r)$ is the free monoid of rank $r$ then mapping $x_1\mapsto f_1,\ld,x_r\mapsto f_r$ makes $S$ into a $W_r$-act. If these transformations are bijective then, thanks to the universal property of the free group $F_r=F(x_1,\ld,x_r)$ of rank $r$, the set $S$ becomes an $F_r$-set. The similarly defined growth of $S$ becomes its invariant, and two $F_r$-sets with different growth are not isomorphic, in the same way as in the case of modules. 

In this paper we consider actions with maximal growth. This means $F_r$- and $W_r$-acts or $\mathcal{A}_r$- and $\mathcal{F}_r$-modules $(r>1)$, whose growth up to equivalence coincides with the growth which is maximal possible, that is, the growth of respective free acts and modules. One of the main examples of $F_r$-acts of maximal growth is the right action of the free group $F_r$ on the set $F_r/H$ of the right cosets of any finitely generated subgroup $H$ of infinite index. In the case of $\mathcal{A}_r$ and $\mathcal{F}_r$, important examples of modules of maximal growth are any infinite-dimensional modules with finite presentation, that is, given by finitely many generators and defining relations. According to \cite{PMC}, any such module has a free submodule of finite codimension, in particular, such module is large. Given an algebra $R$ over a field $\Phi$, a right $R$-module $M$ is called \emph{large} if it has a submodule of finite codimension which can be mapped onto the free $R$-module $R$. Free actions in all four cases and large module in the case of $\mathcal{A}_r$ and $\mathcal{F}_r$ have been studied in our previous paper \cite{BO2}.

We start with the discussion in Section \ref{sGA} of the definition and general properties of the growth of actions over arbitrary monoids and algebras. Then we prove some results which are common (or very close) for all our main examples, like what happens to the growth when we consider subactions or images of actions under the morphisms. We show that  the growth function $g(n)$ of a finitely generated action in the case of $W_r$ and $\mathcal{A}_r$ can be written as $g(n)=\alpha(n)r^n$, where $\alpha(n)$ is a function  converging to a limit $C_0\ge 0$ at infinity. The growth is maximal if and only if $C_0>0$. Similarly, in the case of $F_r$ and $\mathcal{F}_r$ the growth function $g(n)$ of a finitely generated action can be written as $g(n)=\alpha(n)(2r-1)^n$, with the same property for $\alpha(n)$. Note that there are some special properties in the individual cases, and they are discussed individually in the respective sections of the paper.

For example, it is proved in Corollary \ref{cRMG} that in every module $M$ over a free associative algebra or a free group algebra there is a unique maximal submodule $N$ none of whose submodules have maximal growth while in the quotient module $M/N$ every nonzero submodule has maximal growth.

This should be compared with so called bound and unbound modules in the sense of P. M. Cohn's \cite{PMC} where a module $M$ over an algebra $R$ is called \emph{bound} if $\mathrm{Hom}_R(M,R)=\{ 0\}$. Any module $N$ such that $\mathrm{Hom}_R(M,N)=\{ 0\}$ for any bound $M$ is called \emph{unbound}. The unbound modules can also be defined as modules without nonzero bound submodules. The class of all bound modules over a ring $R$ is closed under homomorphic images, extensions and direct limits (in particular, arbitrary sums). It follows that in every $R$-module $M$ there is a unique maximal bound submodule $M_b$ such that the quotient module $M/M_b$ is an unbound module. Thus one can view the class of bound modules as a radical class and that of unbound modules as a semisimple class. 

Our result mentioned above says that in the case where $R$ is a free associative algebra or a free group algebra, both of rank $r>1$, over a field the modules none of whose submodules have maximal growth form a radical class while those in which every submodule has maximal growth form a semisimple class. Since in this case a module is bound if and only if it does not contain $R$ as a direct summand, our radical class is a proper subclass of the class of bound modules while our semisimple class contains all unbound modules. An advantage of our radical class is that, in addition to the closure properties of bound modules, it is also closed under submodules. 

It is interesting that if we replace the maximality of the growth by some other condition, for example the popular exponentiality condition, then the just mentioned result is no more true and one can find examples of this kind in Subsection \ref{sENM}.  

The theorem we just mentioned makes actions where every nontrivial subaction has maximal growth looking like ``free'' actions while those without subactions of maximal growth like ``torsion'' actions. Another common property of actions of maximal growth and free actions is that both types of actions are faithful. However, when we begin examining the properties of actions of maximal growth, we quickly learn that the situation is much more complex than one would expect in the case of free actions. It turns out that actions of maximal growth may satisfy strong finiteness conditions which by their name are supposed to make finitely generated acts (respectively, modules) finite (respectively, finite-dimensional).

A ``popular'' family of finiteness conditions comes down from the Burnside Problems. Let us assume that $R$ is an algebra over a field $\Phi$ and there is a homomorphism $\ve: R\to \Phi$ (the ``augmentation map''). Suppose $\Delta$ is the kernel of $\ve$. We say that a right $R$-module $M$ is \emph{nil} if for each $a$ in a module $M$ and each $u\in\Delta$ there is a number $n=n(a,u)$ such that $au^n=0$. 

In the case of $G$-sets one may speak about the ``periodic'' action of $G$, meaning that all the orbits of the action of each $g\in G$ are finite. The study of such $G$-sets amounts to $G$-sets of the form $G/H$, $H$ a subgroup of the free group $G$, where for any $g\in G$ there is natural $n$ such that $g^n\in H$. One calls such subgroups \emph{Burnside}. If $H$ is normal and Burnside then the factor group $G/H$ is periodic.

It is very easy to produce straightforward examples of finitely generated infinite-dimensional nil-modules over free associative algebras or Burnside subgroups of infinite index in finitely generated free groups if we use widely known example of the negative solution of the Burnside problem in the case of associative algebras or groups. However, as observed in Section \ref{sFAMG}, the growth of the acts associated with these examples is never maximal.

So, using different approaches, in Theorem \ref{tNMMG} we show that given any graded module $M$ of maximal growth over a free associative algebra $R$ there is a submodule $N$ such that $M/N$ is an infinite-dimensional nil-module and still of maximal growth. In Proposition \ref{pNNMG} we show that this fails for any growth less than maximal. In addition to the modules of maximal growth with ``finiteness conditions'', we also describe a procedure which enables us to construct highly transitive acts and simple modules of maximal growth.

In the case of $F_r$-sets, we show in Section \ref{sGAMG}, among other results, that for any finitely generated subgroup $H$ of infinite index in a free nonabelian group $F_r$ there is a Burnside subgroup $K$ of infinite index in $F_r$, such that $H$ a free factor in $K$ and the growth of $F_r$-set $F_r/K$ is still maximal. In addition, for any $k>0$ the action of $F_r$ on $F_r/K$ is $k$-transitive so that $K$ is a maximal subgroup of $F_r$. In the associated cyclic module $M=\Phi(F_r/K)$ over the free group algebra $\mathcal{F}_r$, every element of $F_r$ acts as a locally periodic linear transformation. Additionally, $M$ has a simple submodule $N$ of codimension one. The growth of both $M$ and $N$ is maximal.

In the case of modules over free associative algebras, we produce a number of other examples of modules of maximal growth, as follows. A module $M$ over an algebra $R$ is called \emph{residually finite} if for any nonzero $a\in M$ there is $N\subset M$ such that $a\notin N$ and $\dim M/N <\infty$. Free modules over free algebras are easy examples of residually finite modules. But already the quotient modules of free modules need not be residually finite: take any simple finitely generated module of infinite dimension! (However, any submodule of a residually finite module is residually finite.)  What we manage to produce in Theorem \ref{t991}, is an example of a module $M$ of maximal growth such that any factor-module  $M/P$ is residually finite.  In the setting of Group Theory first example of this kind have been offered in \cite{OO}. 

In the study of $F_r$-sets $S$ of maximal growth our main tool is the Cayley graph of the action. In the case of $F_r$-sets, this graph is the graph of cosets of a subgroup $H$ of $F_r$, denoted by $\mathcal{G}(H)$. For the coset graphs of free groups, J. Stallings in \cite{JS1} introduced the notion of the \emph{core} $\mathcal{C}$ as the subgraph of $\mathcal{G}(H)$ consisting of the origin $o=H$ and all reduced loop starting at $o$. This core is finite if $H$ is finitely generated. We introduce a useful notion of the \emph{deficit} of the core, which is nonzero if and only if the growth of the $F_r$-set $F/H$ is maximal (again $H$ is finitely generated). Suppose we want to embed a subgroup $H$ of infinite index in a Burnside subgroup $K$, such that the growth of $F_r/K$ is still maximal. This can be done by adjoining to $H$, one by one, sufficiently great powers of all elements of $F$. On each step we have a greater subgroup, a new Cayley graph and its new core. We show that each consecutive power can chosen in such a way that the change of the deficit is arbitrarily small. Since the growth function is essentially determined by the deficit, the growth remains to be maximal even after we adjoin the powers of all elements and obtain the Burnside subgroup $K$. Additionally, one can carry out the construction in such a way that for any natural $k$ the action of $F_r$ on $F_r/K$ is $k$-transitive. In particular, $K$ is a maximal subgroup of $F_r$.

In Section \ref{ssTAMG} we associate with every (cyclic) action a closed subset in the measurable ultrametric space $\partial F_r$, respectively, $\partial W_r$, of infinite rays in the Cayley graph of $F_r$ (if we consider $F_r$-sets or $\mathcal{F}_r$-modules) or $W_r$ (if we consider $W_r$-acts or $\mathcal{A}_r$-modules), respectively. We prove that the growth is maximal if and only if the measure of the set is positive. 

Finally, let us emphasize once again that the structure of acts and modules arises on a set $S$ (or a linear space $V$) after we fix on $S$ (or $V$) several (linear) transformations $f_1,\ld,f_r$. This naturally selects a free basis $X=\{ x_1,...,x_r\}$ in the appropriate $\wr$, $\fr$, $\ar$ or $\fgr$ and a map $x_i\mapsto f_i$, $i=1,\ld,r$. The equivalence relation on the growth function introduced by us reflects the requirement that the growth is invariant under the isomorphisms of acts and modules. A coarser equivalences (hence, wider equivalence classes) arise when one considers so called ``semi-isomorphisms'' of acts and modules. For example, in the case of group actions,  we say that $F_r$-sets $S$ and $S^{\prime}$ are \emph{semi-isomorphic} if there is a bijection $f:S\ra T$ and an automorphism $\vp:F_r\ra F_r$ such that $f(x\cir g)=f(x)\cir \vp(g)$, for any $x\in S$ and $g\in F_r$. If we insist that the growth functions should be equivalent under semi-isomorphisms then, as we see, for instance, in \cite{SZ}, all functions with exponential growth fall into the same equivalence class. Using recent results in \cite{KKS} and \cite{KSS}, we give in Subsection \ref{ssSI} examples where $S$ and $S^{\prime}$ are semi-isomorphic, the growth of $S$ is maximal and the growth of $S^{\prime}$ is not.

\section{Growth of action}\label{sGA}

\subsection{Right actions}\label{ssRA}%{\bf 1. Right actions.} 
Here we recall few definitions and facts about the actions. By \emph{monoid} we will understand a semigroup with identity element $1$. A monoid $M$ \emph{acts} on the set $S$ if there is a \emph{structure map} $\mu: S\times M\ra S$, we write $\mu(s,x)=sx$, for $s\in S$ and $x\in M$, such that the following hold for any $s\in S$ and $x,y\in M$
\begin{enumerate}
	\item[\text{(1)}] $s(xy)=(sx)y$,
	\item[\text{(2)}] $s=s1$.
\end{enumerate}
Then $S$ is called a \emph{(right) act over $M$}.

An \emph{algebra} $R$ with $1$ over a field $\Phi$ acts on a linear space $V$ over $\Phi$ if there is a \emph{bilinear} structure map $V\times R\ra V$ satisfying (1) and (2). Of course, in this case $V$ is a (unital, right) $R$-module. In this paper all $M$-acts and $R$-modules will be right, so we silently assume this in all what follows. If a monoid $M$ acts on a set $S$, $R=\Phi M$ is the semigroup algebra of $M$ and $V=\Phi S$ is the linear space with basis $S$, then the action of $M$ on $S$ uniquely extends to the action of $R$ on $V$, and $V$ becomes an $R$-module.

Having in mind that the notions of acts are less familiar then those of modules, we quickly review some further material concerning these objects. Given $M$-acts $S$ and $T$, the map $\vp:S\ra T$ is called a \emph{morphism} of acts if for any $s\in S$ and $m\in M$ one has $\vp(sm)=\vp(s)m$. If $\vp$ has an inverse, $\vp^{-1}$, then $\vp^{-1}$ is also a morphism of $M$-acts, and in this case we say that $\vp$ is an isomorphism of $M$-acts $S$ and $T$. The subset $A$ is a \emph{generating set} of an $M$-act $S$ if $S=AM$. If $\# A=1$ then we call $S$ \emph{cyclic}. Any monoid $M$ is a cyclic act over itself, if one choses $\mu:M\times M\ra M$ to be the product in $M$. The identity element $1$ is the generator of $M$ as an $M$-act. 

Given an $M$-act $L$ with a nonempty generating subset $A$ we say that $A$ is a \emph{basis} of $L$ if for any $a_1,a_2\in A$ and $m_1,m_2\in M$ it follows from $a_1m_1=a_2m_2$ that $a_1=a_2$ and $m_1=m_2$. An act $L$ possessing a (necessarily unique, up to permutation of elements!) basis $A$ is called \emph{free}. If $\# A=s$ then $L$ is called \emph{free of rank $s$} and we write $\rk L=s$. It follows from the definition that any free act of rank $s$, $s$ finite or infinite, is isomorphic to the disjoint union of $s$ copies of the free act of rank 1, which is isomorphic to the $M$-act $M$. If $S$ is an arbitrary $M$-act and $\vp:A\to S$ is an arbitrary map then one can uniquely extend to a morphism of acts $\overline{\vp}:L\to S$. This latter property can be used as a ``more invariant'' definition of free acts.

\subsection{Growth functions}\label{ssGF} 
 Let us fix a finite generating set in a monoid $M$ with an ascending  filtration in $M$: $\{ 1\}=M(0)\subset M(1)\subset\cdots\subset M(n)\subset\cdots$ where $M(1)$ is the fixed generating set, each $M(n)$ is finite, $ M=\bigcup_{n=0}^{\infty}M(n)$, and $M(m)M(n)\subset M(m+n)$, for all $m,n=0,1,\ld$
Let $S$ be an $M$-act and $A$ a finite subset of $S$. We call $\mathcal{B}(A,n)=AM(n)$ the \emph{ball of radius $n$ around} $A$. The set-theoretic difference $\ms{A}{n}=\mb{A}{n}\setminus\mb{A}{n-1}$ will be called the \emph{sphere} of radius $n$ around $A$, $n\ge 1$. If $A$ is a one-element set $A=\{ a\}$ then we simplify our notation, write $\mathcal{B}(A,n)=\mathcal{B}(a,n)$ and call $\mathcal{B}(a,n)$ the \emph{ball of radius $n$ with center $a$}. Similarly, we write $\ms{A}{n}=\ms{a}{n}$. Since $M$ is an act over itself, $\mb{1}{n}=M(n)$ and $\ms{1}{n}=M(n)\setminus M(n-1)$, for all appropriate $n$. Every ball and every sphere are finite subsets of $S$. 

The chain of subsets: $A=\mathcal{B}(A,0)\subset \mathcal{B}(A,1)\subset\cdots\subset \mathcal{B}(A,n)\subset\cdots$ is an ascending filtration in the subact $T=AM$ generated by $A$ in $S$, in the sense that $T=\bigcup_{n=0}^{\infty}\mathcal{B}(A,n)$ and $\mathcal{B}(A,n)M(k)\subset\mathcal{B}(A,n+k)$. The \emph{growth function} $g_{A,T}$ is defined by setting $g_{A,T}(n)=\card \mathcal{B}(A,n)$. If $A=\{ a\}$ is a one-element set then we write $g_{A,T}=g_{a,T}$.

The same approach works in the case of modules over algebras. We only need to replace cardinalities of sets by dimensions of linear spaces. Specifically, let $\Phi$ be a field and $R$ a unital algebra over $\Phi$ with a fixed finite generating set and filtration: $\Phi\,\!.\,\! 1=R(0)\subset R(1)\subset\cdots\subset R(n)\subset\cdots$ where $R(1)$ is spanned by the generating set, each $R(n)$ is finite-dimensional, $R=\bigcup_{n=0}^{\infty}R(n)$, and $R(m)R(n)\subset R(m+n)$, for all $m,n=0,1,\ld$ Let $V$ be a right $R$-module, $A$ a finite subset of $V$. The space $\mathcal{B}(A,n)=AR(n)$ is called the \emph{ball of radius $n$ around} $A$ or the \emph{ball of radius $n$ with center $a$}, if $A=\{ a\}$.  In the latter case we write $\mathcal{B}(A,n)=\mathcal{B}(a,n)$. As in the case of monoids, we have $\mb{1}{n}=R(n)$. Every ball is finite-dimensional.

The chain of subspaces: $A=\mathcal{B}(A,0)\subset \mathcal{B}(A,1)\subset\cdots\subset \mathcal{B}(A,n)\subset\cdots$ is an ascending filtration in the submodule $U=AR$ generated by $A$ in $V$, in the sense that $U=\bigcup_{n=0}^{\infty}\mathcal{B}(A,n)$ and $\mathcal{B}(A,n)R(m)\subset\mathcal{B}(A,n+m)$. The \emph{growth function} $g_{A,U}$ is defined by setting $g_{A,U}(n)=\dim \mathcal{B}(A,n)$.

As mentioned in the Introduction, the action of $r$ transformations of any of the four types we consider gives rise to an ascending filtration in the respective universal algebra $\wax$, $\fgax$, $\fax$ and $\gax$. Let us call these filtrations standard and describe in greater detail.

The standard filtration on the free monoid $W=W(X)$ is given by the word length $|w|$ in the alphabet $X$. Let $X^m$ be the set of all words of length $m$ in $X$. Then $X^m=\ms{1}{m}$ and $W(n)=\mb{1}{n}=\bigcup_{m=0}^n X^m$. In the case of the free group $F=F(X)$ we need to consider the symmetrized set of generators $Y=X\cup X^{-1}$. Let $Y^m_{red}$ be the set of all words of length $m$ in $Y$ which are  reduced, that is, have no subwords $yy^{-1}$, where $y\in Y$. Then $Y^m_{red}=\ms{1}{m}$ and $F(n)=\mb{1}{n}=\bigcup_{m=0}^n Y^m_{red}$. Setting $\mathcal{A}(n)=\Phi W(n)$ and $\mathcal{F}(n)=\Phi F(n)$ defines standard filtrations in the \fa $\mathcal{A}=\fax$ and free group algebra $\mathcal{F}=\gax$, respectively.

The growth function $g_{ 1,W}$ of $W=W(X)$, with $\# X=r$, as the free act over itself can be easily computed considering that $\# X^m=r^m$, for $m=0,1,2,\ld$ Then $g_{1,W}(n)=\# W(n)=1+r+\cdots+r^n$. For the free act $L$ with basis $A$, $\# A=s$, we then have \begin{eqnarray}\label{grfrmo}g_{A,L}(n)=s(1+r+\cdots+r^n).
\end{eqnarray}
In the case of $F=F(X)$, with $\# X=r$, we have that $Y^0_{red}=\{ 1\}$, $Y^1_{red}=Y$ and if $Y^{m-1}_{red}$ has been determined, we will obtain the elements of $Y^{m}_{red}$, $m> 1$, each just once, if we multiply every word $u=u^\prime y\in Y^{m-1}_{red}$ by all letters of $Y$, except $y^{-1}$. Consequently, we have $\# Y^{m}_{red}=(2r-1)(\# Y^{m}_{red})=2r(2r-1)^{m-1}$, if $m\ge 1$. It then follows that $g_{ 1,F}(0)=\# Y^0_{red}=1$, $g_{1,F}(1)=1+(\# Y)=1+2r$, and $g_{1,F}(n)=g_{ 1,F}(n-1)+2r(2r-1)^{n-1}=1+2r(1+(2r-1) +\cdots+(2r-1)^{n-1})$, for $n>1$. For the free $F$-act $L$ with basis $A$, $\# A=s$, we then have
\begin{eqnarray}\label{grfrgr} g_{A,L}(n)=s(1+2r(1+(2r-1)+\cdots+(2r-1)^{n-1})).
\end{eqnarray}

Since the free modules of rank $s$ for $\fax$ and $\gax$ are the linear spaces whose bases are the free acts for $W(X)$ and $F(X)$ and their module structure is just the bilinear extension of the action of these latter monoids, the growth functions of the free module of rank $s$ over the free associative algebra of rank $r$ is given by (\ref{grfrmo}) and over the free group algebra of rank $r$ by (\ref{grfrgr}).

\subsection{Growth functions of actions over free monoids and \fas}\label{ssGFMFAA}

%%%%%%%%%%%%%%%%%%%%%%%%%%%%%%%%%%graphs of actions
We start with an observation which holds in all four cases of actions we study.

\begin{proposition}\label{pNCGF}
Let $g=g_{A,S}$ be the growth function of an action $S$ with finite generating set $A$ over one of $\wr$, $\fr$, $\ar$ or $\fgr$. Set $g(-1)=0$ and define $d(n)=g(n)-g(n-1)$, for $n\ge 0$. Then, for any $n\ge 1$, $d(n)\le rd(n-1)$, in the case of $\wr$ and $\ar$, and $d(1)\le 2rd(0)$, $d(n+1)\le (2r-1)d(n)$, in the case of $\fr$ and $\fgr$.
\end{proposition}

\pp The idea of the proof is the same in all four cases. One has to start with the sequence of balls $A=\mathcal{B}(A,0)\subset \mathcal{B}(A,1)\subset\cdots\subset \mathcal{B}(A,n)\subset\cdots$ in $S$. Using induction on $n$, one can select a subset $\mathcal{E}_n$ in $\mathcal{B}(A,n)$, consisting of some elements of the form $au$, $a\in A$ and $u\in X^n$ ($u\in Y^n_{reg}$, in the case of $\fr$ and $\fgr$), so that $(\bigcup_{m=0}^{n-1}\mathcal{E}_m)\cap\mathcal{E}_n=\emptyset$ and $\mathcal{E}=\bigcup_{n=0}^{\infty}\mathcal{E}_n$ is either $S$, in the case of $\wr$ and $\fr$, or the basis of $S$, in the case of $\ar$ and $\fgr$. In the case of acts we will have $\mathcal{E}_n=\mb{A}{n}\setminus\mb{A}{n-1}$ while in the case of modules, $\mathcal{E}_n$ will be a basis of $\mb{A}{n}$ modulo $\mb{A}{n-1}$. Therefore, $d(n)=\#\mathcal{E}_n$. If one sets $\mathcal{E}_0=A$ and assumes $\mathcal{E}_k$ defined for all $k=1,2,\ld,n-1$ then the elements of $\mathcal{E}_{n}$ can be selected among the products of all elements $au\in\mathcal{E}_{n-1}$ by all elements of $X$, in the case of $\wr$ and $\ar$, and by all elements of $Y$, except for $y^{-1}$ if $u=u'y$, in the case of $\fr$ and $\fgr$. Since $d(n)=\#\mathcal{E}_n$, we have the inequalities claimed in the statement of our proposition.\bx

Notice that in the case of $\wr$ and $\fr$ the subset $\mathcal{E}_n$ is the same as the sphere of radius $n$ centered at $A$.

In the case of $\wr$ and $\ar$ we have a converse to the previous proposition, as follows.

\begin{proposition}\label{AnyGrowth}
Let a nondecreasing function $g$ on the set $\{ 0,1,2,\ld\}$ take positive integral values. Consider $d(n)=g(n)-g(n-1)$ \emph{(}for convenience, we set $g(-1)=0$\emph{)}. Suppose for all $n=1,2,\ld$ we have $d(n)\le rd(n-1)$. Then $g$ is the growth function of a finitely generated act over a free monoid $W=\wax$ and a graded finitely generated module over a \fa $\mathcal{A}_r=\fax$ over a field $\Phi$.
\end{proposition}

\pp It is sufficient to construct a $W_r$-act $S$ act with such growth function. To obtain an $\mathcal{A}_r$-module, one simply can take a linear space $V$ with basis $S$ and naturally extend the action of $W_r$ on $S$ to an action of $\mathcal{A}_r$ on $V$. 

Let us assume $g(0)=m$. Then we can start with a free act $S$ with basis $A=\{ a_1,\ld,a_m\}$. The elements of $S$ are of the form $aw$, where $a\in A$ and $w\in W_r$. Our future act $T$ will be constructed by induction on the degree of $w$, as a subset of $S$ of the form $T=\cup_{n=0}^{\infty}T(n)$, where $T(n)\subset W(n)$, for the balls of radius $n$ around $A$ and $T_n\subset W_n$, for the spheres of radius $n$ around $A$. The action $\cir$ of $W$ on $T$ will appear in the process of construction.  

We start by setting $T_0=S_0=\{ A\}$. Proceeding by induction on $n$, we assume that we have already selected the elements of the ball $T(n)\subset W(n)$ and the sphere $T_n=\{ u_1,\ld,u_{d(n)}\}\subset W_n$ and defined the action of $W$ on $T_{n-1}$. We write $d(n+1)=pr+q$, where $0\le q<r$. Since $d(n+1)\le rd(n)$, we have that $p\le d(n)$, the inequality being strict if $q>0$. We set $u_i\cir x_j=u_ix_j$, for all $i=1,\ld,p$, $j=1,\ld,r$, and add these $u_ix_j$, to $T_{n+1}$. Then we define $u_{p+1}\cir x_1=u_{p+1} x_1, \ld,u_{p+1}\cir x_q=u_{p+1}x_q$ and add the elements $u_{p+1}x_1, \ld,u_{p+1}x_q$ to $T_{n+1}$. Finally,  we set $u_{p+1}\cir x_{q+1}=u_{p+1},\ld, u_{d(n)}\cir x_{r}=u_{d(n)}$. Then the ball $T(n+1)$ will have exactly $g(n+1)$ elements, and the action of $W$ on $T(n)$ has been defined. By induction both the elements of $T$ and the action of $W$ on them have been defined, and by construction, the growth of $T$ is as claimed.\bx

We complete this subsection with 

\begin{problem}\label{op0} Find necessary and sufficient conditions on a function $g(n)$ ensuring that $g(n)$ is the growth function of some action of a free group.
\end{problem}

\subsection{The growth as an invariant.}\label{ssGAI}%{\bf 3. The growth as an invariant.}  
If $S$ is a finitely generated $M$-act with a finite generating set $A$, $M$ a finitely generated monoid, and $B$ is a finite subset of $S$ generating subact $T=BM$ then there is a nonnegative integer $k$ such that $B\subset AM(k)$, and so for any $n$ we have $\mathcal{B}(B,n)\subset \mathcal{B}(A,n+k)$. Then $g_{B,T}(n)\le g_{A,S}(n+k)$. This inequality shows that if we want to produce a well-defined notion of the growth that does not depend on the choice of the finite generating set for $S$, then it is natural to proceed as follows. Given two growth functions $g_{B,T}$ and $g_{A,S}$, we say that $g_{A,S}$ \emph{majorates} $g_{B,T}$ if there is a nonnegative integer $C$ such that $g_{B,T}(n)\le g_{A,S}(n+C)$, for all $n=0,1,2,\ld$. If $T=S$, that is, $B$ is another generating set for $S$, then we say that $g_{B,S}$ is \emph{equivalent} to $g_{A,S}(n+k)$ if each of these two functions majorates the other. 

Speaking formally, let $\mathfrak{F}$ be the set of all nondecreasing functions $f:\{0,\,1,\,2,$ $\ld\}\ra \{1,\,2,\ld\}$. Given $f,g\in\mathfrak{F}$, we say that $f\preceq g$ (and say ``$f$ majorates $g$'') if there is a nonnegative integer $C$ such that $f(n)\le g(n+C)$, for all $n=0,1,2,\ld$. This is a pre-order relation in the sense that it satisfies only the reflexivity and transitivity axioms.  Putting $f\sim g$ if and only if $(f\preceq g)\,\&\,(g\preceq f)$, we obtain an equivalence relation on $\mathfrak{F}$. Finally, on the set $\mathfrak{F}/\!\!\sim$ of equivalence classes $[f]$ of functions $f$ from $\mathfrak{F}$ under $\sim$, we obtain a genuine partial order if we set $[f]\le[g]$ as soon as $f\preceq g$. Notice that $\sim$ is a \emph{congruenc}e in the sense that if $f\sim g$ and $f_1\sim g_1$ then $f+g\sim f_1+g_1$.

Coming back to the acts or modules, notice that their growth functions are in $\mathfrak{F}$. Suppose that $B$ is another finite generating set for an $M$-act $S$ with a finite generating set $A$, then, as we have seen, $g_{B,S}\preceq g_{A,S}$.  From the symmetry of $A$ and $B$, we have $g_{A,S}\preceq g_{B,S}$. So $g_{A,S}\sim g_{B,S}$. The equivalence class $g_S\in\mathfrak{F}/\!\!\sim$ containing all $g_{A,S}$, where $A$ runs through all the finite generating sets of $S$, is called the \emph{growth of an $M$-act $S$}. The growth is an \emph{invariant} of $S$, consequently, two acts with different growths cannot be isomorphic. 

Now suppose that we have a  morphism of $M$-acts $\vp:S\ra T$ and $P=\vp(S)$. Then $B=\vp(A)$ is a finite generating set for $T$ and the images of the balls around $A$ in $S$ are respective balls around $B$ in $P$: $\mathcal{B}(B,n)=\mathcal{B}(\vp(A),n)=\vp(\mathcal{B}(A,n))$. In this case, for any $n=0,1,2,\ld$, $g_{B,P}(n)\le g_{A,S}(n)$ and so $g_{B,P}\preceq g_{A,S}$. Thus, for the growths $g_P$ of $P$ and $g_S$ of $S$ we have $g_P\le g_S$. In particular, for the growth of an $s$-generator act over a monoid $M$ we always have $g_S\le g_L$ if $L$ is the free act of rank $s$.

Similarly, let $T$ be a finitely generated subact of a finitely generated $M$-act $S$. As noted above, for the finite generating sets $A$ and $B$ of $S$ and $T$, respectively, we would have $g_{B,T}\preceq g_{A,S}$. Thus for the growths we would have $g_T\le g_S$.

The case of finitely generated modules $V$ over $\Phi$-algebras $R$ is totally analogous to the case of acts. If $B$ a finite subset in an $R$-module $V$ generated by a finite set $A$, and $U$ is a submodule generated by $B$, then $B\subset AR(k)$, for some nonnegative integer $k$, and so  $\mathcal{B}(B,n)\subset\mathcal{B}(A,n+k)$. It follows that $g_{B,U}\preceq g_{A,V}$. If $U=V$ then by symmetry, $g_{B,V}\preceq g_{A,V}$ and $g_{A,V}\preceq g_{B,V}$. So $g_{A,V}\sim g_{B,V}$, and  so all $g_{A,V}$ are in the same equivalence class $g_V\in\mathfrak{F}/\!\!\sim$ called the \emph{growth of an $R$-module $V$}. As before, two $R$-modules with different growths are non-isomorphic.

Exactly the same argument as above allows one to conclude that if an $R$-module $V$ is a homomorphic image of a finitely generated $R$-module $U$ then $g_V\le g_U$, and if a finitely generated module $U$ is isomorphic to a submodule in a finitely generated $R$-module $V$ then $g_U\le g_V$.

A simple remark following from the formulas (\ref{grfrmo}) and (\ref{grfrgr}) is the following: $g_{\wr}=g_{\ar}=[r^n]$ and $g_{\fr}=g_{\fgr}=[(2r-1)^n]$.

One general result about the growth of modules is as follows.

\begin{proposition}\label{pwm} Let $R=\fax$, $\# X=r\ge 1$, be the free associative algebra over a field $\Phi$, $N$ is a submodule of finite codimension in a finitely generated infinite-dimensional $R$-module  $M$. Then $N$ is also finitely generated and the growth of $M$ is the same as the growth of $N$.
\end{proposition}

\pp We already know that $g_N\le g_M$ where $g_N$, $g_M$ are the growths of $N$ and $M$, respectively. To prove the converse we choose a finite generating set $A$ in $M$ and a finite generating system $C$ in $N$ in accordance with Schreier - Lewin procedure \cite{JL}. It follows then from the rewriting process of elements in $N$ in terms of $C$ that any element in $\mathcal{B}(A,n)\cap N$ is also an element in $\mathcal{B}(C,n)$. Then we have the following chain of inequalities: 
\bea
g_{C,N}(n) &=& \dim \mathcal{B}(C,n) \ge
\dim (\mcb(A,n) \cap N)\\ &=& \dim \mcb(A,n)- \dim (\mcb(A,n) /(\mcb(A,n) \cap N) \ge
g_{A,M}(n) -d,
\eqa
where $d= \dim M/N$. Now, since $N$ is infinite-dimensional, the values of $g_{C,N}$ grow at least by $1$ when we increase the argument by $1$. In this case, $g_{A,M}(n) \le g_{C,N}(n+d)$, that is, $g_{A,M} \preceq g_{C,N}$. It follows that $g_M\le g_N$, as needed.\bx

\section{Maximal growth}\label{ssMG}%{\bf 4. Maximal growth.} 
Using the above notions we can prove the following.
\begin{lemma}\label{lWMG} Let $S$ be a finitely generated act over the free monoid $\wr$ of rank $r>1$. Then
\begin{enumerate}
\item[\emph{(a)}] For any finite set $A$ of generators for $S$ there is $c>0$ such that $g_{A,S}\le cr^n$;
\item[\emph{(b)}] Given any finite set $A$ of generators for $S$, $g_{A,S}$ is majorated by $r^n$.
\end{enumerate}
The same is true for a finitely generated module over a free associative algebra $\ar$ of rank $r$.
\end{lemma}

\pp If $M$ is generated by $s$ elements then $S$ is an image under the morphism of the free $W_r$-act $L$ of rank $s$. So the $g_{A,S}\le g_{A,L}$, and we need to prove both (a) and (b) for $L$. The growth function of $L$ is given by (\ref{grfrmo}). So we have $g_{A,L}(n)=\dfrac{s}{r-1}(r^{n+1}-1)<cr^n$ with $c=\dfrac{sr}{r-1}$, proving (a).  If we take $C$ with $c\le r^C$ then $cr^n\le r^Cr^n = r^{n+C}$, for all $n=0,1,2,\ld$, and so $g_{A,L}(n)$ is majorated by  $r^n$.
\bx

If the reader uses (\ref{grfrgr}) in place of (\ref{grfrmo}) then the same argument proves the following result. 

\begin{lemma}\label{lFMG} Let $S$ be a finitely generated act over the free group $\fr$ of rank $r>1$. Then
\begin{enumerate}
\item[\emph{(a)}] For any finite set $A$ of generators for $S$ there is $c>0$ such that $g_{A,S}\le c(2r-1)^n$;
\item[\emph{(b)}] Given any finite set $A$ of generators for $S$, $g_{A,S}$ is majorated by $(2r-1)^n$.
\end{enumerate}
The same is true for a finitely generated module over a free group algebra $\fgr$ of rank $r>1$. \bx
\end{lemma}

If we fix $r>1$, then it is easily seen from the formulas (\ref{grfrmo}) and (\ref{grfrgr}),  that the function $r^n$ is majorated by the growth function of each free $s$-generator action, for all $s\ge 1$. On the other hand, as seen from Claim (2) of Lemma \ref{lWMG}, the converse is also true. Thus, with $r>1$ fixed, the growth of the free action of rank 1 is maximal among the growths of all finitely generated actions. Thus it makes sense to give the following

\begin{definition}\label{dFGMG}
 In any of the four cases $\wr$, $\fr$, $\ar$, and $\fgr$, we say that the growth of a finitely generated action is maximal if it is the same as the growth of the free action of rank 1.
 \end{definition}

\begin{lemma}\label{lMG} If $r>1$ then the growth of a finitely generated act $S$ (respectively, module $V$) over the free monoid $\wr$ (respectively, \fa $\fr$) of rank $r$  is maximal if and only if there is a finite generating set $A$ in $S$ (respectively, in $V$) and a positive $c>0$ such that $g_{A,S}(n)\ge cr^n$ (respectively, $g_{A,V}(n)\ge cr^n$) for all $n=0,1,2,\ld$. In this case the same inequality, with probably a different positive constant $c$, will hold for any other finite set of generators.
\end{lemma}

In the case of free groups or free group algebras, Lemma \ref{lMG} takes the following form.

\begin{lemma}\label{lfMG} If $r>1$ then the growth of a finitely generated act $S$ (respectively, module $V$) over the free group $\fr$ (respectively, the free group algebra $\fgr$) of rank $r$  is maximal if and only if there is a finite generating set $A$ in $S$ (respectively, in $V$) and a positive $c>0$ such that $g_{A,S}(n)\ge c(2r-1)^n$ (respectively, $g_{A,V}(n)\ge c(2r-1)^n$) for all $n=0,1,2,\ld$. In this case the same inequality, with probably a different positive constant $c$, holds for any other finite set of generators.
\end{lemma}

\pp Since the proof in the case of both lemmas is very similar, we restrict ourselves to the case of a finitely generated act $S$ over the free monoid $W$ of rank $r>1$. First, let us assume that the growth of $S$ is maximal. Then there is a finite generating set $A$ such that the growth function $g=g_{A,S}$ is equivalent to $r^n$. In particular, $r^n\preceq g$. So there is a nonnegative integer $C$ with $r^n\le g(n+C)$ or $g(n+C)\ge r^{-C}r^{n+C}$. As a result, $g(n)\ge r^{-C}r^n$, for all $n\ge C$. If we now choose a positive integer $c$ equal to the minimum of $r^{-C}$ and all $\frac{f(i)}{r^i}$, $i=0,1,\ld,C-1$, then $g(n)\ge cr^n$, for all $n=0,1,2,\ld$ Conversely, if there is $c>0$ with $g(n)=g_{A,S}(n)\ge cr^n$, for all $n=0,1,2,\ld$, then we need to show $[g]=[r^n]$. Indeed, the relation $f\preceq r^n$ is true by Lemma \ref{lWMG}(b). To prove that $r^n\preceq g$ we choose a natural number $C$ so that $cr^C\ge 1$. Then $r^n\le cr^{n+C}\le g(n+C)$ for any $n=0,1,2,\ld$, proving that, $r^n\preceq g$, as claimed, hence $[g]=[r^n]$, and thus the growth of $S$ is maximal.\bx

The reader will easily notice that actually our argument allows us to prove the following. Suppose we are given real $c$ and $C$ with $0<c<C$, and an integer $s>1$. Let $f,g\in \mathfrak{F}$ satisfy $cs^n\le f(n), g(n)\le Cs^n$, for all $n=0,1,2,\ld$ Then $f\sim g$.

One more result which  is important for dealing with the actions of maximal growth is the following.

\begin{lemma}\label{lAlpha} Let $S$ be an act generated by a finite set $A$ over $\wr$, $r>1$, with growth function $g=g_{A,S}$. We can write $g(n)=\alpha(n)r^n$, for a real-valued function $\alpha(n)$. Then $\alpha(n)$ is a function  converging to a finite limit $C_0$ at infinity. The growth is maximal if and only if $C_0 > 0$. The same claim holds if we consider the growth function of a finitely generated module $M$ over the \fa $\ar$, $r>1$.
\end{lemma}

\pp Let $A$ be a finite generating set for $S$, $\wr=\wax$. Then the ball of radius $n+1$ around $A$ can be written as 
\beq{balls}
&&\mathcal{B}(A,n+1)=\mathcal{B}(A,n)\cup\mathcal{B}(A,n)X\nonumber\\&=&\mathcal{B}(A,n-1)\cup\mathcal{B}(A,n-1)X\cup\mathcal{B}(A,n)X=\cdots=A\cup\mathcal{B}(A,n)X.
 \eqe 
So for the growth function $g=g_{A,S}$ we have $g(n+1)\le rg(n)+s$ where $s=\card A$.  Then we will have, for $g(n)$ and $\alpha(n)$ as in the statement of the theorem, 
\bea
\alpha(n+1)r^{n+1}\le s+r\alpha(n)r^n=s+\alpha(n)r^{n+1}\mbox{\ \  or\ \  }\alpha(n+1)-\alpha(n)\le\frac{s}{r^{n+1}}. 
\eqa
Now let us consider a function $\bar{\alpha}(n)=\alpha(n)+\frac{s}{r^n}$. Then from the previous inequality we will obtain
\bea
\bar{\alpha}(n+1)-\bar{\alpha}(n)\le\frac{2s}{r^{n+1}}-\frac{s}{r^{n}}=\frac{s(2-r)}{r^{n+1}}\le 0.
\eqa
Now we have that a positive-valued function $\bar{\alpha}(n)$ is monotonously non - increasing, hence there is $C_0\ge 0$ such that $C_0=\displaystyle{\lim_{n\to\infty}}\bar{\alpha}(n)$. It is easy now that $C_0=\displaystyle{\lim_{n\to\infty}}\alpha(n)$, as claimed. 

The claim that the growth is maximal if and only if $C_0>0$ follows from Lemma \ref{lMG}.

The proof for the modules is the same, except that the cardinalities should be replaced by dimensions and in formula (\ref{balls}) all signs $\bigcup$ should be changed to ``$+$''.
\bx

\begin{remark}\label{rPFA} Using Proposition \ref{AnyGrowth} one can show that for any $\ve>0$ it is possible to produce examples of acts and modules with growth function $\alpha(n)(r-\ve)^n$, such that any nonnegative real number is a limit point for the values of $\alpha(n)$.
\end{remark}

The proof of the following counterpart of the previous lemma for free groups and free group algebras is left to the reader as an easy exercise. Do not forget to replace $X$ by $Y=X\cup X^{-1}$ in formula (\ref{balls}) and $\mathcal{B}(A,i)X$ by the set of expressions $au$ where $a\in A$ and $u\in Y_{red}^j$ where $1 \le j\le i+1$.

\begin{lemma}\label{lGAlpha} Let $S$ be an act finitely generated by a set $A$ over $\fr$, $r>1$, with growth function $g=g_{A,S}$. We write $g(n)=\alpha(n)(2r-1)^n$. Then $\alpha(n)$ is a function  converging to a limit $C_0$ at infinity. The growth is maximal if and only if $C_0 > 0$. The same claim holds if we consider the growth function of a finitely generated module over the free group algebra $\fgr$, $r>1$.\bx
\end{lemma}

An easy consequence of Lemmas \ref{lAlpha} and \ref{lGAlpha} is this.

\begin{proposition}\label{lCAM} The growth of a finitely generated action over any of $\wr$, $\fr$ \emph{(}$\ar$,  $\fgr$\emph{)} of rank $r$ is maximal if and only if it has a cyclic subaction whose growth is maximal.
\end{proposition}

\pp If an action has a subaction of maximal growth then, as we know, the growth of the subaction cannot be greater that the growth of the action, and so the growth of the action must be maximal, too. Conversely, if we have a finitely generated action of maximal growth and none of its cyclic subactions has maximal growth then we can write it as the finite union (sum) of cyclic subactions, each of which is not of maximal growth. The growth function of the union (the sum) does not exceed the sum of the growth functions of the constituent cyclic subactons. By Lemmas \ref{lAlpha} and \ref{lGAlpha}, these latter growth functions have the form $\alpha(n)r^n$ in the case of $\wr$ and $\ar$ or $\alpha(n)(2r-1)^n$ in the case of $\fr$ and $\fgr$, with $\alpha(n)\to 0$. Thus the growth function of the whole of action has the same form, and so the growth of the act is not maximal.\bx

\begin{definition}\label{dMG}
An arbitrary (not necessarily finitely generated) action over any one of $\wr$, $\fr$, $\ar$ or $\fgr$ of rank $r>1$ is called an action of \emph{maximal growth} if it contains a \emph{finitely generated} subaction of maximal growth. 
\end{definition}

Observe that it follows by Proposition \ref{lCAM}, that an \emph{action has maximal growth if and only if it contains a \emph{cyclic} subaction of maximal growth}. 

\begin{corollary}\label{cSNMG} In each of the four cases we consider, every act (module) has the largest subact (submodule) which is not an act (module) of maximal growth. This subact (submodule) is the union (sum) of all subacts (submodules) whose growth is not maximal, or equivalently, of all cyclic subacts (submodules) whose growth is not maximal.
\end{corollary}

\pp This is immediate from Definition \ref{dMG} and Proposition \ref{lCAM}.\bx

In the case of modules over $R=\ar$ or $\fgr$ we will prove more (Theorem \ref{te} below): namely, if we denote by $\mathcal{N}(M)$ the largest submodule of an $R$-module $M$ whose growth is not maximal, then \emph{every} nonzero submodule in $M/\mathcal{N}(M)$  has maximal growth. This makes $\mathcal{N}(M)$ looking like a regular \emph{radical} of modules.

Let us stipulate, for the future, that if an act $S$ and its generating set $A$ are fixed then we will often omit indexes $A$ and $S$ in the notation for the growth function $g_{A,S}(n)$, and simply write $g(n)$. Similar convention will be used for the modules.

\begin{problem}\label{op1}
Is it true that any graded module of maximal growth over the free associative algebra $\ar$ or free group algebra $\fgr$ possesses graded factor-modules of arbitrary possible growth?
\end{problem}

\subsection{Growth and co-growth}\label{ssGC}

According to Grigorchuk \cite{RG}, the definition of amenability for finitely generated groups is equivalent to the following.
Let $F=\fgax$ be the free group of  rank $r$, $N$ a normal subgroup in $F$. The co-growth function for $N$, denoted by $c_N(n)$ is the number of elements in the intersection of $N$ with the ball $F(n)$ in (Cayley graph of) $F$, that is, the number of reduced words in free generators of $F$ of length at most $n$ which are in $N$ (relations for $F/N$). Then  $F/N$ is \emph{amenable} if an only if 
\bea
\lim_{n\ra\infty}\sqrt[n]{c_N(n)}=2r-1.
\eqa
 
Now let $L$ be a free module with basis $A=\{ a_1,\ld,a_s\}$ over one of the algebras $R= \fax$, or $R=\gax$, $r>1$. Let $N$ be a submodule in $L$ and $M = L/N$.  We define the co-growth function $c_{A,N}(n)$ of  $M$  with respect to $A$ by setting $c_{A,N}(n) = \dim (N \cap \mathcal{B}(A,n))$, where $\mathcal{B}(A,n)$ is the ball of radius $n$ in $L$ around $A$. As before in Subsection \ref{ssGAI}, the co-growth functions corresponding to different bases are equivalent, so we can speak about the \emph{co-growth} of $M=L/N$. 
 
\begin{lemma}\label{clN} Let $L$ is the free module with basis $A$ over one of the algebras $R=\ar,\fgr$, $r>1$, $N$ a submodule in $L$, $M=L/N$, $\bA$ is the image of $A$ in $M$ under the natural epimorphism of $L$ onto $M$. Then
 \bea
g_{A,L} = g_{\bA,M}+c_{A,N}.
 \eqa
 \end{lemma}
 
 \pp
Indeed, let $E$ be a basis of $\mathcal{B}(A,n) \cap N$ and $E^{\,\prime}$ its complement to a basis of $L(n)$, where $n=0,1,2,\ld$ Then  \bea
\card E^{\,\prime} = \dim \mathcal{B}(A,n)/(L(n)\cap N) = \dim (\mathcal{B}(A,n)+N)/N = g_{\bA,M}(n).
\eqa
However by definition, $\card E = c_{A,N}(n)$. It follows that $g_{A,L}(n) = (\card E) + (\card E^{\,\prime}) = c_{A,N}(n) + g_{\bA,M}(n)$, as claimed.\bx

\begin{proposition}\label{cCG}
Let $L$ is a free module with basis $A$ over one of the algebras $R=\ar,\fgr$, $r>1$, $N$ a submodule in $L$. If the growth of $M=L/N$ is maximal then there is $0<\theta<1$ such that its co-growth function $c_{A,N}$ satisfies
\bea
\frac{c_{A,N}(n)}{g_{A,L}(n)}  <\theta.
\eqa
If the growth of $M$ is not maximal then
\bea
\lim_{n\ra\infty} \frac{c_{A,N}(n)}{g_{A,L}(n)}=1.
\eqa
\end{proposition}
 
\pp  Since the proofs in $\ar$- and $\fgr$-cases are similar, we restrict ourselves to the case of modules over free associative algebras.
By formula (\ref{grfrmo}), $g_{A,L}(n)<2sr^n$; if $M$ has maximal growth, then by Lemma \ref{lMG}, we have $g_{\bA,M}(n) >cr^{n}$. In this case, $\dfrac{c_{A,N}(n)}{g_{A,L}(n)} =1-\dfrac{g_{\bA,M}(n)}{g_{A,N}(n)}< \theta$, where $\theta = 1- \dfrac{c}{2s}$.

Now suppose the growth is not maximal. Then, as we proved in Lemma \ref{lAlpha}, $g_{\bA,M}(n) =\alpha(n)r^n $, where $\alpha(n)\to 0$. By formula (\ref{grfrmo}) we have $g_{A,L}(n)\ge r^n$. So it becomes obvious that $\dfrac{c_{A,N}(n)}{g_{A,L}(n)}\to 1$, and thus the proof is complete. \bx

\subsection{Faithfulness of actions with maximal growth}\label{sFAMG}

In this subsection we will establish the faithfulness of actions of maximal growth in all for cases we study. In some arguments we will be using the standard total ordering on the set of words, called \emph{ShortLex}. An action is called \emph{faithful} if any two different elements of the acting monoid or algebra act differently. The faithfulness follows quite easily from the following ``folklore'' property of \emph{languages}.

Given an alphabet $X$, any subset $L$ of the free monoid $W(X)$ is called a \emph{language}. One of the ``folklore'' facts about the languages is that if $w$ is a nonempty word in the alphabet $A$ consisting of more than one letter and $L=N(w)$ is a language consisting of words which have no occurrences of $w$ as a subword then the growth of $L$ is exponentially slower than the growth of $W(X)$. If $Y=X\cup X^{-1}$ is a symmetric (``group'') alphabet and the languages consist of reduced words this ``folklore'' property remains valid. In the group case we can refer to a particular paper \cite{AS}, which contains the proof. In the case of monoids, the proof is very short, and we give it for completeness.

We start with giving a precise statement of the ``folklore'' property of languages. Recall that $W(n)$ and $F(n)$ stand for the $n^{\mathrm{th}}$ terms of the standard filtration in $\wax$ and $\fgax$ (see Subsection \ref{ssGF}).

\begin{lemma}\label{lsubword} Let $X$ be an alphabet, $\# X=r>1$, $Y=X\cup X^{-1}$. Suppose $W=\wax$, respectively, $F=\fgax$, is the free monoid, respectively,  the free group, both with basis $X$ and standard filtrations defined in Subsection \ref{ssGF}.
\begin{enumerate}
\item[\emph{(a)}]  Let $u$ be a nonempty word in alphabet $X$, $N(u)$ the set of words in $W$ which do not have $u$ as a subword. Then there are positive numbers $C$ and $\ve$,  such that 
\bea
\card(N(u)\cap W(n))\le C(r-\ve)^n\mbox{ for all }n=1,2\ld .
\eqa
\item[\emph{(b)}] Let $v$ be a nonempty reduced word in alphabet $Y$, $N(v)$ the set of reduced words in $F$ which do not have $v$ as a subword. Then there are positive numbers $C$ and $\ve$,  such that 
\bea
\card(N(v)\cap F(n))\le C(2r-1-\ve)^n\mbox{ for all }n=1,2\ld .
\eqa
\end{enumerate}
\end{lemma}

\pp Let $\deg u = m$. Then the number of monomials of degree $m$ not containing $u$ equals $r^m-1$. Hence the number of monomials of degree $ms$ not containing $u$ is at most $(r^m-1)^s$ because such monomials are products of s factors of degree $m$. It follows that the number of monomials of degree $n$ not containing $u$ is at most
\bea
(r^m-1)^{\frac{n}{m}}r^{m-1}\le C(r^{\prime})^n
\eqa
$\mbox{ where }C = r^{m-1}\mbox{ and } r^{\prime}\mbox{ is the $m$-th root of }r^m-1 : r^{\prime} < r.$
\bx

The results about the faithfulness of actions are valid under weaker restrictions on the growth functions. We say that the function $f:\mathbb{N}\ra\mathbb{N}$ is \emph{subexponential} if 
\bea
\overline{\lim_{n\to\infty}}\sqrt[n]{f(n)}\le 1.
\eqa 

\noindent\textit{1. The case of acts.} An act $S$ over monoid $M$ is called \emph{faithful} if for any two different $u,v\in M$ there is $s\in S$ such that $su\neq sv$.

\begin{proposition}\label{pFAMG} Let $S$ be a finitely generated act over the free monoid $W=\wr$, $r>1$, with the growth function $g(n)$. If  the function\quad $\dfrac{r^n}{g(n)}$\quad  is subexponential then $S$ is faithful. In particular, any act of maximal growth is faithful.
\end{proposition}

\pp Assume by contradiction, that there are two different elements $u$ and $v$ in $W$ that act in the same way. We assume that $u>v$ in ShortLex. We write $S$ as the union of a finite number of cyclic subacts $S=P_1\cup\ldots\cup P_m$. If $a_i$ is a generator of $P_i$, $i=1,\ld,m$, then $\mathcal{B}( a_i,n)\subset a_i(N(u)\cap W(n))$, where $N(u)$ is the set of words in $W$ which do not have $u$ as a subword. Indeed, if $w=w_1uw_2\in W$, then by our assumption, $a_iw=a_iw_1uw_2=a_iw_1vw_2$, and the word $w^\prime=w_1vw_2$ is smaller than $w$ in ShortLex. If $w^\prime\in N(u)$, our claim is proven, otherwise we continue replacing $u$ by $v$. Since the ordering is total, at some step of the process we arrive at $a_i\tilde{w}\in a_iN(u)$. By Lemma \ref{lsubword} we can find constants $C_i>0$ and $\ve_i>0$ such that $g_{a_i,P_i}(n)=\card\, \mathcal{B}(a_i,n)\le C_i(r-\ve_i)^n$  for all $n=1,2\ld$. Since $g(n)\le g_{a_1,P_1}(n)+\cd+g_{a_m,P_m}(n)$, for some $C>0$ and $\ve>0$ we have $g(n)\le C(r-\ve)^n$  for all $n=1,2\ld$. Now for the function $g(n)$ in the statement of the proposition, we obtain
$\sqrt[n]{g(n)}\ge\dfrac{r}{r-\ve}\dfrac{1}{\sqrt[n]{C}}$, so the upper limit in question exceeds $\dfrac{r}{r-\ve}>1$, a contradiction.

Now let $S$ be an act of maximal growth, $u$, $v$ are two different elements of $W$. By definition, $S$ has a finitely generated subact $T$ of maximal growth. By Lemma \ref{lMG}, there is a positive constant $c>0$ such that $g(n)=g_{a,T}(n)\ge cr^n$. For the function $g(n)$, as in the statement of the proposition, we have $g(n)\le \dfrac{1}{c}$, $\sqrt[n]{g(n)}\to 0$, and so by our previous argument $T$ is faithful. It follows that $u$ and $v$ cannot act in the same way on $T$, hence on $v$. So, $S$ is faithful. \bx

Similarly, in the case of acts over the free group $F$, or $F$-sets, we prove the following.

\begin{proposition}\label{pFFSMG} Let $S$ be a finitely generated $F_r$-set over the free group $F_r$ of rank $r>1$, with the growth function $g(n)$. If the function \quad  $\dfrac{(2r-1)^n}{g(n)}$ \quad  is subexponential then $S$ is faithful. In particular, any $F_r$-set of maximal growth is faithful.\bx
\end{proposition}

An important corollary \cite{AS} is as follows.

\begin{corollary}\label{pFANMG} Let $N$ be a nontrivial proper normal subgroup of a free group $F_r$ of rank $r>1$. Then the growth of the $F_r$-set $G=F_r/N$  is not maximal.\bx
\end{corollary}

\noindent\textit{2. The case of modules.} A right module $M$ over a ring $R$ is faithful if the two-sided ideal $\mathrm{Ann}_R M=\{ a\in R\: |\: Ma=0\}$ is zero. We are going to show that every module of \mg is faithful. Again, as in the case of acts, a slightly more general result holds.

\begin{proposition}\label{pFMMG} Let $M$ be a finitely generated module over $R=\ar$ with the growth function $g(n)$. If the function\quad  $\dfrac{r^n}{g(n)}$\quad is subexponential then $M$ is a faithful $R$-module. In particular, any module of maximal growth is faithful.
\end{proposition}

\pp Arguing in the same way as in Proposition \ref{pFAMG}, we may restrict ourselves to the case where $M$ is a cyclic module. Any such module can be written as $M=R/J$ where $J$ is a right ideal of $R$. However, if $M$ is not faithful, $J$ must contain a nonzero two-sided ideal $I$. Since the growth of a factor-module of a module cannot be greater than the growth of the module, any upper bounds we obtain for the growth of $R/I$ will be valid for the growth of $M$. Suppose $0\neq a\in I$. We denote by $u$ the leading term of $a$ with respect to ShortLex. Then any monomial $v$ containing $u$, that is, $v=v_1uv_2$ for some monomials $v_1$ and $v_2$, can be reduced to a smaller monomial modulo $I$ :  
\bea
v_1uv_2 +I = \sum v_1u^{\prime}v_2 +I\mbox{ where }u^{\prime} <u.
\eqa
It follows by induction on ShortLex that the ball of radius $n$ around $\{ 1+I\}$ in $R/I$ is the linear span of all $w+I$ where $w$ is a monomial not containing $u$. By Lemma \ref{lsubword} the number of such monomials is exponentially smaller than the total number of all monomials of degree $r$, and so the function $\dfrac{r^n}{g_{1+I,M/I}(n)}$ is exponential, and, as noted, then $\dfrac{r^n}{g(n)}$ is exponential, in contradiction with our hypotheses.

Now if $M$ is a module of maximal growth generated by a finite set $A$ then by Lemma \ref{lMG} there is a positive constant $c>0$ such that $g(n)=g_{A,M}(n)\ge cr^n$. In this case, as in Proposition \ref{pFAMG}, the function $f(n)=\dfrac{r^n}{g(n)}$ is subexponential and by the above, $M$ is faithful. \bx

Notice that the class of faithful modules is much wider than the class of modules of \mgp Faithful modules can be found among the modules of an arbitrary infinite growth, for instance, modules of linear growth. As an example, one can take the submodule $N$ in the construction of Section \ref{sENM}.  For each $i=1,2,\ld$, this module has a submodule generated by $f_{\vp(i)}x_2R$, whose annihilator is $\Delta^{d_i+1}$. Since $d_i\ra\infty$, the total $\mathrm{Ann}_RM=0$.

In the case of free group algebras we get the following result whose proof follows the lines of Proposition \ref{pFMMG}.

\begin{proposition}\label{pFMMGFGA} Let $M$ be a finitely generated module over $R=\fgr$, $r>1$, with the growth function $g(n)$. If the function of \quad  $\dfrac{(2r-1)^n}{g(n)}$ \quad  is subexponential then $M$ is a faithful $R$-module. In particular, any module of maximal growth is faithful.\bx
\end{proposition}

Notice that in the case of the action of $\fr$, $r>1$, on a linear space $V$, that is, when $V$ becomes an $\fgr$-module, the faithfulness of an $\fgr$-module is a stronger property than the faithfulness of the representation of the group $\fr$. The nontriviality of the kernel in $\fr$, obviously, implies the nontriviality of the kernel in $\fgr$ but not the other way round. For example, $\fgr$ has an ideal of codimension 4 with trivial intersection with $\fr$ as soon as the field $\Phi$ is not locally finite. This follows because $\fr$ can be faithfully represented by $2\times 2$-matrices over $\Phi$.)

\section{Acts of maximal growth}\label{sRCAMG}

In this section we give examples of acts of maximal growth which nevertheless satisfy very restrictive conditions. Some of them readily provide examples of modules over \fas which have maximal growth and satisfy other interesting properties.
 
\subsection{Highly transitive acts of maximal growth}\label{sAMG}

In this section we will construct acts with maximal growth over the free monoid $\wr$ of rank $r>1$, which have interesting additional properties. In the following definition, $S$ is a right act over a monoid $M$.

\begin{definition}$S$  is called \emph{$k$-transitive} if for any $2k$-tuple  $(s_1,...,s_k,s'_1,\ld,s'_k)$ of elements in $S$ such that all $s_1,\ld,s_k$ are pairwise different, there is $m\in M$ such that $s_i\, m = s'_i$, for all $i=1,\ld,k$.
\end{definition}

\begin{theorem}\label{tKTAMG} Let $W_r$ be the free monoid of rank $r>1$.
There exists a right act $S$ over $\wr$ with the following properties
\begin{enumerate} 
\item[\rm(a)] $S$ has maximal growth;
\item[\rm(b)] for any $k$, $S$ is $k$-transitive.
\end{enumerate}
\end{theorem}

\pp  Let $W_r=\wax$ be the free monoid of words in the alphabet $X$, $\# X=r>1$. Choose an infinite language $P\subset W_r$ satisfying the following condition.
   
   \medskip
  
($\dagger$) \textit{If a suffix of $w\in P$ is a prefix of a  $w'\in P$, then $w=w'$.}

\medskip

An example of such a language $P$ in the alphabet $X=\{ x,y\}$ is provided by 
\bea
W =\{ x^2(yx)^t y^2\:|\:\mbox{ where } t=0, 1,... \}.
\eqa

Let us write out all nonempty tuples $v = (v_1,\ld,v_k;\,v_1^{\prime},\ld,v_k^{\prime})$ of elements in $W_r$, for any $k\ge 1$ such that $v_1,\ld,v_k$ are pairwise different and nonempty:
\begin{eqnarray*}
v(1)&=&(v(1,1),\ld,v(1,k(1));\,v(1,1)^{\prime},\ld,v(1,k(1))^{\prime}),\\v(2)&=&(v(2,1),\ld,v(2,k(2));\,v(2,1)^{\prime},\ld,v(2,k(2))^{\prime}),\\&&\ld\\
v(n)&=&(v(n,1),\ld,v(n,k(n));\,v(n,1)^{\prime},\ld,v(n,k(n))^{\prime}),\\&&\ld
\end{eqnarray*}
   
   With each tuple $v(i)$ occurring on the $i^\mathrm{th}$ position in this list we associate a word $w_i\in P$, $i=1,2\ld$. All words $w_i$ are pairwise different, each of length at least  $i + k(i)$.

Now let $U$ be the set of all products $u(i,j)=v(i,j) w_i$, for all possible values $i,j=1,2,\ld$, such that $j\le k(i)$.

By $S$ we denote the set of words without prefixes in $U$.
  
We will make the language $S$ into a right act for $W_r$. To draw distinction between the product in $W_r$ and the action of $W_r$ on $S$, we denote the action by  $\cir$. Since $W_r$ is free with basis $X$, any well-defined action of $X$ on $S$ extends to a well-defined action of $W_r$. 
   
   Let  $s \in S$, $x\in X$. If $sx\notin U$, then we set $s\cir x = sx$.  If $sx = u(i,j)$ then $s \cir x = v'(i,j)$. 
   
 This action is well-defined. Indeed, if we assume $v(i,j)w_i = v(i',j')w_{i'}$ then by Property ($\dagger$), $i = i'$, hence $v(i,j)=v(i,j')$. Since the first $k(i)$ letter in the tuple $v(i)$ are pairwise different, we must have $j=j'$, which proves our claim.

The $k$-transitivity of the action of $W_r$ on $S$, for any $k=1,2,\ld$, can be explained as follows. Let us consider any tuple of elements of $S$ of the form $(v_1,\ld,v_k,v'_1,\ld,v'_k)$ such that all $v_1,\ld, v_k$ are pairwise different. Then there is $i=1,2,\ld$ such that this tuple equals \bea
v(i) = (v(i,1),\ld,v(i,k(i)),v'(i,1),\ld,v'(i, k(i))).
\eqa 
Let us consider  $v(i,j)\cir w_i$. Since $v(i,j)=v_j\in S$, it has no prefixes in $U$. Also, by Property ($\dagger$) none of the proper prefixes of $v(i,j)w_i$ is in $U$. Therefore, in process of acting by all consecutive letters of $w_i$, except the last one, we use the rule of the first kind, $s\cir x =sx$. But when we reach the last letter, we have to apply the rule of the second kind which replaces $u(i,j)w_i$ by $v'(i,j ) = v'_j$. Thus, $v_j\cir w_i = v'_j$, for all $j=1,\ld, k$. This proves the $k$-transitivity of our action.  
   
In order to estimate the growth of $S$, we will estimate the number of words in the $r$-letter alphabet $X$ whose length is at most $n$ and which do not have any of $u(i,1),\ld,u(i,k(i))$ as their prefixes. The length of each $u(i,j)$ is at least $i+j$, by the choice of $w_i$. In this case at most $r^{n-i-j}$ of words of length $\le n$ begins with $u(i,j)$, hence at most $r^{n-i}$ words begins with any of the words $u(i,1),\ld,u(i,k)$. If we perform summation over all $i\le n$, we will see that in the number of elements in the ``ball'' of radius $n$ in $S$ is at least $r^n +r^{n-1}+\cd - r^{n-1} -r^{n-2}-\cd =r^n$. By Lemma \ref{lMG} the growth of $S$ is maximal.\bx

A quick application of this theorem to modules over \fas is the following.

\begin{corollary}\label{cSMMG} Let $\ar$ be the \fa of rank $r>1$ over a field $\Phi$. Then $\ar$ has a simple module whose growth is maximal.
\end{corollary}

\pp Let us choose a linear space $M$ whose basis is the $W_r$-act $S$ from Theorem \ref{tKTAMG}. Here $W_r=\wax\subset\fax=\ar$ is the free monoid of the same rank $r>1$. If we extend the action of $W_r$ on $S$ by linearity to the action of $\ar$ on $M$, then $M$ becomes a right $\ar$-module. Since the growth of $W_r$-act $S$ is the same as the growth of $\ar$-module $M$, we know that $M$ is an $\ar$-module of maximal growth. Let $N$ be a subspace in $M$ consisting of finite $\Phi$-linear combinations $\sum_{u\in S} \lambda_u u$ with $\sum_{u\in S} \lambda_u =0$. Being a $W_r$-invariant subspace, $N$ is an $\ar$-submodule. Since $N$ has codimension 1 in $M$, by Proposition \ref{pwm} we have that $N$ is an $\ar$-module of maximal growth. Let us prove that $N$ is a simple $\ar$-module. Any nonzero element $a\in N$ can be written as the sum of $m$ summands $\lambda_u e_u$ where $\lambda _u\ne \lambda_v$ for some $u,v$ (since $M$ is infinite-dimensional, one of coefficients can always be chosen zero). If we use the $m$-transitivity of the action of $W_r$ on the basis, we will obtain that a submodule $Q$ generated by $a$ also contains the sum $b$ which differs from $a$ by changing places $u$ and $v$. Hence the difference $a-b =\lambda(u-v)$ is an element of $Q$ with $\lambda =\lambda_u-\lambda_v \ne 0$.  From the double transitivity it follows that $u - v \in Q$ now for any $u,v\in S$, and then $Q=N$, proving the simplicity of $N$.\bx

\subsection{Acts with maximal growth and locally finite action of ``skinny'' submonoids}\label{ssLSM}

The aim of this subsection is to construct an act of maximal growth over $W_r$, $r>1$, such that the action of all ``truly'' smaller submonoids is locally finite, even locally nilpotent in the sense we define below.
 
Let $W_r$ be a free monoid of rank $r>1$ and $V$ a submonoid. We will say that $V$ is \emph{fat} in $W_r$ if there are $t\ge 1$ and $h_1,\ld,h_t,g_1,\ld,g_t\in W_r$ such that 
\bea
W_r= h_1 Vg_1\cup...\cup h_tVg_t.
\eqa
If no such  $h_1,\ld,h_t,g_1,\ld,g_t$ exist then the submonoid $V$ is called \emph{skinny}. For example, the submonoid $U$ of all words of even length is fat in $W_r$. On the other hand, if $r>1$ and $V$ consists of all words except the positive powers of $x_1$, then $V$ is skinny. Notice that $V$ is not finitely generated.

To formulate the main result of this subsection, we recall that given a monoid $V$, a $V$-act $S$ is called \emph{locally finite} if for any $s\in S$ the subact $s\cir V$ (we call it the $V$-``orbit'' of $s$, by analogy with the case of group actions) is finite. If, additionally, $S$ has 0, that is, a unique distinguished element 0 such that $0\cir M=0$, we call $S$ \emph{locally nilpotent} if for any $s\in S$ there is a natural $n$ such that  $((s\cir m_1)\cir m_2)\cir \cdots \cir m_n =0$, for all non-identity elements  $m_i \in M$. It is quite obvious that if $V$ is a fat submonoid in $W_r$ whose action on a $W_r$-act $S$ is locally finite then also the action of $W_r$ on $S$ is locally finite. This easily follows since each $W_r$-orbit can we written as $sW_r = \bigcup ((sg_i) V)h_i$. In contrast to this, the following is true.

\begin{proposition}\label{pSKINNY}
For any natural $r>1$ there exists an act $S$ with zero $0$, of maximal growth over the free monoid $W_r$, which is a locally finite, even locally nilpotent $V$-act, for any finitely generated skinny submonoid $V\subset W_r$. In particular, for any $s\in S$ and $u\in W_r$ there is natural $n$ such that $s\cir u^n=0$.
\end{proposition}
 
\pp If $V$ is skinny then, for any $C> 0$ there are words $w\in W_r$, which have the following property

\medskip

\qquad $P(C): $\hspace{1cm} \emph{$w$ cannot be written as $w=hvg$, where $v\in V$ and $|h|+|g| \le C$}. 

\medskip

The converse is also true.

Now let $V$ be a nontrivial finitely generated skinny and $l=l(V)\ge 1$ be such that $V$ can be generated by the words of length $\le  l$. Pick $w\in W_r$ which does not have $P(2l-2)$. Then $w$ is not a subword of any $v\in V$. Indeed, otherwise, $v=v_1\cdots v_s$, where $v_1,\ld, v_s$ are in the generating set for $V$. It follows that $w =v'_{i-1}(v_i\cdots v_j) v'_{j+1}$, where $v'_{i-1}$  (respectively, $v'_{i+1})$ is a proper suffix (prefix) of the word $v_{i-1}$ (respectively, $v_{j+1}$), and so $w$ has $P(2l-2)$.

It follows that $V$ is contained in the set $N(w)$ of all words of $W_r$ without subword $w$. Using Lemma \ref{lsubword}, it is easy to show that there exist $\ve>0$ and natural $n_0$, such that $\#(V\cap W(n))
< (r-\ve)^n$, for all $n \ge n_0$. Here $W(n)$ is the $n^\mathrm{th}$ term of the standard filtration in $W_r$. It should be noted that, conversely, if $V$ grows at such ``slow'' rate then it is skinny, no matter whether the number of generators of $V$ is finite or not; this is obvious considering the difference in the growths of sequences $\{ \# V(n)\}$ and $\{ \# W(n)\}$, where $V(n)=V\cap W(n)$. 

With $l,\ve$ as above, we fix $c>0$ and choose the number $t>n_0$ sufficiently large to satisfy the inequality $(r-\ve)^{t+1}+...+(r-\ve)^{t+l}<cr^t$. 

Now let us introduce a collection of languages attached to $V$, as follows. Choose any $u\in W_r$ and consider the language $L(u,V)$ consisting of all words in $W_r$ which have a prefix of the form $uv$ where $v\in V$ and $t< |v|\le t+l$. It then follows that $\#(L(u,V)\cap W(n))$ is bounded from above  by $cr^n$, for all possible $n$.

Let us enumerate all pairs $(u,V)$, where $u$ is a word and $V$ is a finitely generated skinny  submonoid. For the $i^\mathrm{th}$ pair $(u_i, V_i)$ by the above argument we will have $\# (L(u_i,V_i) \cap W(n))\le \dfrac{c}{2^i}\, r^n$. Therefore the number of words of length $n$ in the language $L$, which is the union of all $L(u_i,V_i)$, grows slower than  $cr^n$. It follows that the number of words of length $n$ in the complement $K$ of $L$ grows faster than $c'r^n$, for some $c'>0$, provided that $c<1$. As a result, we conclude that there is $c'>0$ such that $\# (K \cap W(n))>c'r^n$ for all natural $n$, which implies that the growth of $K$ is maximal.

To define our desired act $S$, we take the union of the complement $K$, and a one-element set $\{ 0\}$, where $0$ is an external element for $S$. Following the same pattern as in Theorem \ref{tKTAMG}, given $s\in K$ and $x\in X$, we set $s\cir x =sx$, if $sx \in K$, and $s\cir x = 0$ otherwise, that is, if $sx\in L)$. Finally, we set $0\cir x=0$. This makes $S$ into a $W_r$-act. By construction, this act has maximal growth.

Now let $V$ be a finitely generated skinny submonoid, $u\in S$ and $u\cir V$ an ``orbit'' of $V$ in $S$. Let $t=t(V)$ and $l=l(V)$ be the numbers appearing in our discussion of the pair $(u,V)$ above. Any nonzero element in $u\cir V$ either belongs to a finite set $\{ u\cir v \,|\, |v|\le t\}$ or has the form of $u\cir v$ where $v = v_1v_2$, where $v_1 \in V$,  $t<|v_1|\le t+l$. From the construction of the language $L$ above, it follows that $uv \in L$ and the $u\cir v = 0$. As a result, any orbit $u\cir V$ is finite. The same argument shows that $u$ acted upon by any product of more than $t$ elements of $V$ produces 0. So $S$ is locally nilpotent with respect to the action of $V$. \bx

We do not have analogues of this result in other cases of actions. In the case of groups, a subgroup $H$ is fat in $F_r$ if 
\bee{eCC}
F=(g_1Hg_1^{-1})(g_1h_1)\cup\ldots\cup (g_sHg_s^{-1})(g_sh_s).
\ene 
B.H.Neumann \cite{BHN} has shown that if a group is covered by finitely many cosets of a set of subgroups then one of the subgroups in the set has finite index in the groups. In our case, it follows from \ref{eCC} that $H$ is fat if $H$ is of finite index in $F_r$. It is quite obvious that if there is a finite orbit for a subgroup $H$ of finite index then the same is true for the whole of $F_r$.

\begin{problem}\label{op7}
Does there exist an $F_r$-set of maximal growth (or an infinite transitive $F_r$-set) such that the restriction of the action to any finitely generated subgroup of infinite index is locally finite?
\end{problem}

  Let us call a subalgebra $B$ with $1$ in $\mathcal{A}_r$ \emph{fat} if $\mathcal{A}_r$ equals to the span $UBV$ of all products $ubv$, $u\in U$, $b\in B$, $v\in V$, for some finite-dimensional subspaces $U,V\subset \mathcal{A}_r$. As an example one may take the linear span of monomials whose degrees are divisible by a fixed number $n$. If $B$ with $1$ is not fat, it is called \emph{skinny}. It is easy to observe that if all cyclic $B$ submodules are finite-dimensional for some fat subalgebra $B$ then $M$ is finite-dimensional. In particular, $M$ cannot be of maximal growth as an $\mathcal{A}_r$-module.

\begin{problem}\label{op8}
  Does there exist an  $\mathcal{A}_r$-module of maximal growth (or a cyclic infinite-dimensional $\mathcal{A}_r$-module) which is locally finite as a module over any skinny subalgebra of $\mathcal{A}_r$?
\end{problem}
%%%%%%%%%%%%%%%%%%%%%%%%%%%%%%%%%%%%%%%%%%%

\section{Modules of maximal growth}\label{sgcfa}

\subsection{Radical defined by growth}\label{ssRDMG}
 Let $R$ be a ring. In \cite[Chapter 5]{PMC}, two classes of modules over $R$, $\mathfrak{M}$ and $\mathfrak{N}$, were called the \emph{annihilators} of each other if
\bea
\mathfrak{N}=\{U\,|\,\mathrm{Hom}(U,V)=0\;\forall\, V\in\mathfrak{M}\},\;
\mathfrak{M}=\{V\,|\,\mathrm{Hom}(U,V)=0\;\forall\, U\in\mathfrak{N}\}.
\eqa
If the classes $\mathfrak{M}$ and $\mathfrak{N}$ are annihilators of each other in the above sense then one calls $\mathfrak{M}$ a \emph{semisimple class} and $\mathfrak{N}$ a \emph{radical class} of $R$-modules.

Let $R$ be one of $\fax$ or $\gax$, with $\# X>1$. Denote by $\mathfrak{M}$ the class of $R$-modules in which every nonzero submodule has maximal growth and by $\mathfrak{N}$ the class of $R$-modules in which no submodule has maximal growth. 
From the general results and definitions of Sections \ref{sGA} we know the following.
\begin{enumerate}
\item If an $R$-module $U$ is mapped  onto an $R$-module $V$ with maximal growth then the growth of $U$ is also maximal;
\item If a submodule $W$ of an $R$-module $U$ has maximal growth then $U$ itself has maximal growth;
\item An $R$-module has maximal growth if and only if one of its finitely generated submodules has maximal growth if and only if one of its cyclic submodules has maximal growth. 
\end{enumerate}

It follows easily that the classes $\mathfrak{M}$ and $\mathfrak{N}$ defined by us are the annihilators of each other and so can serve as respective semisimple and radical classes of $R$-modules. 

As we mentioned in Section \ref{sGA}, in every action there is the largest subaction whose growth is not maximal. So in every module $M$ there is the largest submodule $\mathcal{N}(M)\in \mathfrak{N}$. We now want to show that $M/\mathcal{N}(M)\in \mathfrak{M}$.

\begin{theorem}\label{te} Let $R$ be either a free associative algebra $\fax$ or a free group algebra $\gax$ of rank $r>1$, $M$ an $R$-module of maximal growth, $N$ a submodule of $M$. Then at least one of $N$, $M/N$ is an $R$-module of maximal growth. In other words, the radical class $\mathfrak{N}$ is closed under extensions.
\end{theorem}

\pp Let $\#X=r$. The proof is the same in both cases, with maximal growth $r^n$ in the case of free associative algebras and $(2r-1)^n$ in the case of free group algebras. So we give the proof in the case of free associative algebras. Let us assume that none of $N$, $M/N$ is an $R$-module of maximal growth. Choose a cyclic submodule $uR$ in $M$ whose growth is maximal. Then neither $uR\cap N$ nor $uR/uR\cap N\cong uR+N/N\subset M/N$ is an $R$-module of maximal growth. This allows one to restrict oneself to the case where $M=uR$ is cyclic. If $R(n)$ is the ball of radius $n$ with center $1$ in $R$ and $uR+N/N$ is not of maximal growth, then by Lemma \ref{lAlpha} we will have
\bea
\lim_{n\to\infty}\frac{\dim (uR(n)+N/N)}{r^n}=0.
\eqa
Let us define a subspace $U(m)$ by setting $U(m)=uR(m)\cap N$. Obviously, $U(m)$ is finite-dimensional and so the submodule $U(m)R$ is finitely generated. If for some $m$, $U(m)R$ has \mgc then by Proposition \ref{pwm} we have a cyclic submodule of \mg in $N$, and the proof is complete. Otherwise, by Lemma \ref{lAlpha}, we have
\bea
\lim_{n\to\infty}\frac{\dim (U(m)R(n))}{r^n}=0,
\eqa
for any fixed $m$. Let us consider a subspace $V(m)$ such that $uR(m)=U(m)\oplus V(m)$. Then $\dim V(m)=\dim (uR(m)+N/N)$ and so $(\dim  V(m)/r^m\rightarrow 0$ by the above property (2) of submodules and our assumption about $M/N$. Notice that $\dim uR(k)=\dim U(k)+\dim V(k)$, for all $k$. Choose any $\ve > 0$. Then there is $m>0$ such that $\dim  V(m) <\ve r^m$. Since $U(m)R$ is not of \mg there is $n$ such that $\dim U(m)R(n) < \ve r^n$ for all $n >n_0$.  Now we fix $n\geq n_0$ and let $t$ be any integer with $t >m+n$ such that $r^{-t+m+n} < \ve$. Let $R_i=X^i$, for $i=0,1,\ld$. Then
\bea
uR(t)=uR(m)(R_0+R_1+\cdots+R_{t-m}).
\eqa
Then we will obtain
\begin{eqnarray*}
\dim uR(t)&=&\dim uR(m)(R_0+R_1+\cdots+R_{t-m})\\
&\le& \dim V(m)(R_0+R_1+\cdots+R_{t-m})\\ &+& \dim U(m)(R_n+R_{n+1}+\cdots+R_{t-m})\\
&+&\dim U(m)(R_0 +R_1+\cdots+R_{n-1})\\
&<& \ve r^m(1+r+\cdots+r^{t-m})\\&+&(\ve r^{n}+\cdots+\ve r^{t-m})+\dim uR(m+n-1)\\
&<& \ve r^m(1+r+\cdots+r^{t-m})+\ve r^{n}(1+r+\cdots+r^{t-m-n})\\&+&(1+r+\cdots+r^{m+n-1})\\
&=& \ve r^m\frac{r^{t-m+1}-1}{r-1}+\ve r^n\frac{r^{t-m-n+1}-1}{r-1}+\frac{r^{m+n}-1}{r-1}\\
&<& \ve r^t(r+r^{-m+1}+\ve^{-1} r^{-t+m+n})<
\ve r^t(r+r +1)=\ve(2r+1)r^t.
\end{eqnarray*}
Since $\ve(2r+1)$ is arbitrarily small, we have that
\bea
\lim_{t\rightarrow\infty}\frac{\dim uR(t)}{r^t}=0.
\eqa
and so $M=uR$ is not of \mgc which contradicts our hypothesis.
\bx 

\begin{corollary}\label{cRMG} Let $R$ be either a free associative algebra $\fax$ or a free group algebra $\gax$ of rank $r>1$. Then in every $R$-module $M$ there is the largest submodule $N=\mathcal{N}(M)$ whose growth is not maximal. At the same time, the growth of every nonzero submodule of $M/N$ is maximal.
\end{corollary}

\pp Let $N$ be the largest submodule in $M$ whose growth is not maximal. Its existence has been established in Corollary \ref{cSNMG}. The growth of any submodule $P$ of $M$ which is not in $N$ is maximal. Since $P\cap N$ is not of maximal growth, by Theorem \ref{te}, we must have that $P+N/N\cong P/P\cap N$ is of \mgp Since any nonzero submodule of $M/N$ has this form, our claim follows.\bx

\subsection{Growth and short exact sequences}\label{sENM}

The results in this section should be compared with Theorem \ref{te}: starting with two modules whose growth is ``very slow'' we can construct their extension whose growth is ``arbitrarily fast'' but not maximal, which would be impossible by Theorem \ref{te}. Before we give the precise statement of our claim we recall that according to Lemma \ref{lAlpha}, the growth function of any finitely generated module over the \fa of rank $r>1$ has the form of $\alpha(n)r^n$ where  $\alpha(n)\to C_0\ge 0$ as $n\to \infty$. The growth is maximal if and only if $C_0>0$. The precise statement is as follows.

\begin{proposition}\label{pEMLG} Let $\alpha:\mathbb{N}\ra\mathbb{N}$ be a function satisfying $\alpha(n)\to 0$ as $n\to \infty$. Then there is a cyclic module $M=\ve_1R$ over a \fa $R=\mathcal{A}\langle x_1,\ld,x_r\rangle$, $r> 1$, with a cyclic submodule $N$ such that the growths of $N$ and $M/N$ are linear while for the the growth function of $M$ we have $g_{\ve_1,M}(n)>\alpha(n)r^n$, for all sufficiently large $n$.
\end{proposition}

\pp Given $\alpha(n)$ as mentioned, there is an increasing sequence $d_1,d_2,\ld$ of natural numbers such that $\alpha(n)<\dfrac{1}{r^3}$ for all $n\ge d_1$, $\alpha(n)<\dfrac{1}{r^4}$ for all $n\ge d_2$, etc. We will also need a function $\vp:\mathbb{N}\ra\mathbb{N}$ given by $\vp(1)=1$ and $\vp(n)=4r^{d_n}$, for all $n>1$. 

Let us form a linear basis for $M$ as the union of  
$\{\ve_i\:|\:i\in\mathbb{N}\}$, $\{\eta_i\:|\:i\in\mathbb{N}\}$, and 
$\{\zeta_{\vp(i),u}\:|\:i\in\mathbb{N},\: u\mbox{ any monomial in }R\mbox{ of degree at most }d_i\}$.

We define the action of $R$ on this basis as follows:
\begin{eqnarray*}
\ve_ix_1&=&\ve_{i+1},\: \ve_ix_2=\ve_i+\eta_{\vp(i)},\:\ve_ix_k=0\mbox{ for any }i\in\mathbb{N}\mbox{ and }k=3,\ld,r\\
\eta_jx_1&=&\eta_{j+1},\: \eta_jx_2=\left\{\begin{array}{cc}0&\mbox{ if }j\notin\vp(\mathbb{N})\\\zeta_{j,1}&\mbox{ if }j\in\vp(\mathbb{N})\end{array}\right.,\:
\eta_jx_k=0,\:j\in\mathbb{N},\:k=3,\ld,r\\
\zeta_{\vp(i),u}x_k&=&\left\{\begin{array}{cc}0&\mbox{ if }\deg{u}=d_i\\\zeta_{\vp(i),ux_k}&\mbox{ if }\deg{u}<d_i\end{array}\right.,\mbox{ for any }i\in\mathbb{N}\mbox{ and }k=1,\ld,r.
\end{eqnarray*}

It is an easy remark that $M$ is cyclic, $M=\ve_1R$. Indeed, applying $x_1$ to $\ve_1$ repeatedly, we obtain all the $\ve_i$, $i\in\nat$. Since $\ve_1x_2=\ve_1+\eta_1$ and $\eta_jx_1=\eta_{j+1}$ we have that all $\eta_j$, $j\in\nat$, are also in $\ve_1R$. Next, $\zeta_{\vp(i),1}=\eta_{\vp(i)}x_2$ are in $\ve_1R$ and, finally, $\zeta_{\vp(i),u}=\zeta_{\vp(i),1}u\in \ve_1R$, for all monomials in $R$ of degree at most $d_i$.

If we set $N=\eta_1R$ then the basis of $N$ will consist of all $\eta_i$ and $\zeta_{\vp(i),u}$. As we just noted, $M$ (hence $M/N$) are cyclic modules. Let us evaluate the dimension of $\mb{\ve_1}{n}$, the ball of radius $n$ in $M$, where $d_{i-1}< n\le d_i$, $i=2,3,\ld$. Clearly, $\eta_{\vp(i)}=\ve_ix_2-\ve_i$ is in the ball of radius $i=i(n)$. Hence $\zeta_{\vp(i),1}=\eta_{\vp(i)}x_2$ is in the ball of radius $i+1$. It follows that in the ball of radius $n$ we find all $\zeta_{\vp(i),u}$ where $\deg u\le n-i-1$. All these vectors are nonzero because $n\le d_i$. Their number is at least $r^{n-i-1}$. Hence for the growth function $g=g_{\ve_1,M}$ we have $g(n)\ge r^{n-i-1}$. On the other hand, we have $\alpha(n)<\dfrac{1}{r^{i+1}}$, following from $n\ge d_{i-1}$. Hence $g(n)\ge r^{n-i-1}\ge r^n\alpha(n)$, as needed. 

Let us evaluate the growth of $N$. We consider $\mb{\eta_1}{n}$, the ball of radius $n\ge 4r^{d_2}$ in $N$. It obviously contains $\eta_1,\ld,\eta_{n+1}$ and also some of the vectors of the form $\zeta_{\vp(i),u}$ but only if $\vp(i)\le n$. It follows then that
\bea
g_{\eta_1,N}(n)\le n+1+2(r^{d_1}+\cdots+r^{d_i})\mbox{ where }i\mbox{ is maximal with }\vp(i)\le n. 
\eqa
Let us additionally assume $i\ge 2$. Then it follows from $d_1<d_2<\ld$ that $2(r^{d_1}+\cdots+r^{d_i})<4r^{d_i}=\vp(i)\le n$, and so $g_{\eta_1,N}(n)< 2n+1$, for $n\ge 4r^{d_2}$, that is, this growth is linear.

The factor-module $M/N$ has basis $\{\bar{\ve}_i= \ve_i+N\,|\,i\in\nat\}$ with action $\bar{\ve}_ix_1=\bar{\ve}_{i+1}$, $\bar{\ve}_ix_2=\bar{\ve}_{i}$ and $\bar{\ve}_ix_k=0$, for all $i\in\nat$, $k=2,3,\ld$. Then $\mb{\bar{\ve}_1}{n}$ has basis $\{\bar{\ve}_1,\ld,\bar{\ve}_{n+1}\}$. Thus $g_{\bar{\ve}_1,M/N}(n)=n+1$, that is, the growth of $M/N$ is also linear.\bx

Another interesting extension is provided by an example in Proposition \ref{pNNMG}: if we fix $\alpha(n)$ as just above then there exists a cyclic module $M$ which has a locally finite submodule $N$ such that the growth of $M/N$ is linear and, as in the previous example, the growth of $M$ itself is greater than $\alpha(n)r^n$.

%%%%%%%%%%%%%%%%%%%%%%%%%%%%%%%%%%%%%%%%%%%%%%%%
%\subsection{Simple module of maximal growth}\label{sSMMG}

In the following subsections we give example of modules of maximal growth which satisfy strong finiteness conditions.

\subsection{Nil modules of maximal growth}\label{sNMMG}

We recall that a \fa $R=\fax$ has a standard grading by the subspaces $\Phi X^n$ each of which is spanned by all monomials in $X$ of degree $n$. We will denote these subspaces by $R_n$, $n=0,1,2,\ld$.

A module $M$ over a \fa $R=\fax$ of rank $r>1$ is called \emph{graded} if  $M=\bigoplus_{n=0}^{\infty}M_n$, where each $M_n$ is a subspace and $M_nR_m\subset M_{n+m}$, for all $m,n=0,1,2,\ld$. Now let
$P_n$ be a subspace of dimension $d_n$ in each $M_n$, $n=0,1,\ld$ We consider an $R$-submodule $L$ in $M$ generated by all $P_n$ and set $L_n=L\cap M_n$. Then the following weak analogue of E. Golod's Lemma \cite{G} is true.
 
 \begin{lemma}\label{clG} 
 \bea
 \dim  L_ n \le \sum_{i\le n}  d_i r^{n-i}.
 \eqa
 \end{lemma}
 
\pp
 $L_n$ is the sum of subspaces $P_i R_{n-i}$ with  $i\le n$, whose dimension is at most $d_i r^{n-i}$, proving our claim.
 \bx
 
  \begin{lemma}\label{l1234}  Let $M=aR$ be a cyclic graded module and $C$ a positive number such that $g_{a,M}(n)\ge Cr^n$. Choose $c$ with $0<c< C$. Let $b= a\cdot u$, where $u$ is a monomial in $R$. We consider an arbitrary selection $v_1,\ld,v_k$ of monomials of degree at least 1 in $R$. Then there exists a natural $q$,  and a graded submodule $L$ in $M$, such that the growth function of the graded $R$-module $M/L$ is bounded from below by $cr^n$, and $b(s_1v_1+\cdots+s_kv_k)^q \in L$ for any choice of scalar coefficients $s_1,\ld,s_k$.
\end{lemma}
 
 \pp  We choose $q$, so that $2q^k \le (C-c)  r^q$. 
Then we write the polynomial $(s_1v_1+\cdots+s_kv_k)^q$ as the sum of homogeneous polynomials $w_{m_1,\ld,m_k}$ of degree at least $q$, with coefficients  $s_1^{m_1}\ld s_k^{m_k}$, so that the coefficients inside $w_{m_1,\ld,m_k}$ do not depend on $s_1,\ld, s_k$. Since $m_i \le q$, the number of elements in the set of all $w_{m_1,\ld,m_k}$ is at most $q^k$.
Let us define $L$  as a submodule generated by all  
$b\cdot w_{m_1,\ld,m_k} = au\cdot w_{m_1,\ld,m_k}$. Then, obviously, $b\cdot (s_1v_1+\cdots+s_kv_k)^q \in L$ for all choices of coefficients $s_1,\ld,s_k$.
 
Let us apply Lemma \ref{clG} to the space $P$ spanned by the monomials $w_{m_1,\ld,m_k}$. Let $P_n=P\cap R_n$ and suppose $d_n=\dim P_n$. As mentioned, $\dim P\le q^k$ and so $\sum d_n\le q^k$. Also, $P_n=\{ 0\}$ if $i<q$. Now by Lemma \ref{clG}, the dimension of the $n^\mathrm{th}$ homogeneous component $L_n$ of $L$ is at most 
 \bea 
\sum_{i\ge q}  d_i r^{n-i} \le q^k r^{n-q}.
 \eqa
Now suppose 
\bea
M(m)=\bigoplus_{n=0}^m M_n\mbox{ and }L(m)=L\cap M(m)=\bigoplus_{n=0}^m L_n.
\eqa
Then it follows that 
\bea
\dim L(m) \le  q^k \sum _{m=q}^n r^{m-q} \le 2q^k r^{n-q}.
\eqa
In the factor-module $M/L$ the generator is $\bar{u}=u+L$ and the ball $\mb{\bar{u}}{n}$ equals $M(n)/L(n)$. If we recall the choice of $q$ then for the growth function of $M/L$ we can write: 
\bea
  g_{\bar{u},M/L}(n)=\dim M(n)/L(n) \ge Cr^n - 2q^k r^{n-q}=   r^n\left(C - \frac{2q^k}{r^q}\right)
\ge c r^n,
\eqa 
proving our claim about the lower bound for the growth function of $M/L$. \bx

Now we can prove our result about nil-modules of maximal growth.
 
 \begin{theorem}\label{tNMMG} Let $R$ be a free associative algebra over a field $\Phi$ of rank $>1$ and $M$ a graded $R$-module of maximal growth. Then $M$ has a graded factor-module of maximal growth, which, in addition, is a nil-module.
 \end{theorem}
 \pp Let us enumerate all finite tuples of monomials  $(u,v_1,\ld,v_k)$. By Lemma \ref{l1234}, the growth function $g_M(n)$ for $M$ satisfies  $g_M(n)\ge Cr^n$ for some $C> 0$. We choose any strictly decreasing infinite sequence $c_1>c_2>\ld$ in the interval $(D, C)$ where $D$ is a fixed number, $D\in (0, C)$. 
 
By Lemma \ref{l1234}, one can build a sequence of submodules $L^{(1)}\subset L^{(2)}\subset \ld$ such that if a tuple $(u,v_1,\ld,v_k)$ occurs on the $i^\mathrm{th}$ place then for any choice of coefficients $s_1,\ld,s_k$,  
\bea
au\cdot (s_1v_1+\cdots+s_kv_k)^{n_i}\in L^{(i)}
\mbox{ for some }n_i.
\eqa
The growth function $g_i=g_{M/L^{(i)}}$ of the respective factor-module satisfies $g_i(n)>c_i  r^n$, for all $n=0,1,2,\ld$ and for each $au\in M$ and each $s_1v_1+\cdots+s_kv_k\in R$ there are positive integers $i$ and $n_i$ such that  $au\cdot (s_1v_1+\cdots+s_kv_k)^{n_i}=0$ in $M/L^{(i)}$. Now let us set $L = \bigcup_{i=1}^{\infty} L^{(i)}$ and $N=M/L$. For each $n$ then there is $i$ such that $L(n)=L^{(i)}(n)$. As a result, if $g$ is the growth function of $N$, then for each $n$, $g(n)= g_i(n)\ge c_ir^n>Dr^n$, proving that the growth of the module constructed by us is indeed maximal.

To complete the proof, it remains to notice that any $b\in M$ is a linear combination of the elements of the form $a\cdot u_1,
\ld, a\cdot u_l$, where $u_1,\ld,u_l$ are the monomials; therefore,  $b\cdot (s_1v_1+\cdots+s_kv_k)^ n=0$ in $M/L$ where $n$ is the maximal exponent chosen for the tuples $(u_1, v_1,\ld,v_k),\ldots,(u_l, v_1,\ld,v_k)$
on the respective steps of the construction.

The proof is complete.\bx

As noted in Introduction, the following result cannot be obtained from  Golod - Shafarevich's example of an infinite-dimensional finitely generated nil-algebra because those algebras viewed as modules over free associative algebras do not have maximal growth.
 
 \begin{corollary}\label{cNMMG} Any free associative algebra of rank $r>1$ has a cyclic graded nil module whose growth is maximal.\bx
 \end{corollary}
 
 The conclusion of Theorem \ref{tNMMG} that $M$ with maximal growth has a  nil factor module of maximal growth, that is, the same growth, is not true for any growth which is not maximal. 
 
 \begin{proposition}\label{pNNMG} Let $\alpha(n)$ be a positive real valued function such that $\alpha(n)<1$ and  $\displaystyle\lim_{n\to 0}\alpha(n)=0$. Then there exists a graded cyclic $R$-module $M =aR$ whose growth function $g_M$ satisfies  $g_M(n)\ge  \alpha (n)r^n $, and which does not have nil factor-modules of infinite dimension. Moreover, if $m$ is a natural number and $\vp$ is a homomorphism of $M$ such that $\vp(a)\, x_1^m=0$ then $\vp(M)$ is finite-dimensional.
 \end{proposition}
  
  \pp  Choose a sequence $d_1, d_2,\ld$ of natural numbers so that $\alpha(x)\le\dfrac{1}{r^{i+1}}$ for all $x\ge d_i$, $i=1,2,\ldots$ and a linear space $M$  with a basis, which is the union of two subsets: $\{\ve_i\:|\:i=1,2,\ldots\}$ and \bea
\{\xi_{i,u}\:|\:i=1,2,\ldots;\;u\mbox{ is a monomial in }R\mbox{ of degree at most }d_i\}.
\eqa
The action of the generators $x_1,\ld,x_r$ of $R$ is given by $\ve_i\, x_1 = \ve_{i+1},  \ve_i\, x_2 = \xi_{i,1}$, $\ve_ix_j=0$ for $j>2$, and
\bea  
\xi_{i,u}\, x_k = \left\{\begin{array}{cc}\xi_{i,ux_k}&\mbox{ if }\deg\,u<d_i\\0&\mbox{otherwise}\end{array}.\right.
\eqa

Then $M=\ve_1R$. Let us prove that the growth of $M$ is greater than $\alpha (n)r^n$. Since $\alpha(n)\to 0$ and $\alpha(n) \le 1$, given $n$ there is an integer $j\ge 0$ such that $1/r^j \ge \alpha (n) > 1/r^{j+1}$. By the choice of $d_j$, it follows that $d_j > n$. Now in the case where $n \le j$ we have 
\bea
g(n) \ge 1 = (1/r^n) r^n\ge (1/r^j)r^n \ge \alpha(n)r^n.
\eqa
In the case where $j<n$, we notice that $\xi_{j,1} = \ve_1x_1^{j-1}x_2\in\mb{\ve_1}{n}$, the ball of radius $j$. Hence $\deg u \le n-j$ implies that each $\xi_{j,u}=\xi_{j,1}u\in\mb{\ve_1}{n}$, the ball of radius $n$. It follows that $g(n) >r^{n-j} \ge \alpha(n)r^n$, as needed.

Now if we impose any relation $\ve_1\,x_1^m =0$, then in the factor-module we will have $\ve_{m+1} = \ve_{m+2}=\cdots=0$, hence $\xi_{i,u} = 0$ for  $i > n$. So the module obtained by imposing just one ``nil''-relation of this form is already finite-dimensional, as claimed.\bx

%%%%%%%%%%%%%%%%%%%%%%%%%%%%%%%%%%%%%
\subsection{Residually finite modules of maximal growth}\label{sRFM}

Let $R=\far{r}$, where $r >1$. For each $j=1,\ld,r$ we choose an infinite by $i$ sequence 
$(\alpha_{ji})$ of elements of the field $\Phi$. For any monomial $v=x_{j_1}\cdots x_{j_d}$ of degree $d$, $1\le j_1,\ld,j_d\le r$,  we will define a ``quasi-monomial'' $e_v = (x_{j_1}- \alpha_{j_1,1})\cdots (x_{j_d}-\alpha_{j_d,d})$ of degree $d$. We denote by $\Delta_i$ the (``quasi-augmentation'') right ideal of $R$ generated by all $(x_1-\alpha_{1i }),\ld,(x_r - \alpha_{ri})$.
 
\begin{lemma}\label{l991} The following are true.
\begin{enumerate}
\item[\emph{(a)}] All $e_v$ form a linear basis in $R$.
\item[\emph{(b)}] All ideals $\Delta_i$ are two-sided.
\item[\emph{(c)}] The  quasi-monomials $e_v$ with $\deg v \ge m$ form a basis in the product of ideals $\Delta_1\cdots \Delta_m$.
\end{enumerate}
\end{lemma}
\pp \begin{enumerate}
\item[(a)] Notice that $v = e_v +w$ where $\deg w <d=\deg(v)$. Using induction by degree, we obtain that $R$ is a linear span of all $e_v$. The linear independence follows because the leading terms  (with respect to ShortLex) of $e_v$ is $v$.
 \medskip
\item[(b)] Follows because for each $i$ all $(x_1-\alpha_{1i }),\ld,(x_r - \alpha_{ri})$ generate $R$ as an algebra with 1.
\medskip 
\item[(c)] Since the $i^\mathrm{th}$ factor in the definition of $e_v$ is in $\Delta_i$, we have that $e_v  \in   \Delta_1\cdots \Delta_m$ if $\deg v \ge m$. Conversely, as it follows from (b), any element in $\Delta_1\cdots \Delta_m$ is in the right ideal generated by all $e_v$, $\deg v =d \ge m$.
Since we have
\bee{(*)}
e_v x_j = e_{vx_j} + \alpha_{j,d+1} e_v,
\ene
the subspace spanned by all $e_v$ with $\deg v \ge m$ is a right ideal and all that remains is to refer to our claim (a).
\end{enumerate}\bx

\begin{lemma}\label{l992} Let $\{v_1,v_2,\ld\}$ be a set of monomials in $R$, $\deg v_i = d_i$, for $i=1,2,\ld$. We denote by $I$ the right ideal in $R$ generated by all $e_v$ as above such that $v\in\{v_1,v_2,\ld\}$. Let $g(n)$ be the growth function of $M=R/I$ with respect to the generator $1+I$.  Then
\bea
g(n) > 1+r+\cdots+r^n - 2\sum _{j:\, d_j\le n} r^{n-d_j}.  
\eqa
\end{lemma}

\pp It follows from (\ref{(*)}) that $I$ is the linear span of those exactly $e_w$, for which $w$ begins with one of $v_1,v_2,\ld$  We denote the set of all such $w$ by $V$ and the subset of all monomials of degree at most $n$ in $V$ by $V(n)$. The number of such monomials with $\deg w=l\le n$ is at most $\sum _{j:\, d_j\le l} r^{l-d_j}$, meaning that 
\bea
\card\, V(n)<2 \sum_{j:\, d_j\le n} r^{n-d_j}.
\eqa

It follows from Claim (a) of Lemma \ref{l991} that all the rest quasi-monomials of degree at most $n$ are linearly independent modulo $I$. Since the number of all quasi-monomials of degree at most $n$ equals the number of all monomials of degree at most $n$, that is, $1+r+\cdots+r^n$, and all $e_v+ I$ are in the ball of radius $n$ in $M$, the proof of the lemma is complete.\bx
 
Let us enumerate all quasi-monomials $u_1, u_2,\ld$ in accordance with ShortLex of their leading terms. Then we choose an increasing sequence of natural numbers $1<d_1<d_2<\ld$ such that $\deg u_i\le d_i$ and form a sequence of quasi-monomials 
\bea
w_i = u_i(x_1 -\alpha_{1,\deg(u_i)+1})\cdots(x_1 - \alpha_{1,d_i}).
\eqa

Finally, we denote by $I$ the right ideal of $R$ generated $w_1, w_2,\ld$.

\begin{lemma}\label{l993} The following hold for the module $M=R/I$.
\begin{enumerate}
\item[\emph{(a)}] $M$ has maximal growth.
\item[\emph{(b)}] For each $w \in M$ there are natural numbers $s$ and $t$ such that $s\le t$ and $w(x_1-\alpha_{1,s})\cdots(x_1-\alpha_{1,t})=0$. If  $w \in M_i= M\Delta_1\Delta_2\cdots \Delta_i$ then $s \ge i$.
\end{enumerate}
\end{lemma}

\pp \begin{enumerate}\item[(a)] By Lemma \ref{l992}, for the growth function $g(n)$ of $M$ we have 
\begin{eqnarray*}
  g(n) &>& 1+r+\cdots+r^n - 2\sum _{j:\, d_j\le n} r^{n-d_j}\\ &\ge&
1+r+\cdots+r^n - 2 (r^{n-2} + r^{n-3}+\cdots)\\&=&r^n+r^{n-1}-r^{n-2}\cdots-1>r^n.
\end{eqnarray*}
According to Lemma \ref{lMG}, the growth of $M$ is, indeed, maximal.
\medskip
 \item[(b)] Choose $w \in  \Delta_1\Delta_2\cdots \Delta_i$ for some $i\ge 0$. By Lemma \ref{l991},
$w$ is a linear combination of quasi-monomials  $e_v$, $\deg v \ge i$. For each such $e_v$ we have $e_v = u_j$ for some $j \ge i$. If $s_v = \deg u_j +1, t_v = d_j$, then $e_v(x_1- \alpha_{1,s_v})\cdots(x_1-\alpha_{1,t_v})=u_j(x_1- \alpha_{1,\deg u_j})\cdots(x_1-\alpha_{1,d_j})=w_j \in I$. Let $s$ be the minimum of all $s_v$ and $t$ the maximum of all $t_v$. Since the binomials in $x_1$ commute, for $w$ we obtain $w(x_1-\alpha_{1,s})\cdots(x_1-\alpha_{1,t})\in I$. Since $M_i = (I + \Delta_1\Delta_2\cdots \Delta_i)/I$, the proof is complete.\bx
\end{enumerate}

\begin{lemma}\label{l994} Let every element of the field $\Phi$ occur in the sequence $\alpha(n)$ at most finitely many times, $L$ a module over $R$ and any $w \in L_i = L\Delta_1\Delta_2\cdots \Delta_i$ is annihilated by a quasi-monomial $(x_1-\alpha(s))\cdots (x_1-\alpha(t))$ and $i\le s<t$. Then $\bigcap  L_i =0$.
\end{lemma}
 
\pp Assume $w \in \bigcap L_i$. Then $w$ is annihilated by an infinite sequence of monomials $(x_1-\alpha(s_i))\cdots(x_1-\alpha(t_i))$, where the sequence $s_i$ is unbounded. But in their totality, these monomials in one variable are coprime, by the condition on the sequence $\alpha(i)$. Hence $w=0$, as needed for the proof. \bx

Before we formulate the main result of this section, we recall two definitions.

A module $M$ over an algebra $R$ is called \emph{triangular} if $M$ has a linear basis $\{e_i\:|\:i=1,2,\ld\}$, such that $e_iu$ is an element of the linear span of $\{ e_i,e_{i+1},\ld\}$, for any $u\in R$ and any $i=1,2,\ld$.

A right $R$-module $Q$ is called a \emph{section} of a right $R$-module $M$ if $Q\cong N/P$ where $N,P$ are submodules of $M$ with $P\subset N$. 

\begin{theorem}\label{t991} Let $\Phi$ be an infinite field. Then there exists a cyclic module $M$ over $R=\far{r}$, $r>1$, enjoying the following properties.
  \begin{enumerate}
  \item[\emph{(a)}] $M$ has maximal growth,
\item[\emph{(b)}] $M$ is triangular,
\item[\emph{(c)}] For any factor-module $L$ we have $\cap L_i=\{ 0\}$ where $L_i=L\Delta_1\cdots\Delta_i$,
\item[\emph{(d)}] Any factor-module $L$ of $M$ (including $M$ itself!) is residually finite,
\item[\emph{(e)}] Any simple section $Q$ of $M$ is one-dimensional.
 \end{enumerate}
 \end{theorem}
 
\pp We choose $M=R/I$, as in Lemma \ref{l993}. Note that the choice of $I$ depends on the choice of infinite sequences $\alpha_{1i},\ld,\alpha_{ri}$. For our purposes, we need to assume that none of these sequences has infinite repetitions of the same number. Then we have the following.
\begin{enumerate}
\item[(a)] The growth of $M$ is maximal by Claim (a) in Lemma \ref{l993}.
\medskip
\item[(b)] It follows by (\ref{(*)}) that the basis of the ideal $I$ consists of all $e_v$, such that $v$ begins with a leading term of a quasi-monomial $w_i$. By Claim (a) of Lemma \ref{l991} the images of all the remaining $e_v$ form a basis in $M=R/I$.  If we order the set of these remaining  $e_v$ by ShortLex, then the action of the generators $x_j$ with respect to this basis  is triangular, as readily seen from (\ref{(*)}).
\medskip
\item[(c)] Let $L=M/P$ be a factor-module of $M$. Since the  condition of Lemma \ref{l994} holds for $M_i = M\Delta_1\cdots \Delta_i$, it also holds for $L_i = L\Delta_1\cdots \Delta_i$. Then it follows that $\cap L_i = \{ 0\}$, as claimed.
\medskip
\item[(d)] We continue the argument of the previous claim (c). It remains to explain why all $L/L_i$ are finite-dimensional. If we set $R_i=\Delta_1\cdots \Delta_i$ then $L/L_i$ is a factor-module of $R/R_i$ and so it is sufficient to prove the finite-dimensionality of this latter. However, this easily follows since by Lemma \ref{l991} $R/R_i$ has a basis of some $e_v$ where $\deg v <i$.
 \medskip
 \item[(e)] Let $Q=N/P$ be a simple section. Here $N$ and $P$ are submodules of $M$, $P\subset N$. Let $Q_i=L_i\cap Q$, where, as before, $L_i = L\Delta_1\cdots \Delta_i$. It follows by (c) that $\bigcap Q_i =\{ 0\}$. Applying (d), we find that $Q$ must be finite-dimensional. Since all $Q_i$ are submodules in a simple module $Q$, there is $s\le i$, such that $Q_{s-1}=Q$, but $Q_s= Q_{s-1}\Delta_s = Q\Delta_s =\{ 0\}$, hence $v(x_j - \alpha_{j,s})=0$ for any $j$ and any $v\in Q$. Thus, the one-dimensional space $\Phi v$ is an $R$-module, and so $Q=\Phi v$, and the proof is complete.\bx
\end{enumerate}

\begin{problem}\label{op5}
Is the analogue of Theorem \ref{t991} true in the case of modules over free group algebras?
\end{problem}

The statement of the next open problem reminds the following still open problem due to M. I. Kargapolov \cite[Problem 1.31]{KT}: ``Is a residually finite group with the maximum condition a finite extension of a polycyclic group?''

\begin{problem}\label{op2} Do there exist residually finite Noetherian modules of maximal growth?
\end{problem}

%Similar examples can be produced also in the case of $kF$, the group algebra of a free group $F$ of rank $r>1$. The polynomials $(x-\alpha_1),\ldots,(x-\alpha_m)$ can be replaced by irreducible polynomials of arbitrary degrees, which will allow one to cover also the case of finite fields.

%%%%%%%%%%%%%%%%%%%%%%%%%%%%%%%%%%%%%%%%%
\section{Group actions of maximal growth}\label{sGAMG}

Let $F_r=\fgax$ where $\# X=r$ be the free group of rank $r>1$. In this section we study acts over $F_r$, more commonly called $F_r$-sets.

\subsection{Cayley graphs of $F_r$-sets}\label{ss1} 
We start with a free group  $F$ with a symmetric basis $A$, that is, a union $\{ a_1,\ld,a_r\}\cup\{ a_1^{-1},\ld,a_r^{-1}\}$, where $\{ a_1,\ld,a_r\}$ is a free basis of $F$. Suppose $F$ acts on a set $S$. We introduce a directed Cayley graph of this action, with labeling $\mathcal{G}(S)=(V,E,\mathrm{Lab})$ as follows. For vertices, we set $V=S$. The edges appear as follows. Given $s\in S$ and $a\in A$, there is unique edge $e\in E$ whose \emph{source} $e_{-}$, respectively, \emph{target} $e_{+}$ equals $s$, respectively, $s\cir a$, and $\lb{e}=a$.

So $e_{-}=f_{-}$ and $\lb{e}=\lb{f}$ always imply $e=f$. Two labeled edges $e$ and $f$ are called \emph{inverses} of each other, $f=e^{-1}$, if $e_{+}=f_{-}$, $f_{+}=e_{-}$ and $\lb{f}=\lb{e}^{-1}$.  If $e_{i}, e_{i+1}$ are two consecutive edges in a path $p=e_1\cdots e_n$ then $(e_{i+1})_{-}=(e_{i})_{+}=(e_{i}^{-1})_{-}$. If also $\lb{e_{i+1}}=\lb{e_{i}}^{-1}=\lb{e_i^{-1}}$ then by definition $e_{i+1}=e_{i}^{-1}$. As a result, if a path $p=e_1\cdots e_n$ is \emph{reduced} (that is, has no subpaths $ee^{-1}$) if and only if its label $\lb{p}=\lb{e_1}\cd\lb{e_n}$ is a reduced word in the free group $F$. 
%The degree of each vertex (that is, the number of outgoing edges) equals $2r$.
 
Given a Cayley graph $\mathcal{G}(S)$ as just above, the star of each vertex $v\in \mathcal{G}$ (we denote it by $\mathrm{star}(v)$) has exactly $2r$ outgoing edges labeled by all $a\in A$. We will call such a star \emph{standard.} If all stars of a graph $\mathcal{G}$ are standard
($r$-standard), then their labels define the action of the free group $F$ on $\mathcal{G}$. It is obvious that $\mathcal{G}$ is connected if and only if the action of $F$ on $S$ is transitive. 
 
From now on we assume that  $\cg$ is connected. Let us distinguish one vertex $o$ of graph $\mathcal{G}(S)$ and denote by $H$ the stabilizer of this vertex under the action of $F$ on $\mathcal{G}$. Then the set $V$ of vertices of $\mathcal{G}(S)$ is in one-to-one correspondence with the set of right cosets of $H:   o\cir g \leftrightarrow Hg$, that is, one can identify $\mathcal{G}(S)$ with the graph of right cosets of $H$ in which there is an edge with label $a$ from $Hg$ to $Hg'$ as soon as that $Hga = Hg'$.
 
Conversely, for each subgroup $H$ of $F$ the graph of right cosets $\mathcal{G}(F/H)$ is the Cayley graph for the right action of $F$ on $F/H$ with $H$ as the stabilizer of the coset $H$.
 
If $H\subset H_1$, then there is a well-defined map $\vp:F/H\ra F/H_1$ such that $\vp(Hg) = H_1g$. Obviously, $\vp$ commutes with the natural action of $F$ on $F/H$ and $F/H_1$. Hence $\vp$ induces a surjective label-preserving morphism of graphs $\cg(F/H)\ra\mathcal{G}(F/H_1)$. (Such  morphisms of coset graphs have been introduced and studied by J. Stallings
\cite{JS1} in terms of foldings and covering maps. But for our purposes it is more useful to remember that this morphism is a morphism of $F_r$-sets, that is, if  $S = F/H$ and $S_1 = F_1/H_1$ then we have a surjection $\vp: S\ra S_1$, such that $\vp(s\cir u) = \vp(s)\cir u$ for any $s\in S$ and any $u\in F$.)

%%%%%%%%%%%%%%%%%%%%%%%%%%%%%%%%%%%%%%
Now let $B(n)=\mb{o}{n}$, $n=0,1,\ld$, stand for the \emph{ball} of radius $n$ in $\cg$ with center $o$; each such ball consists of all vertices in the distance at most $n$ from $o$. Similarly one defines the \emph{sphere} $S_n=\ms{o}{n}$. We recall that the growth function $g(n)=g_{o,E}(n)$ of an $F_r$-set $E$ or of the action of $F_r$ on $E$ or else on $\mathcal{G}$ is given by  $g(n)=\# B(n)$. 
										
\begin{lemma}\label{Lemma B} Let $\mathcal{G}$ be the Cayley graph of the action of a free group $F$ of rank $r$ on the set $F/H$. Let us assume that $\# B(n)  \ge  c(2r-1)^n$  for some real $c>0$ and integral $n>0$. Then $\# B(m) \ge \dfrac{2r-2}{2r-1}  c(2r-1)^m$ for all $m$ such that $0<m < n$.
\end{lemma}
										
\pp Proving by contradiction, let us assume that $\# B(m) < \dfrac{2r-2}{2r-1} c (2r-1)^m$ for some $m$ as above. We connect each vertex in $B(n) \setminus B(m)$ with $o$ by a geodesic path $p$. Then $p=p'p''$, where $p'$ has length $\le n-m$ and connects $v$ with a vertex $u$ on the sphere $S_m$. Here $p''$ depends only on the choice of $u$. One can uniquely recover $v$ by $u$ and by the reduced label of $p'$. Now the last letter of the label of $p'$ differs from the inverse of the first letter of the label of $p''$. It follows that the number of possible reduced labels for $p'$ is at most $(2r-1) +\cd+(2r-1)^{n-m}  < \dfrac{2r-1}{2r-2}(2r-1)^{n-m}$. Hence  
\begin{eqnarray*}
\# B(n)&\le&(\# S_m)\dfrac{2r-1}{2r-2}\left(2r-1\right)^{n-m} +\# B(m)\\
&\le&  (\# B(m))\left(1+\dfrac{2r-1}{2r-2}(2r-1)^{n-m}\right)\\
&< &\frac{2r-2}{2r-1} c(2r-1)^m\left(1+\dfrac{2r-1}{2r-2}(2r-1)^{n-m}\right)\\
&\le&\frac{2r-2}{2r-1}  c \left(1+ \dfrac{1}{2r-2}\right)(2r-1)^n.
\end{eqnarray*}
It follows that $\# B(n)<c(2r-1)^n$, a contradiction.\bx
%%%%%%%%%%%%%%%%%%%%%%%%%%%%%%%%%%%%%

\subsection{Core of the Cayley graph}\label{ss2} 
In what follows we identify $\mathcal{G}(S)$ with $\mathcal{G}(F/H)$ and write simply $\mathcal{G}$ (or $\mathcal{G}_1$ in the case where $H_1$ replaces $H$). For any path $p=e_1\ld e_n$ in $\mathcal{G}$ the word $\Lab{p} = \Lab{e_1}$ $\cdots$ $\Lab{e_n}$ is an element of $F$. Each element of $H$ can be read as the label of a unique reduced loop $p$ such that $p_-=p_+=o$. Conversely, by definition of the stabilizer of the vertex $o$, the label of each loop as above is an element of $H$.

We denote by $\mathcal{C}$ the minimal subgraph of the graph $\mathcal{G}$ which contains the vertex $o$ and all reduced loops originating in $o$. One can write $\mathcal{C}=\mathcal{C}(\mathcal{G},o)=\mathcal{C}(H)$, the latter because the choice of a subgroup $H \le F$ is equivalent to the choice of a connected graph $\mathcal{G}$ with standard stars of vertices and a distinguished vertex $o\in  \mathcal{G}$.  Following the above mentioned paper \cite{JS1} we will call $\mathcal{C}(\mathcal{G},o)$ the \emph{core of $\mathcal{G}$}. 

Given a subgroup  $H$ of $F$, let $p_1,\ld,p_s,\ld$ be reduced loops corresponding to the reduced forms of some generators $h_1,\ld, h_s,\ld$ of $H$. Then any element of $H$ can be written as a loop $p$ resulting after all possible cancellations from a product $p(1)\cdots p(t)$ where $p(i)\in \{p_1^{\pm 1},\ld,p_s^{\pm 1},\ld\}$. Therefore, the subgraph $\mathcal{C}$ does not contain any edges other than those in the paths $p_1,\ldots,p_s,\ld$ In particular, $\mathcal{C}$ is finite for a finitely generated subgroup $H$.
 
%\begin{picture}\label{pict1}
%\end{picture}

Since all edges of reduced loops with origin $o$ are in $\mathcal{C}$ and $\mathcal{G}$ is a connected graph, one obtains $\mathcal{G}$ by attaching labeled trees $\mathcal{T}_1, \mathcal{T}_2 ,\ld$ to $\mathcal{C}$ in such a way that one vertex of each tree is attached to one vertex of $\mathcal{C}$.  In each of these trees all stars of its vertices are standard, except the root vertex $o(\mathcal{T}_i)$, which is also a vertex of $\mathcal{C}$. Still, all edges of $\mathcal{T}_i$ whose source is $o(\mathcal{T}_i)$ must have different labels (but their number is at most $2r$). We will call the star of a vertex $v$ (and $v$ itself!) \emph{regular} if the labels of all edges in $\mathrm{star}(v)$ are pairwise different. 

The trees $\mathcal{T}_1, \mathcal{T}_2,\ld$ form a ``forest'' $\mathcal{F}$.
 
Since all vertices in the Cayley graph $\mathcal{G}$ of the action of the free group $F$ of rank $r$ have degree $2r$, the trees will be attached precisely to those vertices $v$ of $\mathcal{C}$ whose degree $d_C(v)$ in $\mathcal{C}$ is less than $2r$. In addition, if $v\in \mathcal{C}$ is identified with $o(\mathcal{T}_i)$ then $d_C(v)+d_{\mathcal{T}_i}(o(\mathcal{T}_i))=2r$, and the labeling of the star of $o(\mathcal{T}_i)$ in $\mathcal{T}_i$ complements the labeling of the star of $v$ in $\mathcal{C}$.

It follows from that the above that the labeled graph $\mathcal{G}$ with a distinguished vertex $o$ can be uniquely (up to isomorphism) recovered by the subgraph $\mathcal{C}$  with the distinguished vertex $o$. From the definition of $\mathcal{C}$ it is also clear how $\mathcal{C}$ determines the subgroup $H$. It follows that every subgroup $H$ in $F$ can be determined by a connected labeled graph $\mathcal{C}$ with the distinguished vertex $o$ such that the stars of all vertices in $\mathcal{C}$ are regular and have degree at least 2, except possibly for $o$, which may have degree 1.
 
If $H$ is a finitely generated subgroup then the number of vertices in $\mathcal{G}$ is finite and so the ``forest'' $\mathcal{F}$ is nonempty if and only if $S$ is infinite, that is, when $[F:H]$ is infinite.
										
%%%%%%%%%%%%%%%%%%%%%%%%%%%%%%%%%%%%%%%%%%%										
										
%%%%%%%%%%%%%%%%%%%%%%%%%%%%%%%%%%%%%%%%%%%
\subsection{One lemma about graphs with bounded degrees of vertices}\label{ss}
We consider graphs each of which has a distinguished vertex. We will associate with each such graph a numerical invariant and show that it is bounded from above by a constant that does not depend on the graph. At this point we do not assume that the graphs are endowed by labeling or even that they are directed. Given such a graph $\Gamma$ and a vertex $v\in V(\Gamma)$, we will denote by $|v|$ the combinatorial distance from $o$ to $v$.
										 
\begin{lemma}\label{Lemma A}  Let $\Gamma$, as above, be an arbitrary connected graph with a distinguished vertex $o$. Fix an arbitrary integer $m >1$. If we have $d_{\Gamma}(v)\le m$ for the degree of every vertex $v\in V(\Gamma)$ then 
\begin{equation}\label{(**)}
\sum_{v\in V(\G)}  (m - d_{\Gamma}(v) )  (m-1)^{-|v|}  \le  m.
\end{equation}
\end{lemma}
\pp Notice that one can write $\G$ as the union of an ascending chain of subgraphs $\G_0 \subset \G_1\subset\ld$ so that 
\begin{enumerate}
\item[(1)] $\G_0 =\{o\}$;
\item[(2)] each subgraph can be obtained from the previous one by adding one edge and at most one vertex; 
\item[(3)] if $v \in \G_i$, then the distance from $v$ to $o$ in $\G_i$ is the same as in $\G$.
\end{enumerate}
Indeed, suppose $i>0$ and the graphs $\G_k$ with $k<i$ have been selected. Then, if possible, adjoin an edge $e\in E(\G)\setminus E(\G_{i-1})$ if $e$ connects two vertices of $\G_{i-1}$ and call the resulting graph $\G_i$ (type (a) transformation). Otherwise, enlarge  $V(\G_{i-1})$ by adding a new vertex $v\in V(\G)$ with minimal distance to $o$ and $E(\G_{i-1})$ by the last edge of a geodesic path from $o$ to $v$; call the resulting graph $\G_i$ (type (b) transformation). If none of (a) or (b) applies, then $\G_i=\G_{i-1}=\G$.  It is easy to check that the chain thus constructed  satisfies all conditions (1) to (3) and that $\G$ is the union of all $\G_i$.
										
It is sufficient to prove (\ref{(**)}) for each $\G_i$ since every finite portion of the left hand side is majorated by the respective sum composed for a subgraph $\G_i$ for $i$ sufficiently big. So we will use induction by $i=0,1,2,\ld$ to prove (\ref{(**)}) for any $\G_i$ in place of $\G$.
										
In the case of $\G_0$ the left hand side of (\ref{(**)}) equals $m(m-1)^0\le m$, as needed. Assume the inequality is true for $\G_{i-1}$ for some $i>0$. Let us consider the two types of transformations used while switching from $\G_{i-1}$ to $\G_i$. If we apply (a) then all the vertices and their distances from $o$ remain the same and at the same time the coefficients $m -\deg(v)$ do not grow. Hence (\ref{(**)}) remains true also for $\G_i$.  If we apply (b), then we have to add to $V(\G_{i-1})$ one new vertex $v$ and an edge $e$ with endpoints $u \in V(\G_{i-1})$ and $v$. We have then that $d_{\G_i}(v)=1$ in $\G_i$, and $d_{\G_i}(u)=d_{\G_{i-1}}(u)+1$. This adds to the left side of (\ref{(**)}) a summand $(m-1) (m-1)^{-|v|}$ corresponding to $v$ and subtracts from the summand corresponding to $u$ the value $(m-1)^{-|u|} = (m-1)^{-|v|+1}$. As a result, while applying (b), the left hand side of (\ref{(**)}) remains unchanged. By induction, our claim is proven.\bx  

\subsection{Deficit of graph}\label{ssDG}										
										
We return to the discussion of Subsection \ref{ss2} and consider the core $\cc$ of a Cayley graph $\cg(F/H)$, where $H$ is a subgroup of a free group $F$ of rank $r$. For any vertex of a graph $\cc$ we define its \emph{deficit} by setting $\df_{\cc} (v) = 2r-d_{\cc}(v)$. This is the number of edges outgoing from $v$ in $\cg$, which are not in $\cc$ but rather belong to the forest $\mathcal{F}$. We introduce also the deficit of $\cc$ by setting 
\bea
\df(C) = \sum _{v \in C} \df_C(v) (2r-1)^{-|v|}.
\eqa
By Lemma \ref{Lemma A} this value is finite and is at most $2r$. To apply Lemma \ref{Lemma A} to $\cc$ we have to remove all edges labeled by inverses of the generators of $F$ and  erase all labels and directions (arrows) on the remaining edges. Notice that $\df(\cc)$ is a measure of the ``infiniteness'' of the index of $H$ in the following sense. If $H$ is a finitely generated subgroup of $F$ then $\df(\cc)=0$ if and only if $H$ is of finite index in $F$, or, in other words, the forest $\mathcal{F}$ is missing.
										
\begin{lemma}\label{Lemma C}
\begin{enumerate}  
\item[\emph{(1)}] If a subgroup $H$ is finitely generated and $g(n)$ is the growth function of the action of $F$ on $F/H$ then there is a function $f(n)$ with $|f(n)|$ bounded such that
$g(n) = \dfrac{\df(C)}{2r-2} (2r-1)^n +f(n)$.

\item[\emph{(2)}] In any case,
\bea
g(n) \ge \frac{\df(C)}{2r-1}  (2r-1)^n
\eqa
for all $n > 0$.
\end{enumerate}
\end{lemma}
										
\pp  To start with, we find the number $h(n)$ of vertices in the ball $B(n)$, outside $\mathcal{C}$. Each vertex is in one of the trees $\mathcal{T}\in\mathcal{F}$. So we first need to compute the number of vertices in $\mathcal{T} \cap B(n)$ which are not in $\cc$.
										 
Let $\mathcal{T}$ be attached to $\mathcal{C}$ at the vertex $v$ whose distance from $o$ is $m$. Then the distance from $v' \in \mathcal{T}$  to  $o$  equals  $m + \dist_{\mathcal{T}}(v,v')$. Since the degree of every vertex in $\mathcal{T}$, except $v$, is $2r$, there are $\Def{\mathcal{C}}{v}$ vertices in $\mathcal{T}$ whose distance is 1 from $v$, hence $m+1$ from  $o$. Next, there are $\Def{\mathcal{C}}{v} (2r-1)$ vertices in the distance 2 from $v$, hence $m+2$ from  $o$. Finally, there are $\Def{\mathcal{C}}{v}(2r-1)^{n-m-1}$ vertices in the distance $n-m$ from $v$, hence $n$ from $o$.
As a result,  if $n>m$, then $B(n)\cap \mathcal{T}$ has $\Def{\mathcal{C}}{v} \dfrac{(2r-1)^{n-m}-1}{2(r-1)}$ vertices outside $\cc$. Otherwise the set of such vertices is empty.
										
Performing summation over all vertices of $\cc$ such that $|v|< n$, we will obtain: 
\begin{eqnarray}\label{(***)}  &&h(n) = \frac{1}{2r-2} \sum_{v\in \cc, |v|<n} \df_{\cc}(v)((2r-1)^{n-|v|})-1)\nonumber\\
&&=\frac{1}{2r-2}  \left((2r-1)^n \sum _{v\in \cc, |v|<n} \df_{\cc} (v) (2r-1)^{-|v|} - \sum_{v\in \cc , |v|<n} \df_C (v) \right)
\end{eqnarray}
										
Notice that $\df_{\cc}(v) \le 2r-2$ for all vertices except $o$, in which case $\df_{\cc}(o) \le 2r-1$. Therefore, 
\bea
\sum_{v\in \cc , |v|<n} \frac{\df_C (v)}{2r-2}\le \#(B(n) \cap\cc) + \frac{1}{2r-1}.
\eqa
Using this and (\ref{(***)}), we obtain
\begin{eqnarray}\label{(****)}
&&g(n) =h(n) + \# (B(n)\cap \cc)\nonumber\\&& \ge  \frac{1}{2r-2} \left((2r-1)^n \sum _{v\in\cc, |v|<n} \df_{\cc} (v) (2r-2)^{-|v|} -\frac{1}{2r-1}\right)
\end{eqnarray}
										
If $\cc$ is finite, then the first sum on the right side of (\ref{(****)}) is $\dfrac {\df(C)}{2r-2} (2r-1)^n$, for all $n$ greater than the diameter of $\cc$.  Thus the proof of Claim (1) is complete.
										
To prove Claim (2), we set 
\bea
D_n =   \sum _{v\in\mathcal{C},|v|<n} \df_{\cc} (v) (2r-1)^{-|v|}.
\eqa
Then $0\le D_n \le 2r$, by Lemma \ref{Lemma A}, and for any $\ve >0$ we have $0\le \df(\cc) - D_n <\ve$, if only $n > N(\ve)$. The right hand side in (\ref{(****)}), for such large $n$ is greater than $\dfrac{1}{2r-2}(\df(\cc) - 2\ve) (2r-1)^n$, following since the subtracted term in (\ref{(****)}) is bounded by a constant. By Lemma \ref{Lemma B},  $g(n) \ge \dfrac {1}{2r-1} (\df (\cc)-2\ve)(2r-1)^n$, where now $n$ is any positive integer. Since $\ve$ can be chosen arbitrarily small, our Claim follows.\bx

\subsection{Elementary graphs}\label{ssEG}

%%%%%%%%%%%%%%%%%%%%%%
Before we formulate our next lemma, we briefly recall the construction of the Schreier basis for a subgroup $H$ of a free group. Let $\cg=\cg(H)$ be the Cayley graph of the $F$-set $F/H$. Using Zorn's Lemma, one can select in $\cg$ a maximal subtree $\mathcal{T}$. For any two vertices $u,v\in\mathcal{G}$ there is a unique reduced path in $\mathcal{T}$ that goes from $u$ to $v$. For some considerations it is important that $\mathcal{T}$ can be chosen so that the length of the path from $o$ to $v$ within $\mathcal{T}$ is the shortest among all the lengths of the paths from $o$ to $v$ within $\cg$. Such maximal subtrees are called \emph{geodesic}. With $\mathcal{T}$ fixed, the set of edges $E(\mathcal{G})$ splits into two subsets: the edges in $E(\mathcal{T})$, called \emph{tree} edges, and the edges in $E(\mathcal{G})\setminus E(\mathcal{T})$, called \emph{non-tree} edges. Suppose $e$ is a non-tree edge. Let $p$ be the reduced path on $\ct$ from $H$ to $e_{-}$ while $q$ the reduced path on $\ct$ from $H$ to $e_{+}$. Suppose $v=\lb{p}$, $a=\lb{e}$ and $w=\lb{q}$; then  $v a w^{-1}$ is called a \emph{Schreier generator} for $H$. Since $e^{-1}$ is also a non-tree edge, the Schreier generator defined by $e^{-1}$ is the inverse of the Schreier generator defined by $e$. The collection $B=B(\mathcal{T})$ of all Schreier generators is known to be a symmetric basis of the free group $H$ \cite{LS}.
										
\begin{lemma}\label{lff} Let $H \subset H_1$ be subgroups of a free group $F$ such that the restriction of the induced morphism of the graphs of actions $\cg\ra\cg_1$ to the core $\cc$ is an injective embedding of $\cc$ to $\cc_1$. Then $H$ is a free factor of $H_1$.
\end{lemma}
\pp We will use the construction of the Schreier basis for $H$ described just before the statement of this lemma. 
										
Notice that a maximal subtree of $\cg$ consists of a maximal subtree of the core $\cc$ and the forest $\mathcal{F}$. From our hypotheses using Zorn's Lemma it follows that a maximal subtree of $\cc$ can be extended to a maximal subtree of $\cc_1$. As a result, the set of all non-tree edges of $\cc$ (and $\cg$) can be included in the set of all non-tree edges of  $\cc_1$ (and $\cg$). By the construction of the Schreier system of free generators for $H$ and $H_1$ it follows that some free basis of $H$ can be included in a free basis of $H_1$, proving our lemma.\bx

%%%%%%%%%%%%%%%%%%%%%%%
\medskip									
Let us call a graph  with labeling \emph{elementary} if it has one of the following forms. 
\begin{enumerate}
\item A simple cycle with one distinguished vertex;
\item A cycle and a simple arc (the \emph{cycle with leg} of the graph) glued together by one vertex, with a distinguished vertex of degree 1;
\item A simple arc with two endpoints distinguished.
\end{enumerate}
The restriction on the labeling is as follows:  all stars in the graph are regular in the sense that in every star different edges have different labels.
										
In each elementary graph there is a reduced \emph{distinguished path} $q$ as follows. In the case (1), it starts at the distinguished vertex and goes around the cycle until it reaches the distinguished vertex, where it terminates. In the case (2), it starts at the distinguished point of the leg, goes in the direction of the cycle, then goes around the cycle until it reaches the leg and then goes along the leg toward the distinguished vertex, where it terminates. In the case (3), it starts at one distinguished point and goes toward the second distinguished point, where it terminates. 
										 
By \emph{attaching an elementary graph} to a graph $\Gamma$ with labeling and with regular stars we mean identifying a distinguished vertex (or two vertices) with one vertex of $\Gamma$ (or with any two different vertices of $\Gamma$, respectively), so that in the resulting graph all stars remain regular.
\begin{lemma}\label{Lemma D}  Let $H$ and $H_1$ be finitely generated subgroups of a free group $F$ such that the core $\mathcal{C}_1= \mathcal{C}(H_1)$ is obtained from the core $\mathcal{C}=\mathcal{C}(H)$ by attaching an elementary graph with distinguished path $q$ of length $l$. Then $H_1=H*\langle g\rangle$, where $g$ is the label of an arbitrary reduced loop with source $o$ which contains each edge from $q$ exactly once. For the deficits of the cores $\cc$ and $\cc_1$ one has
\bea
0 \le  \df(\cc)-\df(\cc_1)  \le (2r-1)^{2-\frac{l}{2}}.
\eqa
\end{lemma}
										
\pp Using the proof of Lemma \ref{lff}, we find a free basis $\{ g_1, h_1,h_2,\ld\}$ in $H_1$, where $\{ h_1, h_2,\ld\}$ is a free basis of $H$ and $g_1$ is the label of a path containing $q$, which includes each edge of $q$ exactly once. In this case $g = u_1g_1u_2$ for some $u_1, u_2 \in H$, hence $H_1 = H*\langle g\rangle$, as claimed.
										
To simplify the proof of the main claim, we will not be assuming that the stars in $\cc$ and $\cc_1$ are regular (see Case (2)) and even allow the multiplicity of the stars being greater than $2r$ (the same as their deficits being negative!).
										 
\textbf{Case (1).}  Suppose we attach a simple arc or a cycle, so that $|v_0|=m$ and $|v_l|=k$, where $v_0=q_-$ and $v_l=q_+$. When we switch from $\cc$ to $\cc_1$, the contribution of these vertices to the deficits of respective graphs goes down by $(2r-1)^{-m}$ and $(2r-1)^{-k}$, respectively. The total decrease from these two vertices is $(2r-1)^{-m} +(2r-1)^{-k}$, even in the case where $v_0=v_l$.  At the same time, there are $l-1$ vertices with deficit $2r-2$ in $\cc_1 \setminus \cc$, which provide a positive contribution when computing $\df(\cc_1)$.
										
Let $v_1,\ld, v_{l-1}$ be the vertices in $\cc_1 \setminus \cc$. It is easy to observe that $|v_i|=\min (m+i, k+(l-i))$ so that the contribution of $v_i$ to $\df(\cc_1)$ equals  \bea(2r-2)(2r-1)^{-\min (m+i, k+(l-i))}.\eqa  The sum of all these contributions to $\df(\cc_1)$ is less than $(2r-2)\sum_{i =1}^{\infty} ((2r-1)^{-m-i} +(2r-1)^{-k-i}) =(2r-1)^{-m} +(2r-1)^{-k}$, which yields the desired inequality $\df(\cc) - \df(\cc_1) > 0$, in this case.
										
For the proof of the second inequality, we may assume $0\le  m-k\le l$. Then for $s =
[(l+k-m)/2]$, each vertex $v_j\in\{ v_1,\ld, v_{s-1}\}$ is at the distance $j+m$ from $o$, and at the same time, each $v_j\in\{ v_s,...v_{l-1}\}$ is at the distance $k+(l-j)$. Therefore, the total contribution of all these vertices to $\df(\cc_1)$ is at least 
\begin{eqnarray*}
&&(2r-2)\left(\sum_{i =1}^{s-1} (2r-1)^{-m-i} +\sum_{i=s}^{l-1}(2r-1)^{i-l-k}\right)\\
&=&(2r-2)\left(\sum_{i =1}^{s-1} (2r-1)^{-m-i}+ \sum_{j=k+1}^{k+l-s}(2r-1)^{-j}\right).
\end{eqnarray*}
If we subtract this value from 
\bea(2r-1)^{-m} +(2r-1)^{-k} =
(2r-2)\left(\sum_{i=1}^{\infty} (2r-1)^{-m-i}+\sum_{j=k+1 }^{\infty} (2r-1)^{-j}\right),
\eqa
then we obtain
\begin{eqnarray*}
&&(2r-2)\left(\sum_{i=s}^{\infty} (2r-1)^{-m-i}+\sum_{j=k+l-s+1 }^{\infty} (2r-1)^{-j}\right)\\
& <&(2r-1)^ {-m-s+1} +(2r-1)^{-k-l+s}.
\end{eqnarray*}
Since $s = [(l+k-m)/2]$, this sum is less than $(2r-1)^{(-l/2) + 2}$, proving the second inequality, in this case.

\begin{figure}[ht]
	\centering
  \includegraphics[height=45mm]{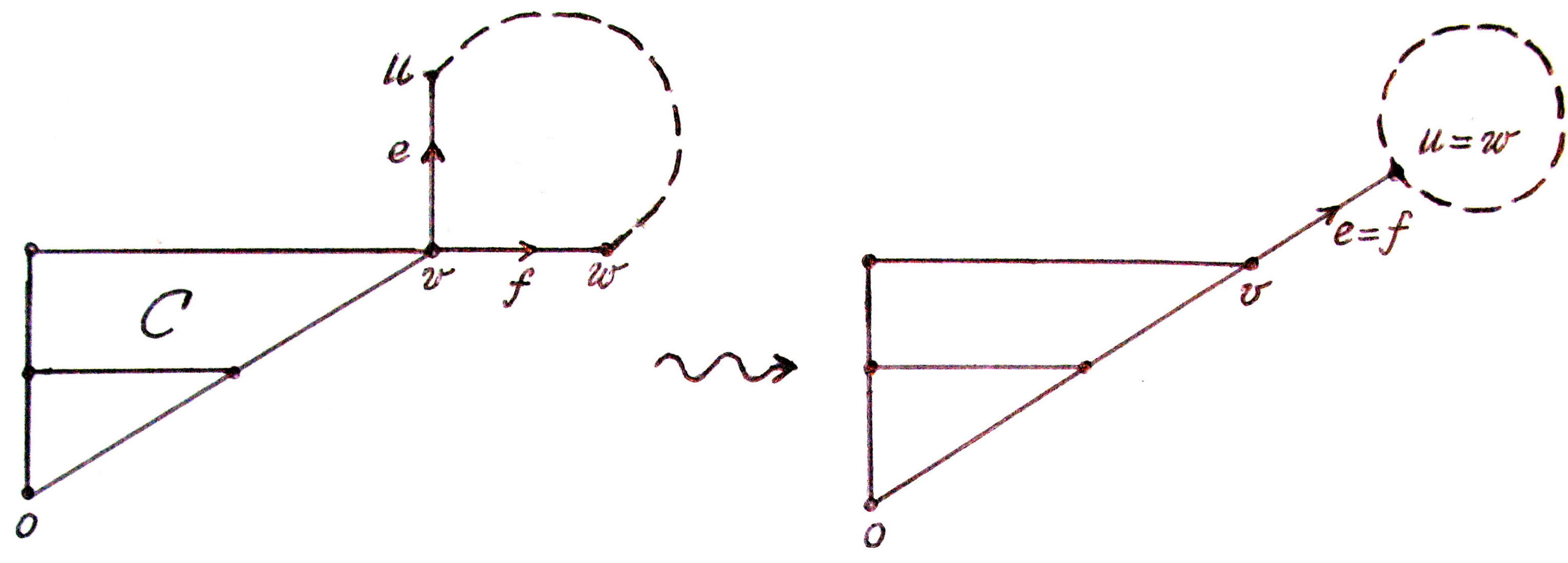}
 \caption{Transition from loop with leg}
\end{figure}
 
\textbf{Case (2).}  Suppose we attach a ``cycle with leg''. This can be done in two stages: first, we attach a cycle and second, we identify, one by one, several edges of the cycle to form a leg. If we prove that the appropriate operations of the second type do not change the deficit of the graph, then this case will be completely reduced to the previous one.
										
We start with by attaching at one vertex $v\in\cc$ of a cycle whose label is the same as the label of the distinguished path $q$. Let us see what happens if we identify two edges $e$ and $f$ of a cycle,  both incident to $v$. Let $d$ be the distance of $v$ from $o$. After the identification, the contribution of $v$ to the deficit of the graph increases by $(2r-1)^{-d}$. Now we need to measure the change of the deficit produced by the endpoints $u$ and $w$ of $e$ and $f$, respectively. Before the identification, the contribution from these two vertices of degree 2 was $2(2r-2)(2r-1)^{-d-1}$. After the identification, we obtain one vertex of degree 3 contributing the value of $(2r-3)(2r-1)^{-d-1}$ into the common deficit. The decrease obtained is $(2(2r-2)- (2r-3))(2r-1)^{-d-1} = (2r-1)^{-d}$, which matches the increase produced by $v$. Since the contribution of other vertices remains unchanged, the proof is now complete also in Case (2).\bx

\subsection{Adding a power of element}\label{ssAPE}

Given an element $g$ of a group $G$, we denote by $\langle g\rangle$ the cyclic subgroup of $G$ generated by $g$. Our main result in this subsection is the following.										

\begin{theorem}\label{t001n} Let $H$ be a finitely generated subgroup in a free group $F$ of finite rank $r>1$ and $g\in F$ ($g\neq e$) such that $H\cap\langle g\rangle$ is trivial. Then for any $\ve\in(0,\df(\cc))$ there exists a natural $n =n(g,H)$ and a subgroup $H_1$ of infinite index in $F$ containing $g^n$ such that $H_1 = H *\langle g^n\rangle$ and $\df(\cc)-\df(\cc_1)\le\ve$, where $\cc$ and $\cc_1$ are the cores of the Cayley graphs $\cg(F/H)$ and $\cg(F/H_1)$, respectively.
\end{theorem}
 
\pp  Let us write $g$ as the reduced product $uwu^{-1}$ where $w$ is a cyclically reduced word. In this case all integral powers $g^n = uw^n u^{-1}$ will also be reduced. Since all cosets $Hg^i$ are pairwise different for different values of $i$, all cosets $Huw^ i$ are pairwise different, as well.

Consider the core $\mathcal{C} = \mathcal{C}(H)$ in the Cayley graph $\mathcal{G}$ of the action of $F$ on $F/H$. Since $H$ is finitely generated, this graph is finite and for any $l>0$ one can choose vertices $v_i, v_j$ of the $\mathcal{G}$ corresponding to some $Huw^i$ and $Huw^j$, with $i>0, j<0$, in the distance greater than $l$ from $\mathcal{C}$.

Let $p_i, p_j$ be the shortest paths from the distinguished vertex $o$ to $v_ i$ and $v_j$. They look like $p_k=p'_kp''_k$, $k=i,j$ , where $p'_k$ is a subpath in $\mathcal{C}$ while $p''_k$ is a subpath of the forest of the graph $\mathcal{G}$ and $|p''_k| > l/2$ for any given $l$. 

We are going to construct a new graph $\mathcal{C}_1$ by adding paths  $p''_i$ and $p''_j$ to $\mathcal{C}$ and by identifying their target vertices.  In the graph thus obtained, the label on the loop $s=p^{''}_i(p^{''}_j)^{-1}$ with origin $o$  reads as the word $uw^i w^{-j}u^{-1}$, which is reduced and graphically equal to $uw^{i-j}u^{-1}$, following since $i, j$ have different signs and $w$ is cyclically reduced. It follows that different edges in $\mathcal{C}_1$ with the same source have different labels. This allows us to enlarge $\mathcal{C}_1$ by attaching few trees to the ``deficient'' vertices of $\mathcal{C}_1$,  in accordance with Subsection \ref{ss2}, and to obtain a labeled graph $\mathcal{G}_1$ in which the stars of all vertices are standard. Therefore, $\cg_1$ is the graph of the action of $F$ on $F/H_1$ where $H\subset H_1$, because $\mathcal{C} \subset \mathcal{C}_1$. 
%When we switched from $\mathcal{C}$ to $\mathcal{C}_1$ we added at least one vertex of degree less than $2r$ and to this we attached in $\mathcal{G}_1$ an infinite tree. It follows then that the new graph obtained is infinite provided that $l \ge 2$.
 
Finally, $g^n \in H_1$, for $n=i-j>0$, following because the label of a certain loop of the new graph $\mathcal{C}_1$, which is a subgraph of $\mathcal{G}_1$, can be read as the word $uw^nu^{-1}$ Moreover, by Lemma \ref{Lemma D}, $H_1=H*\langle g^n\rangle$.

To compare the deficits of the cores $\cc$ and $\cc_1$, we have to choose the above number $l$ so that $(2r-1)^{2-l/2} < \ve$. By Lemma \ref{Lemma D}, $\df(\cc)-\df(\cc_1) \le \ve$, as needed. Finally, from this inequality and our assumption about $\ve$ it follows that $\df(\cc_1)>0$, proving that the index of $H_1$ in $F$ is infinite.\bx
 
 \bigskip
 
 As mentioned in Introduction, a subgroup $K$ of a free group $F$ is called \emph{Burnside} if $g^n\in K$ for every element $g\in F$, where $n=n(g)$ is a positive integer. This property is equivalent to the fact that for any $g \in F$ all orbits of the right action of the subgroup $\langle g\rangle$ on $F/K$ are finite. Indeed,
 \bea
 (K u) v ^m = K(uv^m u^{-1})u = Ku\Longleftrightarrow     (uvu^{-1})^m \in K.
 \eqa

\begin{corollary}\label{c0021n} Any finitely generated subgroup $H$ of infinite index in the free group $F$ of rank $r >1$ is contained as a free factor in a Burnside subgroup $K$. One can choose $K$ with the maximal growth of action of $F$ on $F/K$. It follows that there exists a transitive action of $F$, with maximal growth and with finite orbits for each element $g\in F$, which factors through the action of $F$ on $F/H$.
\end{corollary}

\pp Let $g_1, g_2,\ld$ be the list of all elements in $F$. If we apply the first claim of Theorem \ref{t001n}, then we obtain an increasing chain of finitely generated subgroups of infinite index $H=H_0\subset H_1\subset H_2\subset \ld$ such that $H_{i-1}$ is a free factor of $H_i$, for all $i=1,2,\ld$, $g_i^{m_i}\in H_i$ for a sequence of positive integers $m_1, m_2,\ld$ (If $\langle g_i\rangle \cap H_{i-1} \ne {1}$ then we have $H_i=H_{i-1}$.) In this case $K=\bigcup H_i$ is a Burnside subgroup in $F$ and $H$ is a free factor of $K$. This subgroup has infinite index because otherwise it would be finitely generated and hence should coincide with one of the subgroups $H_i$, a contradiction. Thus the proof of the first claim is complete.

It follows from our remark before Lemma \ref{Lemma C} that the deficit of the action of $F$ on $F/H$ is positive, say, equal $c>0$. The second claim of Theorem \ref{t001n} enables us to choose the numbers $m_1, m_2,\ld$ co that $\df(\cc_i) > c- \dfrac{c}{3}- \cdots-\dfrac{c}{3^i}>\dfrac{c}{2}$ for the core $\cc_i$ of the action of $F$ on $F/H_i$.
According to Claim  (2) in Lemma \ref{Lemma C}, for the growth function $f_i(n)$ of the action of $F$ on $F/H_i$, we have $f_i(n) > \dfrac{c}{4r} (2r-1)^n$, for any $i, n$. Using Lemma \ref{Lemma D}, we obtain the sequence of integral inequalities $f_1(n)\ge f_2(n)\ge\ld$ which, obviously, must stabilize for any fixed $n$. As a result, for the limit function $\bar{f}$ we obtain $\bar{f}(n) \ge \dfrac{c}{4r} (2r-1)^n$.  Now $\bar{f}$ is the growth function for $F/K$ which follows since for any two words $g, g'$ of length at most $n$ we obtain $Kg=Kg'$  if and only if $H_ig=H_ig'$, for some $i=i(n)$. This proves the second claim of this corollary.\bx

%%%%%%%%%%%%%%%%%%%%%%%%%%%%%%%%%%%
\subsection{A construction based on generic properties}\label{ss7}
We will say that \emph{almost all reduced words  $w=w(x_1,\ld,x_r)$, $r>1$, have a certain property $\mathcal{P}$}, if the ratio of the number of all words of length $m$ without  $\mathcal{P}$ to the number $N_r(m)$ of all reduced words of length $m$ in $x_1,\ld, x_r$ tends to zero when $m \ra\infty$. Similarly one defines the properties which hold almost for all $k$-tuples for $k >1$.

For example, it is well-known that almost all words have the uniqueness property  of the occurrence of long subwords. More precisely, if one fixes any $\lambda >0$, then in almost any word $w$ any of its subwords of length $\ge \lambda |w|$ has a unique occurrence in $w$: $w$ is graphically equal to $uvu'$, for a unique pair of words $(u, u')$. 

Hint: 
\begin{enumerate}
\item[(a)] If $|v|\ge \lambda |w|$, and
$v$ occurs in $w$ in two different ways then one can find in $w$ two disjoint occurrences of a subword $v'$ of length $\ge (\lambda/3 ) |w|$, that is,  $w \equiv u_1v'u_2v'u_3$ ($\equiv$ is the graphic equality) . \item[(b)] The number of such words $w$ with fixed lengths $|u_1|, |u_2|, |u_3|$ is  exponentially small (with respect to $n=|w|$) when compared to $N_r(n)$ because $w$ is uniquely defined by the quadruple of words $(u_1,v',u_2,u_3)$ such that the sum of their lengths is at most $\le (1-\lambda/3)n$. 
\item[(c)] The number of different triples of lengths $|u_1|, |u_2|, |u_3|$ is polynomial in $n$.
\end{enumerate}

Similarly, for any fixed natural $k$ and any $\lambda  > 0$ almost all $k$-tuples of reduced words $(w_1,\ld, w_k)$ have the following property of uniqueness of the occurrence of long subwords : if $v$ is a subword of length $\ge \lambda |w_k|$ in one of $w_k^{\eta}$, $\eta =\pm 1$, then the words $u$ and $u'$ in the decomposition $w_k^{\eta}\equiv uvu'$ are uniquely defined by $v$, and $v$ does not occur in $w_l^{\delta}$ if $(l,\delta)\ne (k,\eta)$. We call this \emph{property $\mathcal{P}(k, \lambda)$}.

According to Subsection \ref{ss2}, there is a vertex $v$ with nonzero deficit in $\mathcal{C}=\mathcal{C}(H,o)$ if the subgroup $H\subset  F$ is finitely generated and has infinite index in $F$. It follows (see \cite{A}) that almost any reduced word $w$ has no subwords of length $\ge |w|/2$, which one can read on the paths of the graph $\mathcal{C}$. We call this \emph{property $\mathcal{P}$}[A].

%%%%%%%%%%%%%%%%%%%%%%%%%%%%%%%%%%%

Suppose now that again we have a finitely generated subgroup $H$ of infinite index in a free group $F$ of rank $r>1$. Let us choose a natural $k$ and fix in the graph $\mathcal{G}$ of action of $F$ on $F/H$ two $k$-tuples of the vertices $(v_1,\ld, v_k)$ and $(v'_1,\ld,v'_k)$, such that $v_i\ne v_j$ and $v'_i\ne v'_j$, for any $i\ne j$. Let us draw in $\mathcal{G}$ some reduced paths $p_1,\ld,p_k$ from $v_1,\ld,v_k$ with the same label $w$, and from $v'_1,\ld,v'_k$ the paths $p'_1,\ld, p'_k$ with label $w'$. The construction depends on the choice of reduced words $w$ and $w'$.
 
Let us note that the path $p_i=y_iz_iq_i$ where each $y_i$ is a path from $v_i$ through the forest $\mathcal{F}$ to $\mathcal{C}$; it is uniquely defined by $v_i$. Also $z_i$ is a path on $\mathcal{C}$ and again $q_i$ is a path through $\mathcal{F}$. Note that some of these paths may be missing. In a similar manner one defines the decompositions $p_i^{\prime}=y'_iz'_iq'_i$. Since the $y$-parts are fixed, it follows from the property $\mathcal{P}([A])$ that for almost any pair of words $w$ and $w'$ of length $\le m$ we will have $|q_i| > (2/5) |w| > m/3$ and $|q'_i| > (2/5) |w'|>m/3$. (One has to keep in mind that for almost any pair of words of length $\le m$ it is true that $|w|, |w'| > 5m/6$.)
 
Since the paths $q_i$ are entirely in the forest, for any pair of paths $q_i, q_j$ we either have that they have no vertices in common or $q_i = s(i,j)t(i,j)$, where $s=s(i,j) = s(j,i)$, while $t(i,j)$ and $t(j,i)$ have no edges in common. Since $s$ is a subpath of different paths $p_i$ and $p_j$, with common label $\Lab{s}$ (recall that $(p_i)_-
\ne (p_j)_-$),  and the stars of vertices in $\mathcal{G}$ are regular,  $\Lab{s}$ has two different occurrences in $w$. It follows by $\mathcal{P}(k,1/12)$, that for almost all $w$ we have $|s(i,j)|<|w|/12\le m/12$, and hence 
$|t(i,j)|=|q_i|-|s(i,j)|>(1/3-1/12)m=m/4$.  Similarly,  $|t'(i,j)|>m/4$ for almost all $w'$.
 
 In a similar manner one can compare the common parts $\bar{s}(i,j)$ in the paths $q_i$ and $q'_j$, for any $i,j$. Their labels produce the same subwords in $w$ and $w'$. It follows by the property $\mathcal{P}(2, 1/12)$ that these common parts have length $<m/12$ for almost any pair of words $w, w'$ of length $\le m$.
 
It follows that our construction (depending on the words $w, w'$) has the following property:

(*) \emph{For almost any pair of words $w,w'$ of length $\le m$, all paths $p_i, p'_i$ end with the subpaths $t_i, t'_i$ of length $\ge m/4$ that follow through the forest $\mathcal{T}$ in the direction FROM the subgraph $\mathcal{C}$ in such a way that no two of these $2k$ subpaths have common vertices}.

In the statement of the following theorem that uses (*), given a subgroup $H$ (respectively, $H_1$) of a free group $F$, we denote by $\cc$ (respectively, $\cc_1$) the core of the graph of the action of $F$ on $F/H$ (respectively, $F/H_1$).

\begin{theorem}\label{t002n} Let $H$ be a finitely generated subgroup of infinite index in a free group $F$ of rank $r>1$.
Let $(Hg_1,\ld,Hg_k)$ and $(Hg'_1,\ldots,Hg'_k)$ \emph{(}$k\ge 1$\emph{)} be two $k$-tuples of pairwise different cosets. Then for any $\ve\in(0,\df(\cc))$ there are in $F$ a finitely generated subgroup $H_1$ of infinite index in $F$ with $0 \le  \df(\cc) - \df(\cc_1)\le \ve$, and an element $b\in F$, such that $H_1g_ib = H_1g'_i$, for any $i=1,\ldots,k$. Additionally, $H$ is a free factor of $H_1$.
\end{theorem}
\pp  Let two $k$-tuples $(v_1,\ld,v_k)$ and $(v'_1,\ld,v'_k)$ of pairwise different vertices of the graph $\mathcal{G}$ of action of $F$ on $F/H$ correspond to the $k$-tuples of the cosets $(Hg_1,\ld,Hg_k)$ and $(Hg_1^{\prime},\ld,Hg_k^{\prime})$. We select two reduced words $w$ and $w'$ and perform the construction in $\mathcal{G}$, according to the construction  preceding the statement of this theorem. Let us denote by $\Gamma$ the minimal connected subgraph in $\mathcal{G}$ containing $\mathcal{C}$ and all paths $p_i, p'_i$.

By Property (*), the words $w$ and $w'$ can be selected in such a way that $|t_i|, |t'_i|>m/4$, $m>2l$ for an arbitrarily large $l$, and 
all vertices in $t_i, t'_i$ have degree at most $2$ in $\Gamma$, while the endpoints $o_i, o'_i$ have degree 1.

For the last edges $e_i$ of the paths $t_i$ the label is one and the same letter $x \in {x_1^{\pm 1},...x_r^{\pm 1}}$, since this letter is the last in the word $w$. Similarly, all $e'_i$ have the same label $x'$. We will chose $w$ and $w'$ with different last letters. This is possible because, if necessary, we can always make the word $w^{\prime}$ longer by attaching an appropriate letter at the end.

By definition, the graph $\mathcal{C}_1$ can be obtained from $\Gamma$ by identifying the endpoints $o_i$ and $o'_i$ of the paths $p_i$ and $p'_i$. Since $x\neq x^{\prime}$, in $\mathcal{C}_1$ no two different edges with the same source have equal labels. It is also obvious that the degrees of all vertices in $\mathcal{C}_1$ are at least two, with possible exception of $o$. It follows that attaching a forest to $\mathcal{C}_1$ (see Subsection \ref{ss2}) results in a graph $\mathcal{G}_1$ with labeling such that the stars of all vertices are standard. Then $\mathcal{G}_1$ is the graph of action of $F$ on $F/H_1$, where $\mathcal{C}_1=\mathcal{C}(H_1)$; since $\mathcal{C} \subset \mathcal{C}_1$, it follows that $H \subset H_1$, and by Lemma \ref{Lemma D} it follows that $H$ is a free factor of
$H_1$. The vertices $v_1,\ld,v'_k$ are also in  $\mathcal{C}_1$ but here they define the cosets $H_1g_1,\ld,H_1g'_k$. Since $v_i$ and $v'_i$ are connected in $\mathcal{G}_1$ by the paths $r_i$ with label $ww'^{-1}$, which do not depend on $i$, for all $i=1,\ld,k$ we obtain: $H_1g_i b = H_1g'_i$ where $b=ww'^{-1}$.

Finally, in order to estimate the difference $\df(\cc)-\df(\cc_1)$, we will apply Lemma \ref{Lemma D} $k$ times. Namely, when on the $i^\mathrm{th}$ step we identify $o_i$ with $o'_i$, we attach to the core of the graph, arising after gluing together $i-1$ previous pairs of vertices, an elementary subpath of length greater than $m/4+m/4\ge l$, because by the properties of the paths $t_i, t'_i$, all their vertices keep their degrees after the first $i-1$ steps (``nothing sticks to them''). By Lemma \ref{Lemma D}, after $k$ steps we will obtain
\bea
0 \le  \df(\cc) - \df(\cc_1) \le  k(2r-1)^{2- l/2 }<\ve,
\eqa
for  $l$ sufficiently big. As in Theorem \ref{t001n}, from the latter inequality and our assumption about $\ve$ it follows that $\df(\cc_1)>0$, proving that the index of $H_1$ in $F$ is infinite.\bx

\subsection{Main corollaries}\label{ssMC}

\begin{corollary}\label{c0022n} Any finitely generated subgroup $H$ of infinite index in the free group $F$ of rank $r>1$, is contained as a free factor in a subgroup $K$ of infinite index in $F$ such that for any natural $k$ the right action of $F$ on $F/K$ is $k$-transitive. In particular, $K$ is a maximal subgroup in $F$. One can choose $K$ in such a way that the growth of the action of $F$ on $F/K$ is  maximal.
\end{corollary} 

\pp Let us enumerate the pairs of tuples $(g_1,...g_k)$, $(g'_1,\ld,g'_k)$ of pairwise different elements in $F$, for all natural $k$. Set $H_0 = H$ and assume the sequence of subgroups $H_0 \subset H_1\subset\ld \subset H_{i-1}$ already constructed, where $H_{i-1}$ is a free factor in $H_i$, for all $i=1,2,\ld$. If for the $i^\mathrm{th}$ pair of tuples all classes $H_{i-1}g_j$ are pairwise different for $g_1,..,g_k$ and all $H_{i-1}g'_j$ also different then we proceed as follows and if not, we set $H_i=H_{i-1}$. By Theorem \ref{t002n} there is a finitely generated subgroup $H_i$ of infinite index in $F$ such that  $H_{i-1}$ is a free factor of $H_i$, and element $b\in F$ such that $H_i g_j b = H_ig'_j$, for all $j=1,\ld,k$. 

We set $K= \bigcup H_i$. Then $K$ is of infinite index in $F$ because all $ H_i$ are of infinite index and $H$ is a free factor of $K$. For any two $k$-tuples of pairwise different cosets $Kg_1,\ld,Kg_k$ and $Kg'_1,\ld, Kg'_k$ the $k$-tuples $H_{i-1}g_1,\ld,H_{i-1}g_k$ and $H_{i-1}g'_1,\ld, H_{i-1}g'_k$ also consist of pairwise different elements each for all $i$. Hence, for some $i$ and some $b\in F$, we will obtain $H_ig_jb=H_ig'_j$, for all $j=1,\ld,k$. It is immediate then that $Kg_jb=Kg'_j$, for all $j$ and the $k$-transitivity of the action of $F$ on $F/K$ follows.

Let us recall that any $2$-transitive action of any group $F$ is primitive and that an action is primitive if and only if the stabilizer $K$ of any point is a maximal subgroup in the acting group (that is, if $K\subset L\subset F$ then $L=K$ or $L= F$ where $L$ is a subgroup of $F$). It follows that the subgroup $K$ obtained is maximal in $F$.
 
The maximality of the growth for $F/K$ follows in the same way as in Corollary \ref{c0021n} but one has to replace reference to Theorem \ref{t001n} by reference  to Theorem \ref{t002n}. Now the proof is complete.\bx
  
One can alternate the steps in the proofs of Corollaries \ref{c0021n} and \ref{c0022n} to obtain the following.

\begin{corollary}\label{c003n} Any finitely generated subgroup $H$ of infinite index in a free group $F$ of rank $r>1$, is a free factor in a Burnside subgroup $K$ of infinite index such that for any natural $k$ the action of $F$ on $F/K$ is $k$-transitive. One can choose $K$ so that the growth of the action of $F$ on $F/K$ is maximal.\bx
\end{corollary}

The results obtained so far have an application to the modules of maximal growth as follows.

\begin{corollary}\label{c004n} Let $\Phi$ be a field. Then there is a module $M$ of maximal growth over the free group $F=\fgax$, or, equivalently, over the free group algebra $R=\gax$, both  of rank $r>1$, satisfying the following additional properties.
\begin{enumerate}
\item[\emph{(a)}] The module $M$ is monomial, that is, induced from a trivial one-dimensional module of a subgroup of $F$;
\item[\emph{(b)}] The module $M$ has a simple submodule $N$ of codimension $1$ (hence the growth of $N$ is also maximal);
\item[\emph{(c)}] The modules $M$ and $N$ are \emph{periodic} in the sense that for any $a\in M$ and $g\in F$ there is a positive $m=m(a,g)$ such that $ag^m=a$ (In other words, the orbits of the action of any cyclic subgroup $\langle g\rangle$ of $F$ on $M$ are finite.)
\end{enumerate}
\end{corollary}

\pp Let us choose a linear space $M$ with basis $\{ e_i, i\in F/K\}$ where $K$ is the subgroup from the previous corollary. We expand by linearity the action of $F$ on $F/K$ to $M$. It is obvious that $M$ is induced from $1$-dimensional trivial $K$-module $L$: 
$M = L \otimes_{\Phi K} \Phi F$. 
 
The growth of $M$ is  maximal since the growth of the action of $F$ on $F/K$ is maximal.
  
As in Corollary \ref{cNMMG}, we obtain a simple module of maximal growth over $\gax$ if we consider the subspace $N$ in $M$ consisting of all finite $\Phi$-linear combinations $\sum \lambda_i e_i$ with $\sum \lambda_i =0$. The proof of simplicity and the maximality of the growth are exactly the same as well.
   
Now let $a=\sum \lambda_i e_i \in M$ and $g\in \mathcal{G}$. Since the subgroup $K$ is Burnside, the $\langle g\rangle$-orbit of each $e_i\in F/H$ is finite and to obtain the equality $a g^m =a$ one has to set $m$ equal the least common multiple of all $m_i$ such that $e_i g^{m_i} =e_i$ for all vectors in the decomposition of $a$ as above.
 
Now the proof is complete.\bx

\section{Other properties of maximal growth}\label{sOPMG}

\subsection{Topological approach to maximal growth}\label{ssTAMG}

We will use the notation and some facts from of Subsections \ref{ss1} and \ref{ssEG}.

Given a free group $F_r$ with fixed symmetric basis (alphabet) $B=A\cup A^{-1}$,  $A=\{ a_1,\ld,a_r\}$, we denote by $\partial F_r$ the set of all functions $w:\{ 1,2,\ld\}\ra B$ such that $f(n)\neq f(n-1)^{-1}$, for any $n=2,3,\ld$. One can view each such function as a \emph{semi-infinite} to the ``right'', or \emph{right-infinite} reduced word in the alphabet $B$.  One can also view $\partial F_r$ as the set of (extended) labels of infinite reduced rays in the Cayley graph of $F_r$ originating in $1$.
 
It is well-known that $\partial F_r$ can be turned into a metric space if one defines the distance between $w_1, w_2\in \partial F_r$ as $\dfrac{1}{(2r-1)^n}$, where $n$ is the length of the maximal common prefix $u$ of $w_1$ and  $w_2$. As a result, $\partial F_r$ becomes an \emph{ultrametric compact} space. A basis of topology in $\partial F_r$ is given by the open subsets 
\bea
O_u=\{ s\in \partial F_r\,|\,u\mbox{ is a prefix of }s\}.
\eqa
 One can equip $\partial F_r$ with a \emph{countably additive measure} $\mu$ such that   
\bea
\mu(\partial F_r)=1\mbox{ and }\mu(O_u) = \frac{(2r-1)^{-|u|+1}}{2r}\mbox{ for }|u|\ge 1.
 \eqa
 Now suppose $\cg$ is the Cayley graph of the action of $F_r$ on $F_r/H$ and $\ct$ is a maximal subtree in $\cg$. We denote by $Y=Y(\ct)$ the subset of $\partial F_r$ consisting of right-infinite words that can be read on the infinite rays of $\ct$ originating in $o$.  The set $Y$ is closed in $\partial F_r$. Indeed, if $w$ is a limit point for $Y$ then $Y$ includes the words which have arbitrarily long common prefixes with $w$. It follows that all prefixes of $w$ are in $\ct$ hence $w\in Y$. We may conclude that $Y$ is \emph{measurable}.
 
\begin{theorem}\label{lMEASURE} If the growth of $\cg$ is maximal then the measure of the subset $Y=Y(\ct)$, $\ct$ a fixed maximal geodesic Schreier subtree  in $\cg$, is positive. Conversely, if there is a maximal geodesic Schreier subtree $\ct$ in $\cg$ such that the measure of the respective set $Y$ is positive then the growth of $\cg$ is maximal.
\end{theorem}

\pp Suppose the growth is maximal. Then by Lemma \ref{lMG} there is $c>0$ such that for any ball $B(n)=\mb{o}{n}$ with center $o=H$ and radius $n>0$ we must have $\# B(n) \ge c(2r-1)^n$.  

Suppose $\mu(Y)< c^2\dfrac{2r-1}{10r}$. Then there exists $\ve<c^2\dfrac{2r-1}{10r}$, such that we can choose a countable open covering of the set $Y$ by the subsets $O_u$, with total measure $<\ve$. Being a closed subset of a compact set $\partial F_r$, our set $Y$ is itself compact and so we may assume that our covering is finite, say, $Y\subset O_{u_1}\cup\cdots\cup O_{u_t}$. By our condition on the total measure, we have $\sum (2r-1)^{-|u_i|}<2r \dfrac{\ve}{2r-1}< c^2/5$.
 
  Now let us consider all possible spheres $S_n=\ms{o}{n}$ of radius $n$ with center $o$ in $\cg$, $n=1,2,\ld$. 
  
 \textbf{Case 1.} There is $n>0$ such that the number of elements in the set $V_n$ of vertices on $S_n$ which belong to the infinite rays in $\ct$ is bounded by $\dfrac{c}{2} (\# S_n) = cr(2r-1)^{n-1}$. Then the vertices $x$ of $S_n \setminus V_n$ can only belong to the finite branches of $\ct$ starting from $o$.  It follows from the compactness principle that the number of such branches can only be finite because otherwise $x$ would also belong to an infinite ray of $\ct$.  It follows that if $m>n$ then we would have $\# S_m \le (\# V)(2r-1)^{m-n} +C$, where $C$ does not depend on $m$.
In this case, $\# S_m \le cr (2r-1)^{m-1}+C$, and $\# B(m)  \le c (2r-1)^m$, for all sufficiently large $m$.  This contradicts our choice of $c$.

  \textbf{Case 2.} For any $n\ge 1$, we have that $\# V_n >  \dfrac{c}{2} (\# S_n) = cr(2r-1)^{n-1}$. Each vertex of $V_n$ is an element of both $S_n$ and an infinite reduced in $\ct$, hence a ray in one of $O_{u_i}$. For each $i\le t$, the number of such rays is at most $\max (1, (2r-1)^{n-|u_i|})$, because $u_i$ is the prefix of the (infinite) label of each such ray.  As a result, $\# V_n \le  t + (2r-1)^n K$, where $K = \sum (2r-1)^{-|u_i|}<2r\dfrac{\ve}{2r-1}$. 
 
 It follows that, for all spheres whose radius is sufficiently large, we have $\# S_n\le 2\dfrac{\# V_n}{c}  < 3\ve \dfrac{(2r-1)^n}{c}$, hence for the balls we have $\# B(n) < 5\dfrac{\ve(2r-1)^n}{c} < c(2r-1)^n$, by the choice of $\ve$.  Again, we have a contradiction with the maximality of the growth. Thus we have shown that if the growth is maximal then the $Y$ is the set of positive measure.
 
 Conversely, suppose $\mu (Y) = \mu > 0$.  Let us set $s_n=\# S_n$.  Then the number of different prefixes of length $n$ in the word from $Y$ is at most $s_n$ hence $Y$ can be covered by $s_n$ different subsets $O_u$ with $\mu(O_u)<(2r-1)^{-n}$. Therefore,  $\mu=\mu(Y) < s_n(2r-1)^{-n}$. Hence $s_n > \mu (2r-1)^n$.  Then also $\# B(n) > \mu(2r-1)^n$, for all $n>0$, and hence the growth of $\cg$ is maximal by Lemma \ref{lMG}.  \bx
 
 Similar topological characterization works also in the case of cyclic acts over a free monoid $W=\wax$, $\# X=r>1$. Again, we have to consider the ultrametric space $\partial W$, with $2r-1$ replaced by $r$ when we define the metric and the measure. For a cyclic act $S$ over $W$ we can define a directed graph with labelling, in the same way as we defined the graph $\cg(S)$ in Subsection \ref{ss1}. One can choose in $\cg$ a directed geodesic Schreier subtree $\ct$ (all edges are directed from the fixed vertex $o$, corresponding to the generator of the act), select a measurable subset $Y(\ct)$ and proceed in the way described above for the groups.  The conclusion is the same: \emph{the growth of $S$ is maximal if and only if $\mu(Y(\ct))>0$}. 
 
In the case of modules over a \fac respectively, a free group algebra $R=\Phi M$, $M$ the free monoid, respectively, the free group, the topological characterization of cyclic modules of maximal growth completely reduces to the two cases described above. Indeed, the ambient ultrametric topological space for $R$ is $\partial M$, $M$ as above. Now let us consider a cyclic $R$-module $V=R/I$, where $I$ a right ideal of $R$. By \cite{JL}, there is a Schreier system $\mathcal{S}$ of monomials which serve as representatives of elements of $R$ modulo $I$. The system $\mathcal{S}$ is prefix closed and geodesic in the sense that none of its terms can be written modulo $I$ as a linear combination of lesser monomials with respect to ShortLex. One can view $\mathcal{S}$ as the set of labels written on the branches of a tree $\ct$. As previously, this gives rise to the subset $Y(\ct)\subset\partial W$, and we can proceed in the same way as previously with $\partial M$, $M=\wax$ or $M=\fgax$. Again, \emph{the growth of $R/I$ is maximal if and only if $\mu(Y(\ct))>0$.} Naturally, in the proofs one has to replace the cardinality $\#$ by dimension $\dim$.

\subsection{Growth and semi-isomorphisms of $F_r$-sets}\label{ssSI}

As it was mentioned in subsection \ref{ssGAI}, the growth is invariant under the isomorphisms of $F_r$-sets. Still, it is reasonable to ask what happens if we replace isomorphisms of $F_r$-sets by more general ``semi-isomorphisms''. Given $F_r$-sets $S$ and $S'$, a bijection $f: S\ra S'$ is called a  \emph{semi-isomorphism} of $F_r$-sets if there is an automorphism $\vp : F_r\ra F_r$ such that $f(x\cir g)=f(x)\cir \vp(g)$ for any $x\in S$ and $g\in F_r$. In other words, we would like to know what happens to the growth if in our definition of the growth functions in subsection \ref{ssGF} we replace a filtration of $F_r$ associated with one free basis by a filtration associated with another free basis. A simple example shows that the growth may change. Suppose $S=\mathbb{Z}$. We define an action of $F_2=F(a,b)$ on $S$  by $m\cir a=m\cir b=m+1$, for any $m\in S$. Then $\# \mb{0}{n}=2n+1$ and the growth is equivalent to $2n$. However, if we replace the free basis $\{ a,b\}$ by $\{a,ab\}$ then the ball of radius $n$ will contain $4n+1$ numbers and so the growth is equivalent to $4n$. As we noted in Introduction, $2n$ is not equivalent to $4n$. Similar examples work in all four cases considered by us in this paper. (As it is known, a semi-linear isomorphism of modules need not be an isomorphism.)

An example of an $F_r$-set whose growth is maximal for one free basis of $F_r$ and not maximal for another is by far less obvious. The goal of this subsection is to show that the notion of maximality of the growth does depend on the choice of a free basis in $F_r$. This means that when we speak about the maximality of the growth we have to keep in mind a free basis of $F_r$. However, this does not blur the notion of maximality of the growth of an action since usually a free group comes with its free basis fixed, as it happened in the situation described in the Introduction when the action of $F_r$ was defined by the map $a_i\ra \mathcal{A}_i$, where $\{ a_1,\ld a_r\}$ is a free basis of $F_r$ and $\mathcal{A}_i$ are fixed transformations of the set $S$.

Let us denote by $w_v$ the number of different occurrences of a word $v$ in a reduced word $w$ and set $s_n = \# S_n  = 2r(2r-1)^{n-1}$ the number of elements of length $n \ge 1$ in $F_r$. We we will need the following result, Proposition 5.3 from a recent paper \cite{KSS}, based on the
Large Deviation Theory.
     
\begin{lemma}\label{K-S-S} Let $r$ be an integer, $r>1$. Then the following are true
\emph{(1)} For any $\ve >0$ and any letter in the symmetric alphabet $\{ a_1,\ld,a_r,a_1^{-1},\ld, a_r^{-1}\}$ it is true that
\bea    
    \lim   \frac{\#\left\{ w\in F_r\,\left|\, |w|=n\mbox{ \emph{and} }\frac{w_{a}}{n} \in \left(\frac{1}{2r} - \ve, \frac{1}{2r}+\ve\right)\right.\right\}}{s_n} =1.
\eqa    
\emph{(2)} For any two letters $a, b$ such that $b \ne a^{-1}$, and any $\ve >0$ it is true that
\bea     
    \lim   \frac{\#\left\{ w\in F_r\,\left|\, |w|=n\mbox{ \emph{and} }\frac{w_{ab}}{n} \in \left(\frac{1}{2r(2r-1)} - \ve, \frac{1}{2r(2r-1)}+\ve\right)\right.\right\}}{s_n} =1 .
 \eqa    
    In addition, the rate of convergence in both limits is exponential.
\end{lemma}

Given positive $\ve$ and natural $l$, we denote by $Z_{\ve, l}(n)$ the set of all reduced words $w$ of length $n\ge l$ such that for any $m \in [l,  n]$ and any prefix $v$ of length $m$ of $w$, the following inequalities hold:
\bea 
\frac{v_u}{m}  \in \left(\frac{1}{2r}    -  \ve, \frac{1}{2r}    +  \ve\right)
\eqa
as soon as $|u|=1$ and  
 \bea  
 \frac{v_u}{m}  \in \left(\frac{1}{2r(2r-1)}    -  \ve, \frac{1}{2r(2r-1)}    +  \ve\right)
 \eqa
as soon as $u$ is a reduced word with $|u|=2$.
     
\begin{lemma}\label{l'} For any $\ve >0$ there is $l= l(\ve, r)$ such that the ratio $\dfrac{\# Z_{\epsilon,  l}(n)}{s_n}$ is a monotonously decreasing function of $n$ whose limit is a positive number $\nu = \nu (r, \ve)$.
\end{lemma}  
      \pp  It is obvious that the function is decreasing since $s_{n+1}=(2r-1)s_n$ and each word in $Z_{\ve, l}(n+1)$ is the product of a word in $Z_{\ve, l}(n)$ by one of $2r-1$ letters.

Next, Let $Y_{\ve}(n)$ be the set of all words $w$ of length $n$ such that the above inequalities hold for $\frac{w_u}{n}$ only, not necessarily for all possible subwords of various lengths $m$. Finally, set $X_{\ve}(n) = S_n \setminus Y_{\ve}(n)$.

    Using Lemma \ref{K-S-S}, the sequence $\dfrac{\# X_{\ve} (n)}{s_n}$ converges to 0 at exponential rate. In other words, $\dfrac{\# X_{\ve} (n)}{s_n}$ is bounded from above by $2^{-\delta n}$, for some positive $\delta =\delta(\ve,r)$ and all $n\ge l$.

Let us increase $l$ just chosen so that, additionally, $\sigma = \sum_{i=l}^{\infty} 2^{-\delta i} < 1$.

Let us prove that for any $n\ge l$ we always have $\dfrac{\# Z_{\ve,l}(n)}{s_n} > 1-\sigma_n$, where $\sigma_n = \sum_{i=l}^n 2^{-\delta i} <\sigma$. Once this is done, the second claim of our lemma will hold with $\nu = 1-\sigma$.
     
  If $n= l$ then the inequality in question will follow from $Z_{\ve,l}(l)=Y_{\epsilon}(l) = S_n \setminus X_{\epsilon}(l)$ and by the choice of $l$. For the induction step from $n$ to $n+1$, we notice that some of the products of the words in $Z(n)=Z_{\ve,l}(n)$ on the right by one letter are not in $Z(n+1)$;  but then they will land in $X(n+1)=X(\ve, n+1)$. Using the induction hypothesis and the equation $s_{n+1}=(2r-1)s_n$, we will obtain the following:
\bea    
       && \frac{\# Z(n+1)}{s_{n+1}} \ge \frac{(2r-1)(\# Z(n)) - \# X(n+1)}{s_{n+1}} =
    \frac{\# Z(n)}{s_n}- \frac{\# X(n+1)}{s_{n+1}}\\&& \ge 1-\sigma_n - 2^{-\delta(n+1)} = 1 - \sigma_{n+1}.
 \eqa  
Now the proof is complete. \bx  

\begin{remark}\label{rC1}
It follows by the definition of the numbers $\sigma$ and $\nu$ that $\nu$ can be made arbitrarily close to $1$, if we choose $l$ appropriately.
\end{remark}

 Now let us fix $\ve< \dfrac{2r-3}{6r(2r-1)}$ and choose $l$ in agreement with Lemma \ref{l'}. Using an approach similar to that used in \cite{KKS}, we check that the stretching coefficient $\lambda$ of the Nielsen automorphism $\vp :  a_1\ra a_1$, $a_2\ra a_1a_2$, $a_i\ra a_i (i>2)$ is strictly greater than $1$ on any word of any set $Z(n)$, $n\ge l$. Namely, the following is true.
 
\begin{lemma}\label{l''} There exists $\lambda >1$ such that, for the automorphism $\vp$ just defined and any $w\in Z(n)=Z_{\epsilon,l}(n)$ \emph{(}$n\ge l$, where $l$ is chosen above\emph{)}, we have $|\vp(w)|>\lambda |w|$.
\end{lemma}
     
\pp We write $w = x_1\cdots x_n$, where $x_1,\ld,x_n$ are not necessarily different letters of the symmetrized alphabet $\{ a_1^{\pm 1},\ld,a_r^{\pm 1}\}$. Applying $\vp$ to every letter $w$, we will obtain a not necessarily reduced product $v = \vp(x_1)\cdots \vp(x_n)$.  By definition of $Z(n)$, $w$ has at least $\left(\dfrac{1}{2r}-\ve\right)\!n$ entries of $a_1$ (and at least the same number of entries of $a_1^{-1})$. Thus the number of factors of the form $a_1a_2$ (or  $a_2^{-1} a_1^{-1}$) in the given decomposition of $v$ is greater than $\left(\dfrac{1}{2r} -\ve\right)\!n$. It follows that before we apply any cancellations, the length of $v$ is greater than $n + \left(\dfrac{1}{r}-2\ve\right)\!n$.
 
To reduce $v$ we only need to apply cancellations of the form $a_1^{-1}(a_1a_2)\ra a_2$,  $(a_2^{-1}a_1^{-1})a_1\ra a_2$ because after them no further cancellations are possible. These cancellations correspond to the occurrences of the 2-letter words $a_1^{-1}a_2$ and $a_2^{-1} a_1$ in the original word $w$. By definition of $Z(n)$, the number of occurrences of these $2$-letter words in $w$ is less than $\left(\dfrac{1}{2r(2r-1)} +\ve\right)\!n$. Hence, after all cancellations has been performed in $w$, the number of letters that vanish is less than $4\left(\dfrac{1}{2r(2r-1)}+\ve\right)\!n$.  As a result, the length of the reduced word $\vp(w)$ is greater than
\bea     
    n+\left(\dfrac{1}{r}-2\ve\right)\!n - \left(\dfrac{2}{r(2r-1)}+4\ve\right)\!n = \left(1+ \dfrac{2r-3}{r(2r-1)}- 6\ve\right)\!n.
    \eqa  
    By our choice of $\ve$, the coefficient $\lambda = 1+ \dfrac{2r-3}{r(2r-1)}- 6\:\!\ve$ is greater than $1$, and the proof of the lemma is complete.\bx
    
\begin{remark}\label{rPREFIX}
From the description of cancellations in the word $v$, it follows that any prefix of a reduced form of $\vp(w)$ can be obtained from the image $\vp(w_1)$, where $w_1$ is a prefix of $w$, by multiplication on the right by a word whose length is at most $1$.
\end{remark}

 Now let $\ve$ and $l$ are chosen as just before Lemma \ref{l''}. Set $Z =\bigcup_{n=0}^{\infty} Z(n)$, where $Z(n)$ in the case where $n<l$ is defined as the set of all reduced words of length $n$. Let us denote by $V$ the set of all reduced words in $F_r$ representing the elements in $\vp(Z)$.
 
 Let $B$ be any symmetric subset in the fixed symmetric generating set $\{ a_1^{\pm 1},\ld,$ $a_r^{\pm 1}\}$ of $F_r$, that is, $B$ is closed under inverses of its elements. By definition, the subset $V_B$ will be obtained if we replace in each word from $V$ all occurrences of letters from $B$ by their inverses. We will ``symmetrize'' $V$ to obtain a set $U$ as the union of all $V_B$. Finally, let $\bU$ be the closure of $U$ under taking prefixes. We denote by $V(n)$, $U(n)$, $\bU(n)$ and $F_r(n)$ the sets of all elements of length at most $n$ in each of the sets $V$, $U$, $\bU$ and $F_r$, respectively. 

\begin{lemma}\label{l'''} The ratio $\dfrac{\# \bU(n)}{\# F_r(n)}$ converges to zero at exponential rate.
 \end{lemma}
     
\pp By Lemma \ref{l''}, any sufficiently long word $w\in V$, say with $|w|=n\ge n_0$, for some $n_0$, is of the form $\vp(v)$,  where $\lambda |v| <|w|$, t.e $v \in F_r \left(\dfrac{n}{\lambda}\right)$.
    Thus $V(n)$ has less than $2(2r -1)^{n/\lambda} +C$ words, for some constant $C=C(r,n_0)$, and at the same time, $\# F_r(n) > (2r-1)^n$. It follows that the ratio $\dfrac{\# V(n)}{\# F_r(n)}$ exponentially fast goes down to $0$. The same will happen to the ratio $\dfrac{\# U(n)}{\# F_r(n)}$, since replacing some letters by their inverses does not change the lenghth of the words.
    
    In what follows we will show that any word in $\bU$ can be obtained as a result of multiplying a word in $U$ by a word whose length is at most 1. Thus $\# \bU(n) \le 2r(\# U(n))$, and the claim of the lemma follows from the estimate of a previous paragraph.
     
 Obviously the above property of $\bU$ follows from the same property for $\bV$ (the closure of $V$ by prefixes) when we compare it with $V$. Finally, for $V$ and $\bV$ this property follows from Remark \ref{rPREFIX} because, since by definition of $Z$, this set is closed under prefixes. Now the proof is complete.\bx
 
 Since the set $\bU$ is closed under prefixes, in the Cayley graph of $F_r$, with respect to generators  $\{ a_1^{\pm 1},\ld,a_r^{\pm 1}\}$ there is a subtree $\mathcal{T}$ such that the set of labels of reduced paths $o - o^{\,\prime}$ of $\ct$ from the origin $o$ to an arbitrary vertex $o'$ is precisely $\bU$.  
Now let us consider the labelled graph $\ct$ without its connection to the ambient Cayley graph of $F_r$. Then some of the vertices of $\mathcal{T}$ may carry a nonzero deficit, that is, for such a vertex, there are no outgoing edges labelled by some $x\in\{ a_1^{\pm 1},\ld,a_r^{\pm 1}\}$ and no incoming edges with inverse labels. Such a vertex will be called \emph{$x$-deficit}. Notice that thanks to the symmetry of conditions on the words in $U$ and $\bU$ with respect to the change $x\leftrightarrow x^{-1}$, on each level $n$ (on the sphere $\ct(n)$ of $\ct$) the number of $x$-deficit vertices is the same as the number of $x^{-1}$-deficit vertices.

Now, by adjoining only edges, but not vertices, we will embed  $\mathcal{T}$ (with all distances to $o$ preserved) in a graph of transitive action of $F_r$ on the same set of vertices (which we identify with $\bU$). Namely, for each $x$-deficit vertex $v$ on an arbitrary level $n$ we will find an $x^{-1}$-deficit vertex $v'$ on the same level and draw an edge with label $x$ from $v$ to $v'$ and  an edge with label $x^{-1}$  from  $v'$ to $v$. We keep doing this until we obtain a connected graph $\mathcal{G}$ all of whose vertices are standard in the sense of Subsection \ref{ss1} and then, as it is mentioned there,  $\mathcal{G}$ is the graph of transitive action of $F_r$ on the set $F_r/H$ of right cosets of the stabilizer $H$ of the vertex $o$ in $F_r$.

The subtree $\mathcal{T}$ of the graph $\mathcal{G}$ is a geodesic Schreier subtree. This follows because additional edges connecting vertices of the same level cannot change distances to $o$. Therefore, the reduced words written on all paths $o - o^{\,\prime}$ in  $\mathcal{T}$ form a geodesic transversal for the subgroup $H$. The growth of $\mathcal{G}$ is not maximal by Lemma \ref{l'} and \ref{lMG}; moreover, by Lemma \ref{l'} it can be majorated by a function of the form $(2r-1-\gamma)^n$, for some $\gamma>0$.
     
 Let us keep the same right action of $F_r$ on $F_r/H$, but now consider the graph $\mathcal{G}'$ of this action with respect to a new basis $\{ b_1=a_1,
    b_2= a_1a_2 , b_i = a_i (i>2)\}$ of the group $F_r$.  Then $V(\mathcal{G}')=V(\mathcal{G})$, but the edges (and their labelling) are different.
    
    Let $R(n)$ be the subset of those vertices $o'$ in  $\mathcal{G}'$ (or in $\mathcal{G}$) with which $o$ is connected in $\mathcal{G}$ by a path with label $w \in V$, where $w$ is the reduced form of $\vp(v)$, $v\in Z(n)$.
    
    It is obvious that after the change of our basis to $\{ b_1^{\pm 1}, \ld,b_r^{\pm 1}\}$ any vertex $o'\in R(n)$ will be connected with $o$ in $\mathcal{G}'$  by a path whose label is a $b$-copy of the word $v^{-1}$,
because $b_i =  \vp(a_i)$, for $i=1,..,r$. Since $v \in Z(n)$, for the estimation of the growth of $\mathcal{G}'$ from below it is sufficient to obtain the lower bound for the sequence of numbers $\# R(n)$. By construction, the vertices $o', o''$ of $\mathcal{G}$  are different when $w', w'' \in V$ are different, hence when $v', v''\in Z(n)$ are different. It follows that $\# R(n) \ge \# Z(n)$. Hence the growth of $\mathcal{G}$ is not smaller then the growth of the sequence of numbers $\# Z(n)$. By Lemma \ref{l'} and Lemma \ref{lMG}, this growth is maximal.

The comparison of the estimates obtained by us for the graphs $\mathcal{G}$ and $\mathcal{G}'$ allows us to draw the following conclusion.

\begin{theorem}\label{tMNMG} For any $r >1$ there is a transitive action of the free group $F_r$, whose growth with respect to one free basis of $F_r$ is maximal and with respect to another basis of $F_r$ it is not just maximal but actually bounded from above for sufficiently large values of $n$ by a function of the form  $(2r-1-\gamma)^n$, for some $\gamma >0$. \bx 
\end{theorem}

\begin{remark}\label{rSA} It is easy to check that the growth of the action of $F_r$ remains to be the same when the change of bases is performed by a permutation of elements of a basis or by an inner automorphism of $F_r$. In \cite{KKS} the products of the respective automorphisms are called \emph{simple}. At the same time, on the basis of a more thorough use of \cite{KSS} and \cite{KKS}, one can expand the effect outlined in Theorem \ref{tMNMG} to any two basis which are connected by the automorphisms which are not simple.
\end{remark}
     
It follows from Remark \ref{rPREFIX} that the action of the automorphism $\vp$ well defines a homeomorphism $\tilde\vp$ on the space $\partial F_r$ of all infinite reduced words in the alphabet $\{ a_1^{\pm 1},\ld,a_r^{\pm 1}\}$. For an arbitrary automorphism this follows from the so called  Bounded  Cancellation Lemma in \cite{DC}. Since the labels of infinite rays in the Schreier subtree of the transitive action form a closed subset in the above space, we can apply Theorems \ref{lMEASURE} and  \ref{tMNMG}  to obtain the following result, which has no immediate connection to the action of $F_r$.

 \begin{corollary}\label{cNRA} For the homeomorphism $\tilde\vp$ of the space $\partial F_r$ induced by the above Nielsen automorphism $\vp$ of $F_r$ there is a \emph{closed} subset $W$ in $\partial F_r$ with positive measure such that the measure of $\tilde\vp(W)$ equals 0.
 \end{corollary}
 \begin{problem}\label{op4} Is there an analogue of Theorem \ref{tMNMG} in the case of modules over free associative algebras? 
 \end{problem}

\begin{center}
\textbf{Acknowledgment}
\end{center}
The authors are grateful to G. Arjantseva, G. Bergman, L. Bokut', R. Grigorchuk, I. Kapovich, D. Osin, V. Petrogradsky, A. Sambusetti and V. Shpilrain, with whom they discussed several aspects of the present paper. 

%%%%%%%%%%%%%%%%%%%%%%%%%%%%%%%%%%%%% 


\begin{thebibliography}{00}
\bibitem{A}Arzhantseva, G. N.,
\emph{On groups in which subgroups with a fixed number of generators are free},
Fundam. Prikl. Mat. \textbf{3} (1997), no. 3, 675--683.
\bibitem{BO} Bahturin, Y.; A. Olshanskii, \textit{Large restricted Lie algebras}, J. Algebra, \textbf{310} (2007), 413 - 427.
\bibitem{BO2} Bahturin, Y.; A. Olshanskii, \textit{ Schreier rewriting beyond the classical setting}, Science in China. Series A: Mathematics, to appear; arxiv:0811.1336v1 [math.RA].
\bibitem{PMC} Cohn, P.M.,  \textsc{Free Rings and Their Relations}, 2nd edition, Academic Press, London, 1985.
\bibitem{DC} Cooper, D.,  \emph{Automorphisms  of  free groups have finitely generated fixed point sets}, J.Algebra,  \textbf{111}(1987),  435-456.
\bibitem{G} Golod, E.S. On nil-algebras and residually finite $p$-groups (Russian), Izv. Akad. Nauk SSSR, Ser. Mat. \textbf{28}(1964), 273 - 276.
\bibitem{RG} Grigorchuk R., Symmetric random walks on discrete groups. In: Multicomplex random systems, Nauka, Moscow, 1978, p. 132 - 152.
\bibitem{KKS} Kaimanovich, V., Kapovich, I., Shupp, P. \textit{The subadditive ergodic theorem and generic stretching factors of group homomorphisms}, Israel. Math. J. 157(2007), 1 - 46.
\bibitem{KSS} Kapovich, I., Shpilrain V., Schupp, P.E.   \textit{Generic properties of Whitehead's Algorithm and isomorphism rigidity of random one-relator groups}, Pacific. J. Math., \textbf{223}(2006), 113-140. 
\bibitem{MA} Kilp, M.; Knauer, U.; Mikhalev, A., Monoids, Acts and Categories; Walter de Gruyter: Berlin, 2000.
\bibitem{KT} Kourovka Notebook: Unsolved Problems in Group Theory. AMS Translations Ser. 2, \textbf{121}(1983), 112pp.
\bibitem{JL} Lewin, Jacques., \textit{Free modules over free algebras and free group algebras: the Schreier technique}, Trans. AMS, \textbf{145}(1969), 455-465.
\bibitem{LS} Lyndon, R., and Schupp, P., Combinatorial Group Theory, Springer-Verlag, 2001.
\bibitem{BHN} Neumann, B.H., \textit{Groups covered by finitely many cosets}, Publ. Math. \textbf{3}(1954), 227-242.
\bibitem{OO} Olshanskii, A. Yu.; Osin, D.V.,  \textit{Large groups and their periodic quotients}, Proc. Amer. Math. Soc., 136 (2008), 753 - 759.
\bibitem{AS} Sambusetti, A., \emph{Growth tightness of free and amalgamated products},  Ann. Sci. \' Ecole Norm. Sup. 4 s\' erie, \textbf{35}(2002), 477-488. 
\bibitem{JS1} Stallings, J.R., \emph{Topology of finite graphs}, Invent. Math., \textbf{71} (1983), 551 - 565.
\bibitem{SZ} Stephenson, D.R.; Zhang, J.J., \emph{Growth of graded Noetherian rings}, Proc. Amer. Math. Soc. \textbf{125}(1997), 1593-1605. 
\end{thebibliography}
\end{document}